\documentclass{amsart}       
\usepackage[margin=1.2in]{geometry}
 \usepackage{quiver}
 \usepackage{tikz-cd}

\usepackage{hyperref}
\usepackage{graphicx}
\usepackage{fixltx2e}                    
\usepackage[utf8]{inputenc}            
\usepackage{nccmath, amssymb,latexsym, dsfont, amsfonts, amsbsy, amsmath, mathrsfs, dsfont, stmaryrd}            
\usepackage{mathtools}                   
\usepackage{caption}
\usepackage{bm}                          
\usepackage{enumerate}                   
\usepackage{verbatim}                    
\usepackage{url}                         
\usepackage{tabularx}

\usepackage[normalem]{ulem}

\usepackage{microtype}  
\usepackage[all, knot]{xy}
        \xyoption{arc} 
        \xyoption{web}                 
\makeatletter                            
\def\MT@register@subst@font{\MT@exp@one@n\MT@in@clist\font@name\MT@font@list
 \ifMT@inlist@\else\xdef\MT@font@list{\MT@font@list\font@name,}\fi}
\makeatother

\usepackage{color}

\newcommand{\muwave}[1]{\text{\uwave{$#1$}}}

\usepackage{scalerel,stackengine}
\stackMath
\newcommand\what[1]{%
\savestack{\tmpbox}{\stretchto{%
  \scaleto{%
    \scalerel*[\widthof{\ensuremath{#1}}]{\kern-.6pt\bigwedge\kern-.6pt}%
    {\rule[-\textheight/2]{1ex}{\textheight}}
  }{\textheight}%
}{0.5ex}}%
\stackon[1pt]{#1}{\tmpbox}%
}

\newcommand{\termDef}[1]{\textbf{\textit{#1}}}

\newcommand{\?}{\ensuremath{\mkern0.4\thinmuskip}}   
\newcommand{\with}{\?\?\?\&\?\?\?}

\let\leq=\leqslant
\let\geq=\geqslant
\let\Box=\square                            

\newcommand{\depth}{\mathbf{\mathtt{d}}}

\newcommand{\Diamondbox}{\mathrel{\text{\rlap{$\Diamond$}}\Box}}

\tikzcdset{scale cd/.style={every label/.append style={scale=#1},
    cells={nodes={scale=#1}}}}

\let\alg=\bm                                    
\let\class=\mathbb                             
\let\mod=\mathfrak							
	
\let\alg=\mathbf

\let\epsilon=\varepsilon
\let\Lambda\varLambda
\let\Gamma\varGamma
\let\Delta\varDelta
\let\Lambda\varLambda
\let\Omega\varOmega
\let\Theta\varTheta
\let\Xi\varXi
\let\Pi\varPi
\let\Sigma\varSigma

\theoremstyle{plain}
\newtheorem{theorem}{Theorem}[section]
\newtheorem{lemma}[theorem]{Lemma}
\newtheorem{corollary}[theorem]{Corollary}
\newtheorem{proposition}[theorem]{Proposition}

\newtheorem{obs}[theorem]{Observation}
\newtheorem*{claim*}{\bf Claim}
\theoremstyle{definition}

\newtheorem{definition}[theorem]{Definition}
\theoremstyle{remark}

\newtheorem{exa}{\bf Example}

\usepackage{multirow}

\begin{document}
\title[On the local modal product logic]{On the local modal product logic: standard completeness and decidability}

\author[A. Vidal]{Amanda Vidal}

\maketitle

\vspace{-0.8cm}
\begin{center}
Institute of Computer Science of the Czech Academy of Sciences\\
Pod Vodárenskou věží 271/2
182 00 Praha 8, Czech Republic \\
\textsf{amanda@cs.cas.cz}
\end{center}
\vspace{0.2cm}

\thispagestyle{empty}

\begin{abstract}

We study local consequence relations in modal extensions of product logic over Kripke models with either valued (fuzzy) or crisp accessibility relations. In both settings, we consider semantics over the full class of product algebras as well as over the standard product algebra on $[0,1]$.

Our main result is a constructive reduction of these modal logics to propositional product logic. As consequences, we prove that all the resulting systems are decidable and standard complete, i.e., the local consequence relation over all product algebras coincides with the one induced by the standard product algebra. In the valued-accessibility case, our methods strengthen previous results on decidability by extending them from theoremhood to arbitrary local consequence relations, and covering standard completeness. In the crisp case, the techniques are substantially different and yield, to the best of our knowledge, the first decidability and standard completeness results for local modal product logics with crisp accessibility relations.


%
\end{abstract}

\section{Introduction}

Modal logic is broadly understood as the study of different kinds of \textit{modalities}, or modes of truth, such as alethic (“necessarily”), epistemic (“it is known that”), deontic (“it ought to be”), temporal (“it is always the case that”), and provability (“it is provable that”), among others. Although these operators arise in different contexts, their common logical behavior justifies grouping them under the general label of “modal”. From a technical perspective, a logic may be regarded as modal whenever it includes operators that are not truth-functional, in the sense that their evaluation is not determined solely by the truth values of their arguments in the algebra associated with the underlying propositional logic. In practice, however, first-order logics are not usually classified as modal; rather, the term is typically reserved for systems lying between propositional and first-order logic in expressive power. This viewpoint naturally motivates the study of modal operators over non-classical propositional logics, which constitutes the general setting of the present work.

In this paper, we investigate modal extensions of product logic, one of the three principal systems commonly referred to as fuzzy logics \cite{Ha98}. Product logic is characterized by its completeness with respect to the standard product algebra, whose universe is the interval $[0,1]$ endowed with real multiplication as monoidal operation. This semantics makes product logic especially suitable for modeling certain forms of graded reasoning. Introducing modal operators into this framework raises a number of fundamental questions concerning the definition of the resulting systems, their axiomatizability, completeness, and computational properties.

The main contributions of this paper concern two of these questions: standard completeness and decidability. We study both modal systems based on Kripke frames with fuzzy (valued) accessibility relations and those based on classical (crisp) accessibility relations. For both semantics, we prove \textit{standard completeness}; namely, we show that the local consequence relation defined over the full variety of product algebras coincides with the one obtained by restricting attention to the standard product algebra. Moreover, we establish decidability of the corresponding local consequence relations. Both results are obtained through constructive reductions of the modal systems to propositional product logic (Theorems \ref{th:all} and \ref{th:allValued}).

The notion of modal many-valued logic considered here follows the tradition initiated by Fitting \cite{Fi92a, Fi92b} and Hájek \cite{HaHa96, Ha98}, and differs substantially from other approaches to modal substructural logics, such as those studied in \cite{O93, Re93b, Kam02}. In our setting, modal logics are interpreted over Kripke models in which both the accessibility relation and the valuation of propositional variables take values in suitable algebras at each world. The modal operators $\Box$ and $\Diamond$ generalize their classical counterparts and are, in general, not inter-definable \cite{RoVi20}.

These semantics naturally give rise to two notions of entailment: local and global derivability. In the local setting, entailment is evaluated world-wise, whereas in the global setting the premises are required to hold throughout the entire model. The relationship between these two notions and the first-order semantics of the underlying many-valued logic parallels the classical case, as well as the situation for Fisher--Servi intuitionistic modal logic \cite{FS77}.

The axiomatizability and computational properties of global modal logics over product and \L ukasiewicz logics were studied in detail in \cite{Vi22}. There it was shown that the systems defined over standard algebras with crisp accessibility relations are undecidable and not even recursively enumerable. These results align with the corresponding phenomena in first-order logic \cite{Scar62, Ra83, Ra83c}, and sharply contrast with the classical modal setting. In particular, they imply that the global logics determined by the standard algebras differ from those defined over the full variety of algebras, since the latter are recursively enumerable by virtue of the axiomatizability of the associated first-order logic. 

In contrast, comparatively little is known about local modal logics. For product logic, it is known that the transitive modal extension over the standard algebra with crisp accessibility is undecidable \cite{Vi20}. On the other hand, the set of theorems of the modal logic based on valued accessibility over the standard algebra is decidable \cite{CeEs22}. This latter result provides one of the motivations for the present work. In Section \ref{sec:valued}, we extend the techniques developed in \cite{CeEs22} to prove completeness with respect to the standard product algebra for the full class of product algebras (Theorem \ref{th:standardCompVal}), together with decidability of the corresponding local consequence relation (Theorem \ref{th:decVal}). As we later explain, however, these techniques do not extend to the crisp case, which is especially significant because crisp accessibility corresponds precisely to the classical notion of Kripke frame.

Looking at local modal logics based on the other principal fuzzy logics offers only limited guidance, since those systems are considerably simpler in several respects. For example, first-order \L ukasiewicz logic over arbitrary chains is complete with respect to witnessed models \cite[Theorem~8]{CiH06}. As a consequence, local modal \L ukasiewicz logics are complete with respect to finite models, and in particular with respect to models whose size is bounded by the modal complexity of the formulas involved. Decidability and standard completeness therefore follow rather directly by reduction to propositional \L ukasiewicz logic. Nevertheless, no finite axiomatization of these modal systems is currently known.

By contrast, local modal logics based on G\"odel logic share some of the difficulties present in the product case, most notably the failure of the finite model property. Still, standard completeness follows relatively easily from the analogous first-order result \cite[Theorem~5.3.3]{Ha98}, unlike in the product setting. Furthermore, both local G\"odel modal logics—with crisp and with valued accessibility relations—are known to be decidable \cite{CaMeRo17}. However, the methods used there rely heavily on algebraic properties specific to G\"odel algebras, and attempts to adapt them to product logic have so far been unsuccessful.

The paper is organized as follows. Section \ref{sec:Prelim} introduces the necessary background, covering topics from propositional logic to modal and first-order logics. Section \ref{sec:statements} presents in detail the main results on modal product logic with crisp accessibility, namely Theorem \ref{th:all}, together with its two principal consequences: the decidability and standard completeness of $K\Pi$. Sections \ref{sec:soundness} and \ref{sec:completeness} are devoted to the proofs of the main result stated in the previous section. Finally, Section \ref{sec:valued} establishes Theorem \ref{th:allValued}, an analogue of the main theorem from Section \ref{sec:statements} for modal product logic with valued accessibility.

For the reader's convenience, some proofs have been deferred to appendices corresponding to each section, in order to keep the presentation of the main results and their proofs as readable as possible, without unnecessary technical interruptions. Proofs are included there either because they consist of relatively direct calculations or because, while technically involved, their detailed presentation is less essential for the overall understanding of the paper.

\section{Preliminaries}\label{sec:Prelim}

Before proceeding with more specific notions, we remark that in this paper a logic is identified with a consequence relation \cite{Fo16}, as opposed to a set of formulas. 
While the second approach is more common in the literature of modal logics \cite{ChZa97}, we opt for the former definition because the differences between local and global modal logics are lost if only the tautologies of the logic are considered. In particular, in this paper we will work with finitary logics, namely, those that are determined by the entailments from finite sets of premises. Consequently, when we say that a logic $\models$ is decidable, we mean that the set of entailments following from finite sets of premises is decidable; namely, that the set
$\{\langle \Gamma, \varphi \rangle\colon \mid \Gamma\mid  < \omega, \Gamma \models \varphi\}$ is decidable.

Let us start by introducing the underlying propositional product logic and some of its characteristics.

The variety of \termDef{product algebras} \cite{HaGoEs96,CiT00} $\class{P}$ is formed by algebras with structure $A = \langle A,  \with, \rightarrow,  \bot \rangle$ that of type $\langle 2,2, 0\rangle$ such that, upon defining the derived operations
\begin{align*}
\top \coloneqq& \bot \rightarrow \bot;\  &  \neg x \coloneqq& x \rightarrow \bot; \\
x \wedge y \coloneqq& x \with (x \rightarrow y);\ &   x \vee y \coloneqq& ((x \rightarrow y) \rightarrow y) \wedge ((y \rightarrow x) \rightarrow x);
\end{align*}
 satisfies the following axioms:
 \begin{enumerate}
 \item[(A1)] $\langle A, \with, \top\rangle$ is an abelian monoid (that is, $\with$ is associative, commutative and $x \with \top = x$ for all $x$);
 
  \item[(A2)] $\langle A, \wedge, \vee, \bot, \top \rangle$ is a lattice with smallest element $\bot$ and greatest element $\top$;
  \item[(A3)] $x \with (y \vee z) = (x \with y) \vee (x \vee z), x \with (y \wedge z) = (x \with y) \wedge (x \with z)$;
  \item[(A4)] $(x \with y) \rightarrow z = x \rightarrow (y \rightarrow z)$;
   \item[(A5)] $(x \rightarrow y) \vee (y \rightarrow x) = \top$;
   \item[(A6)] $x \rightarrow x = \top$;
   \item[(A7)] $(\neg \neg z \with ((x \with z) \rightarrow (y \with z))) \rightarrow (x \rightarrow y) = \top$;
   \item[(A8)] $x \wedge \neg x = \bot$.
  \end{enumerate}
  As usual, we will also define $x \leftrightarrow y \coloneqq (x \rightarrow y) \wedge (y \rightarrow x)$ and $\varphi^n \coloneqq \varphi \with \ldots \with \varphi$ $n$ times.
%
%

A formula $\varphi$ follows from a set of premises $\Gamma$ (written in the language of product algebras, over a countable set of variables $\mathcal{V}$)
in a product algebra $\alg{A}$ (in symbols, $\Gamma \models_{\alg{A}} \varphi$) whenever there is a finite $\Gamma_0 \subseteq_\omega\Gamma$ such that, for any homomorphism $h$ from the set of formulas in the language of product algebras into $\alg{A}$, if $h(\gamma) = \top$ for all $\gamma \in \Gamma_0$, then $h(\varphi) = \top$ too. We say that formula $\varphi$ follows from a set of premises $\Gamma$ in \termDef{product logic}, and we write $\Gamma \models_\Pi \varphi$, whenever $\Gamma \models_{\alg{A}} \varphi$ for any $\alg{A} \in \class{P}$.

It immediately follows from the subdirect representation Theorem (see e.g., \cite{BuSa00}) that product logic is complete with respect to the finitely subdirectly irreducible elements in the variety of product algebras, which are the linearly ordered ones (to which we will also refer-to by the term \textit{chains}). Not only this, but it is known that product logic is also complete with respect to the so-called \textit{standard} product algebra $[0,1]_\Pi$, the one whose universe is the real unit interval $[0,1]$ and whose conjunction is the usual real product operation \footnote{i.e., $\top = 1$, $\bot = 0$ and for $a,b \in [0,1]$,  $a \with b = a \cdot b$ and $a \rightarrow b = \begin{cases} 1 &\hbox{ if }a \leq b, \\ \frac{b}{a} &\hbox{ otherwise.} \end{cases}$}  \cite{Ha98}[Theorem~4.1.13]. Namely, $\Gamma \models_\Pi \varphi$ if and only if $\Gamma \models_{[0,1]_\Pi} \varphi$, for arbitrary set of formulas $\Gamma \cup \{\varphi\}$\footnote{This is not the case if the entailment is defined as $\Gamma \models_{\alg{A}} \varphi$ iff $h(\Gamma) \subseteq \{1\}$ implies $h(\varphi) = 1$ for any homomorphism, $h$, since for infinite sets of premises the logic over $[0,1]$ and that over the generated variety differ.}.

The Baaz-Monteiro $\Delta$ operator allows for a convenient expansion of the previous language with a unary symbol, which will be of use in this work. For the definition of the operator over arbitrary product algebras see for instance \cite[Definition~2.4.6]{Ha98}. The most relevant feature of $\Delta$ that we will be using in this work is that, for any linearly ordered product algebra $\alg{A}$ extended with $\Delta$, and for any $a \in A$, the following holds:
\[\Delta(a) = \begin{cases} \top &\hbox{ if }a = \top,\\ \bot &\hbox{ otherwise.}\end{cases}\]

We will denote by $\models_{\Pi_\Delta}$ the logic arising from product algebras expanded by the $\Delta$ operation, defined in the same fashion as $\models_\Pi$. It is the case that also $\models_{\Pi_\Delta}$ is complete with respect to chains in $\class{P}$ expanded with $\Delta$ (as above), and also, with respect to the standard product algebra expanded with $\Delta$ \cite[Theorem 4.1.13]{Ha98}. Furthermore, while $\models_{\Pi}$ enjoys only a weak version of the Deduction Theorem, called \textit{local} Deduction theorem, $\models_{\Pi_\Delta}$ enjoys a strong version of this property. Formally:
\begin{itemize}
\item $\Gamma, \varphi \models_{\Pi} \psi$ if and only if there is some $n\in \omega$ for which $\Gamma \models_{\Pi} \varphi^n \rightarrow \psi$ (local D.T).
\item $\Gamma, \varphi \models_{\Pi_\Delta} \psi$ if and only if $\Gamma \models_{\Pi_\Delta} \Delta \varphi \rightarrow \psi$ (D.T).
\end{itemize}

It is known that the set of theorems of $\models_\Pi$ and of $\models_{\Pi_\Delta}$ are decidable \cite[Lemma 6.2.16]{Ha98}, \cite[Theorem 6.2.1]{Ha11}. We will use in this work that the logics $\models_{\Pi}$ and $\models_{\Pi_\Delta}$ are decidable (as consequence relations), which is a direct result of the above D.T. for $\Pi_\Delta$.\footnote{ It is immediate for $\models_{\Pi_\Delta}$, and follows for $\models_{\Pi}$ because $\Gamma, \psi \models_{\Pi} \varphi$ if and only if $\Gamma \models_{\Pi_\Delta} \Delta \psi \rightarrow \varphi$, for $\Gamma, \varphi, \psi$ formulas without $\Delta$.}

In the following sections, we will denote by $Fm_\Pi$ and $Fm_{\Pi\Delta}$ to, respectively, the set of formulas in the language of product algebras (over a countable set of variables $\mathcal{V}$), and that of formulas in the language of product algebras with $\Delta$.

\subsection{Modal logics}

In an analogous way on how modal logic is built from classical logic, we will consider a language expansion of the one of product logic by means of two unary operators, $\Box$ and $\Diamond$, and denote this set of formulas by $Fm$. We will sometimes use the symbol $\heartsuit$ to denote instinctively any of the previous two modal operators.
Let us begin by introducing several purely syntactical notions.

 We define the set of (modal) \termDef{subformulas} of a given formula $SFm()$ inductively as usual. Namely,
\begin{align*}
SFm(p) \coloneqq& \{p\} & \text{ for }p \text{ propositional variable or constant,}\\
SFm(\varphi \star \psi) \coloneqq& \{\varphi \star \psi\} \cup SFm(\varphi) \cup SFm(\psi) & \text{ for }\star \text{ propositional connective,}\\
SFm(\heartsuit\varphi) \coloneqq& \{\heartsuit \varphi\} \cup SFm(\varphi). & 
\end{align*} 
We will say a set of (modal) formulas $\Sigma$ is \termDef{closed} whenever it is closed under $SFm()$.

On the other hand, we define the set of \termDef{propositional subformulas} $PSFm()$ \footnote{The name "propositional subformulas" refers to the fact that we only go down up to propositional operations. Observe that this set is not formed by propositional formulas.} by letting
\begin{align*}
PSFm(p) \coloneqq& \{p\} & \text{ for }p \text{ propositional variable or constant,}\\
PSFm(\varphi \star \psi) \coloneqq& \{\varphi \star \psi\} \cup PSFm(\varphi) \cup PSFm(\psi) & \text{ for }\star \text{ propositional connective,}\\
PSFm(\heartsuit\varphi) \coloneqq& \{\heartsuit \varphi\}.  & 
\end{align*} 
 For $\Psi$ a set of modal formulas, we let $(P)SFm(\Psi) \coloneqq \bigcup_{\psi \in \Psi} (P)SFm(\psi)$. Clearly if $\Psi$ is a set of propositional formulas, $SFm(\Psi) = PSFm(\Psi)$, and it is closed (as defined before) whenever it is closed in the propositional sense (i.e., all subformulas belong to it).

%

\begin{definition}
Let $\psi$ be a modal formula. Its \termDef{modal depth} $md(\psi)$ is recursively defined by
\begin{align*}
md(p) &\coloneqq 0, & \text{ for }p \text{ variable or constant,} \\
md(\varphi \star \phi), &\coloneqq  max\{md(\varphi), md(\phi)\}& \text{ for }\star \text{ propositional connective,}\\
md(\heartsuit \psi) &\coloneqq  1 + md(\psi).& 
\end{align*}
For $\Psi$ a finite set of modal formulas, we let $md(\Psi) \coloneqq max\{md(\psi)\colon \psi \in \Psi\}.$
\end{definition} 

Let us introduce a last syntactical definition, which allows us to split in levels the closure of a set of modal formulas. 
 \begin{definition}\label{def:modalLevels}
 Let $\Psi$ be a finite set of modal formulas.  For $0 \leq i < md(\Psi)$ we let
 \begin{align*}
 \Psi_0 \coloneqq PSFm(\Psi), & & 
 \Psi_{i+1} \coloneqq \bigcup_{\heartsuit \psi \in \Psi_i} PSFm(\psi).
 \end{align*} 
 
 \end{definition}

 The modal logic we will work with in this paper will be the one arising, in a semantical spirit, from a suitable generalization of (classical) Kripke models -the usual semantics of modal logic- evaluated over product algebras. 

\begin{definition}Let $\alg{A}$ be a product algebra with universe $A$.
\begin{itemize}
\item An \termDef{$A$-valued Kripke frame} is a pair $\langle W, R \rangle$ where $W$ is a non-empty set (of worlds) and $R \colon W \times W \rightarrow A$.

 We will say that the frame is \termDef{crisp} -or classical- whenever $R$ is a binary accessibility relation in $W$ (i.e., $R \subseteq W \times W$), or, equivalently, if $R \colon W \times W \rightarrow \{\bot, \top\}$.\footnote{The convention is then that $R(v,w) = \top$ if and only if $\langle v, w \rangle \in R$.} In this latter case, we will write $Rvw$ instead of $\langle v,w \rangle \in R$, and, for a fixed $v$, whenever the model is clear from the context, we will sometimes write $Rvw$ as a shortening for "$w \in W$ such that $Rvw$".

\item An \termDef{$\alg{A}$-Kripke model} $\mod{M}$ is a pair $\langle \mod{F}, e\rangle$ where $\mod{F}$ is an $A$-valued Kripke frame, and $e \colon W \times \mathcal{V} \rightarrow A$.  The evaluation $e$ is extended to $W \times Fm$ as a world-wise propositional homomorphism\footnote{i.e., $e(v, \bot) = \bot, e(v, \top) = \top$, and $e(v, \varphi \star \psi) \coloneqq e(v,\varphi) \star e(v, \psi)$ for any propositional connective $\star$.}, and
\[e(v, \Box \varphi) \coloneqq \bigwedge \{R(v,w) \rightarrow e(w, \varphi) \colon w \in W\}, \qquad e(v, \Diamond \varphi) \coloneqq \bigvee \{R(v,w) \with e(w, \varphi)\colon w \in W\}.\]
If any of the above infima/suprema do not exist, the value is undefined, and this condition is inherited. The model is said to be \termDef{safe} whenever, for any world $v$ in its universe, and any formula $\varphi$, the value $e(v, \varphi)$ is defined.

We will say that the model is \termDef{crisp} whenever its underlying frame is crisp. In this case, it is clear that the evaluation of the modal formulas can be rewritten as follows:
\[e(v, \Box \varphi) = \bigwedge \{e(w, \varphi) \colon Rvw\}, \qquad e(v, \Diamond \varphi) = \bigvee \{e(w, \varphi) \colon Rvw\}.\]
\end{itemize}

\end{definition}
We will usually assume the notation $\mod{M} = \langle W, R, e\rangle$ when talking about (valued) Kripke models. If we are using several models at the same time, we will use symbols $\mod{M}, \mod{N}$, or subindexes/superindexes, and refer to its elements by $W[\mod{M}], R[\mod{M}]$ and $e[\mod{M}]$.

Given a class of product algebras $\class{A}$ we denote by $\class{M}_{\class{A}}$ to the class of safe $\alg{A}$-models, for each $\alg{A} \in \class{A}$. If $\class{A}$ is a singleton $\{\alg{A}\}$, we will write $\class{M}_{\alg{A}}$.

\begin{definition}\label{def:tree}
A \termDef{tree} is a frame $\mod{T} = \langle W, R\rangle$ (or a model with such a frame as underlying structure, it will be clear from the context) for which there is a splitting $W = \bigcup_{n \in \omega} W_n$, with $W_n \bigcap W_m = \emptyset$ for any $n \neq m$, and such that:
\begin{itemize}
\item $W_0$ is a singleton (the \textit{root} of the tree),
\item If $R(v,w) > 0$ (or $Rvw$ if we are considering a crisp frame/model), then there is some $n \in \omega$ for which $v \in W_n$ and $w \in W_{n+1}$,
\item If $w \in W_{n+1}$ then there is a single $v \in W_n$ such that $R(v,w) > 0$ (resp. $Rvw$),
\item If $W_n = \emptyset$ then $W_m = \emptyset$ for all $m \geq n$. 
\end{itemize}
\end{definition}
%
%

Given a world $v$ in a tree $\mod{T}$ with root $r$, we refer to the \termDef{depth} of $v$ in $\mod{T}$, and denote by $\depth(v)_{\mod{T}}$, to the minimum distance from the root to that world, or $\omega$ if this is infinite. Namely, 
\[\depth(v)_\mod{T} \coloneqq \begin{cases} min\{n \in \omega\colon \exists w_0 = r, \ldots, w_n = v \text{ s.t. } R(w_i, w_{i+1}) > 0 \forall i < n\} &\hbox{ if this value exists, } \\ \omega &\hbox{otherwise.} \end{cases}  \]
We will remove the reference to $\mod{T}$ in the previous definition, and simply write $\depth(v)$, if this is clear from the context. 
For $v \in W$ in a tree and different from the root, we denote by $par(v)$ its parent node, i.e., the only $u \in W$ such that $Ruv$.

The \termDef{subframe  of $\mod{F}$ (or submodel of $\mod{M}$) generated by one world $w$}, denoted by $\mod{F}_w$ (resp. $\mod{M}_w$) is the one obtained by restricting the original frame or model to the worlds accessible from $w$. We will restrict these to the worlds accessible in a finite number of steps.\footnote{We will see below that our logics are complete with respect to models with finite depth, so this definition will be enough for us.} Formally:
\[W_w \coloneqq \{w\} \cup \bigcup_{n \in \omega} \{u \in W\colon \exists u_1, \ldots u_n = u \text{ s.t.} R(u_i, u_{i+1}) > 0 \text{ for each } 1 \leq i < n\}\]
and $\mod{F}_w$ (resp. $\mod{M}_w$) is the corresponding frame (model) restricted to the previous universe.

Given a class of product algebras, there are two modal logics arising from its class of models in a natural way: the local and the global ones. Along this work we will focus in the local logics, defined in the following way.
\begin{definition}
Given a \textit{finite} set of formulas $\Gamma \cup \{\varphi\} \subset Fm$ and a model $\mod{M}$, we say that $\varphi$ follows (locally) from $\Gamma$ in $\mod{M}$ (in symbols $\Gamma \models_{\mod{M}} \varphi$) whenever for each $v \in W$, if $e(v, \Gamma) \subseteq \{\top\}$ then $e(v, \varphi) = \top$ too. 

More in general, given an arbitrary set of formulas $\Gamma \cup \{\varphi\} \subseteq Fm$ and a class of algebras $\class{A}$, we say that $\varphi$ follows from $\Gamma$ in the class of safe models arising from $\class{A}$ (in symbols, $\Gamma \Vdash_{\class{M}_{\class{A}}} \varphi$) whenever there is a finite set $\Gamma_0 \subseteq_{\omega} \Gamma$ such that, for each $\mod{M}\in \class{M}_{\class{A}}$, $\Gamma_0 \models_{\mod{M}} \varphi$. 
\end{definition}
Clearly, for $\Gamma$ finite, if $\Gamma \not \Vdash_{\mod{M}} \varphi$, there is a world $v$ where the above condition is falsified (i.e., $e(v, \Gamma) \subseteq \{\top\}$ and $e(v, \varphi) < \top$). In this case, sometimes we will write $\Gamma \not \Vdash_{\mod{M},v} \varphi$.

In order to lighten the notation, we will denote by $\mathrm{M}\Pi$ the class $\class{M}_{\class{P}}$, and we will refer to them as \termDef{valued (product Kripke) models}, and by $\mathrm{M}[0,1]_{\Pi}$ the class $\class{M}_{[0,1]_\Pi}$, which we will refer to as \termDef{valued standard (product Kripke) models}. Similarly, we will denote by $\mathrm{K}\Pi$ the class of crisp models in 
$\mathrm{M}\Pi$, and refer to them as simply \termDef{(product Kripke) models}, and by 
$\mathrm{K}[0,1]_{\Pi}$ the class of crisp models in 
$\mathrm{M}[0,1]_{\Pi}$, and refer to them as \termDef{standard (product Kripke) models}.\footnote{We use the terminology ``valued product Kripke model" and ``product Kripke model", instead of talking about crisp product Kripke models, because the paper mainly works with this second kind of models and the terminology is therefore lighter.}
Accordingly, in this paper we will work with the logics 
$\Vdash_{\mathrm{M}\Pi}, \Vdash_{\mathrm{M}[0,1]_\Pi}$ and $\Vdash_{\mathrm{K}\Pi}, \Vdash_{\mathrm{K}[0,1]_\Pi}$.
Below, we will talk about valued models, and of course these also comprise models (i.e., crisp models, as introduced before).

The last concept we will introduce over modal product logics is that of the witness of a formula beginning with a modality in a world of a model. In modal (classical) logic this notion is not much used because it does not provide, in general, any additional information, but we will see it does in modal product logic.
 
\begin{definition}\label{def:witness}
Let $\mod{M}$ be a valued Kripke model, $v\in W$ and $\varphi \in Fm$. We let the \termDef{witnesses} of $\Box \varphi$ and $\Diamond \varphi$ in $v$ (in the model $\mod{M}$) be defined as follows:
\begin{itemize}
\item $wits_{\mod{M}}(v, \Box \varphi) \coloneqq \{w \in W\colon e(v,\Box \varphi) = R(v,w) \rightarrow e(w, \varphi)\},$
\item $wits_{\mod{M}}(v, \Diamond \varphi) \coloneqq \{w \in W\colon e(v,\Diamond \varphi) = R(v,w) \with e(w, \varphi)\}.$
\end{itemize}
\end{definition}
We will omit the subscript or references to the model $\mod{M}$ if it is clear from the context.
Observe that, if the model is crisp, the above definitions can be rephrased, for $\heartsuit \in \{\Box, \Diamond\}$ as
\[wits_{\mod{M}}(v, \heartsuit \varphi) = \{w \in W\colon Rvw \text{ and } e(v, \heartsuit \varphi) = e(w, \varphi)\}.\]

For a valued model $\mod{M}$ and a world $v \in W$ and a set of formulas $\Sigma$, we let the set of \termDef{unwitnessed formulas of $\Sigma$ at $v$ in $\mod{M}$} be the set 
\[uWit_{\mod{M}}(v, \Sigma) \coloneqq \{\heartsuit \psi \in \Sigma\colon wit_\mod{M}(v, \heartsuit \psi) = \emptyset\}.\]

A valued model $\mod{M}$ is \termDef{witnessed} whenever for every formula  $\heartsuit \psi \in Fm$, and every world $v \in W$, $wits_{\mod{M}}(v, \heartsuit \varphi) \neq \emptyset$. 
A valued model $\mod{M}$ is \termDef{quasi-witnessed} whenever for any $v\in W$ and $\varphi \in Fm$, $wits_{\mod{M}}(v, \Diamond \varphi) \neq \emptyset$ and either $wits_{\mod{M}}(v, \Box \varphi) \neq \emptyset$ or $e(v, \Box \varphi) = \bot$.

It is easy to see that, as it happens in classical modal logic (e.g., see \cite[Prop.~2.15]{BlRiVe01}), the above modal logics are complete with respect to valued trees with finite depth.  Let us properly formulate the direction we are interested in, which we will use along the paper without further notice. The proof, a natural generalization of the classical case, is included in the Appendix corresponding to this section (Appendix 1).
\begin{proposition}\label{prop:compTreesFiniteDepth}
Let $\Gamma \cup \{\varphi\} \subseteq_{\omega} Fm$, $\mathbf{A}$ a product algebra, $\mod{M} \in \class{M}_{\mathbf{A}}$ and $v \in W$ be such that $\Gamma \not \Vdash_{\mod{M}, v} \varphi$. Then there is a tree $\mod{T} \in  \class{M}_{\mathbf{A}}$ with root $v'$ for which
$\Gamma \not \Vdash_{\mod{T}, v'} \varphi$. Moreover, if $\mod{M}$ was witnessed, quasi-witnessed or crisp, the same condition holds for the resulting tree.
 \end{proposition}

Again, as in the classical case, in order to preserve the values of $\Upsilon = \Gamma \cup \{\varphi\}$ in the root of $\mod{T}$, it is only necessary to preserve the values of the formulas in $\Upsilon_i$ at each world at depth $i$. Since by definition, $\Upsilon_i$ will be empty for any $i > md(\Upsilon)$, we can crop the previous tree at depth $md(\Upsilon)$, obtaining the following corollary.

\begin{corollary}\label{cor:LocalFiniteDepth}
Let $\Gamma \cup \{\varphi\} \subseteq_{\omega} Fm$, $\mathbf{A}$ a product algebra, $\mod{M} \in \class{M}_{\mathbf{A}}$ and $v \in W$ be such that $\Gamma \not \Vdash_{\mod{M}, v} \varphi$. Then there is a tree of finite depth $\mod{T} \in  \class{M}_{\mathbf{A}}$ with root $v'$ for which
$\Gamma \not \Vdash_{\mod{T}, v'} \varphi$. Moreover, if $\mod{M}$ was witnessed, quasi-witnessed or crisp, the same condition holds for the resulting tree.
\end{corollary}

Let us conclude this subsection with a notion that will be useful later on, consisting on a very simple syntactic translation from the modal language to the propositional one over an extended set of (propositional) variables. 
\begin{definition}\label{def:modalToPropFormulas}
Given a (modal) formula $\varphi$ and an arbitrary symbol $\alpha$, let $Vars(\varphi)$ be the set of variables appearing in $\varphi$. we denote by $\varphi_\alpha$ the propositional formula in the language $\{x_\alpha\colon x \in Vars(\varphi)\} \cup \{(\heartsuit \chi)_\alpha \colon \heartsuit \chi \in PSFm(\varphi)\}$ inductively defined as follows:
\begin{align*}
x_\alpha \coloneqq x_\alpha \text{ for }x \in Vars(\varphi); \quad 
(\psi \star \chi)_\alpha \coloneqq \psi_\alpha \star \chi_\alpha \text{ for }\star \text{ propositional connective}; \quad 
(\heartsuit \chi)_\alpha \coloneqq (\heartsuit \chi)_\alpha.
\end{align*}
\end{definition}

For example, given the formula $\varphi \equiv y \with (\Box x \rightarrow (x \rightarrow \Diamond (y \with x)))$ and the symbol $\langle 0 \rangle$, we have that $\varphi_{\langle 0 \rangle} \equiv y_{\langle 0 \rangle} \with ((\Box x)_{\langle 0 \rangle} \rightarrow (x_{\langle 0 \rangle} \rightarrow (\Diamond (y \with x))_{\langle 0 \rangle}))$.

\subsection{Relation of modal and predicate logics}\label{subsec:prel-predLogic}

The modal product logics defined above can be seen as a certain fragment of the corresponding predicate product logics (namely, over $[0,1]_{\Pi}$ or over the whole variety of product algebras). Since we will be interested only in how certain results from the predicate product logic can be used in modal product logics, we will only do a brief introduction of the relevant definitions and results. For more details on first order product logic, we refer the reader to Section 5.1 of \cite{Ha98}. 

A predicate language consists of a non-empty set of predicates, each together with a positive natural number -its arity- and a (possibly empty) set of object constants. Predicates will be  mostly denoted by $P,Q,R\ldots $ and constants by $c,d,\ldots$. Logical symbols are object variables $x,y,z,\ldots$, connectives $\with, \rightarrow \wedge, \vee$, truth constants $\top, \bot$ and quantifiers $\forall, \exists$. Terms are object variables and object constants. 

Atomic formulas are either truth constants or have the form $P(t_1, \ldots, t_n)$ where $P $ is a predicate of arity $n$ and $t_1,\ldots, t_n$ are terms. If $\varphi, \psi$ are formulas and $x$ is an object variable then $\varphi \star \psi$, $\forall x \varphi$ and $\exists x \varphi$ are formulas, for $\star$ propositional connective; each formula results from atomic formulas by iterated use of this rule.

Let $\mathcal{P}$ be a predicate language and $\alg{A}$ a product algebra. An \termDef{$\alg{A}$-structure} or model $\mod{M}=\langle M, \mathcal{P}^{\mod{M}}\rangle$ for $\mathcal{P}$ has a non-empty domain $M$, for each predicate $P \in \mathcal{P}$ of arity $n$, a mapping $P^{\mod{M}}\colon M^n \rightarrow A$, and for each object constant $c \in \mathcal{P}$, $c^{\mod{M}} \in M$. 
An $\mod{M}$-evaluation of object variables is a mapping $e$ assigning to each object variable $x$ an element $e(x) \in M$.

For $e$ an evaluation, $x$ an object variable and $m \in M$, we denote by $e[x\mapsto m]$ the evaluation $e'$ such that $e'(y) = e(y)$ for each object variable $y$ different from $x$, and $e'(x) = m$. 
For a structure $\mod{M}$, we let $\vert x \vert_{\mod{M},e} \coloneqq e(x)$, and $\vert c \vert_{\mod{M},e} \coloneqq c^{\mod{M}}$. 
We then define the truth value $\vert \varphi \vert_{\mod{M},e}$ of a formula $\varphi$ in $\mod{M}$ under $e$ as follows:
\begin{align*}
\vert P(t_1, \ldots, t_n)\vert_{\mod{M},e} \coloneqq& P^{\mod{M}}(\vert t_1\vert_{\mod{M},e}, \ldots, \vert t_n\vert_{\mod{M},e}); & 
\vert (\forall x)\varphi\vert_{\mod{M},e} \coloneqq& \bigwedge_{m \in M}  \vert \varphi \vert_{\mod{M},e[x\mapsto m]};\\
\vert \varphi \star \psi \vert_{\mod{M},e} \coloneqq& \vert \varphi \vert_{\mod{M},e} \star \vert \psi \vert_{\mod{M},e} \text{ for }\star \text{ propositional connective}; & 
\vert (\exists x)\varphi\vert_{\mod{M},e} \coloneqq& \bigvee_{m \in M}  \vert \varphi \vert_{\mod{M},e[x\mapsto m]}.
\end{align*}

As for Kripke models, the above infima and suprema might not exist, leaving the value of that formula undefined.
The structure $\mod{M}$ is \termDef{safe} for the evaluation $e$ if the truth values of every formula are defined. It is witnessed (for $e$) whenever for each formula $\varphi(x)$ in one free variable $x$ there exist 
$m,n \in M$ such that $\vert (\forall x)\varphi\vert_{\mod{M},e} = \vert \varphi \vert_{\mod{M},e[x\mapsto m]}$ and $\vert (\exists x)\varphi\vert_{\mod{M},e} = \vert \varphi \vert_{\mod{M},e[x\mapsto n]}$. It is quasi-witnessed (for $e$) whenever, for each formula $(\exists x) \varphi$  there is $n \in M$ such that $\vert (\exists x)\varphi\vert_{\mod{M},e} = \vert \varphi \vert_{\mod{M},e[x\mapsto n]}$ (i.e., witnessed for all existential formulas) and, for each formula $(\forall x) \varphi$, either there is $m$ such that $\vert (\forall x)\varphi\vert_{\mod{M},e} = \vert \varphi \vert_{\mod{M}ev[x\mapsto m]}$ (i.e., witnessed) or $\vert (\forall x)\varphi\vert_{\mod{M},e} = \bot$.

Observe that all the previous notions, from the definition of the evaluation of the quantifiers, to the notions of safe, witnessed and quasi-witnessed models are analogous to the ones given for valued product Kripke models.
Furthermore, observe that if $\varphi$ is a sentence (i.e., it has no free variables), its truth value in a model does not depend on the evaluation. Similarly, if $\varphi$ has only one free variable $x$, its truth value depends only on the evaluation of that free variable. Accordingly, in this latter case we will write 
$\vert \varphi(x) \vert_{\mod{M},[x \mapsto m]}$ to denote the value of that formula in $\mod{M}$
 when the only free variable $x$ is mapped to $m$.

	We can use the usual translation from modal to predicate logics in order to get a reduction from the modal logics over a class of algebras to certain deductions over the corresponding predicate logic. 
	For $\varphi$ modal formula, consider the predicate language $\{R/2\}, \{E_p/1\colon p\text{ variable in }\varphi\}$. For arbitrary $i \in \omega$, let us define the translation $\langle \varphi, x_i \rangle^{\sharp}$ recursively by letting
	\begin{itemize}
		\item  $\langle p, x_i \rangle^\sharp \coloneqq E_p(x_i)$;
		\item $\langle \varphi \star \psi, x_i \rangle^\sharp \coloneqq \langle \varphi, x_i \rangle^\sharp \star \langle \psi, x_i \rangle^\sharp$ for $\star$ propositional connective;
		\item $\langle \Box \varphi, x_i \rangle^\sharp \coloneqq \forall x_{i+1} R(x_i, x_{i+1}) \rightarrow \langle \varphi, x_{i+1} \rangle^\sharp$;
		\item $\langle \Diamond \varphi, x_i \rangle^\sharp \coloneqq \exists x_{i+1} R(x_i, x_{i+1}) \with \langle \varphi, x_{i+1} \rangle^\sharp$.
	\end{itemize}
	
	For $\mod{M} = \langle W, R, e\rangle$ a valued Kripke model over an algebra $\alg{A}$ we let $\mod{M}^\times$ be the first order model (over $\alg{A}$) in the language  $\{R/2\} \cup \{E_p/1\colon p\text{ variable in }\mathcal{V}\}$ given by\footnote{We write $M^\times, R^\times$ and $E_p^\times$ to denote respectively $M^{\mod{M}^\times}, R^{\mod{M}^\times}$ and $E_p^{\mod{M}^\times}$.}:
	\[ M^\times \coloneqq W, \qquad  R^\times(u,v) \coloneqq R(u,v), \qquad  E^\times_p(u) \coloneqq e(u,p) \text{ for each }p \in \mathcal{V}.\]

If $\mod{M}$ is a (crisp) Kripke model, we will let $R^\times(u,v) \coloneqq \top$ if $Ruv$ and $\bot$ otherwise.

	Similarly, for $\mod{N} = \langle M^{\mod{N}}, R^{\mod{N}}, \{E_p^{\mod{N}}\colon p \text{ variable in }\mathcal{V}\} \rangle$ a (first order) structure over the algebra $\alg{A}$ in the language $\{R/2\} \cup \{E_p/1\colon p\text{ variable in }\mathcal{V}\}$ we let $\mod{N}_\times$ be the valued Kripke model over $\alg{A}$ given by
	\[W_\times \coloneqq M^{\mod{N}}, \qquad R_\times(u,v) \coloneqq R^{\mod{N}}(u,v), \qquad e_\times(u,p) \coloneqq E^{\mod{N}}_p(u) \]

	If $\mod{N}$ validates the sentence $\forall x,y R(x,y) \vee \neg R(x,y)$ we will call it \termDef{R-crisp}, and we let $\mod{N}_\times$ be the crisp Kripke model where the universe and evaluation are as above, and $R_\times u,v$ if and only if $ R^{\mod{N}}(u,v) = \top$.

The next observation can be checked easily by induction on the formulas.
	 
	\begin{obs}\label{obs:languageModalToFO} Let $\varphi$ be a modal formula over variables $\mathcal{V}$. Then,
	\begin{enumerate}
	\item for $\mod{M}$ a 
	Kripke model over $\alg{A}$, and $v$ in $W$ it holds that \[\vert \langle \varphi, x_0\rangle^\sharp\vert_{\mod{M}^\times, [x_0 \mapsto v]} = e(v,\varphi).\]
	Furthermore, if $\mod{M}$ is crisp, so is $\mod{M}^\times$.
	\item for $\mod{N}$ a 
	first order model in the language $\{R/2, \{E_p/1\colon p\text{ variable in }\mathcal{V}\}\}$, it holds that \[e_\times(v, \varphi) = \vert \langle \varphi, x_0\rangle^\sharp\vert_{\mod{N}, [x_0 \mapsto v]}.\] 
	Furthermore, if $\mod{N}$ is witnessed, quasi-witnessed or crisp, so is $\mod{N}_\times$.
	\end{enumerate}
	\end{obs}	
	
	We introduce the following translation of modal formulas (and analogously, of sets of formulas) into sentences on the predicate language  $\{R/2, \{E_p/1\colon p\text{ variable in }\varphi\}, c/0\}$:
	\[\varphi^l \coloneqq \langle \varphi, x_0\rangle^\sharp[x_0/c]\]
As before, 	for a set of sentences $\Gamma \cup \{\varphi\}$, we write $\Gamma \Vdash_{\mod{N}} \varphi$ to denote that either $\vert \Gamma \vert_{\mod{N}} < 1$ or $\vert \varphi \vert_{\mod{N}} = 1$.

	The above observation immediately implies the following.
	
\begin{proposition}\label{prop:modalToFo}
Let $\Gamma \cup \{\varphi\}$ be a finite set of modal formulas. 
For $\mod{M}$ a  Kripke model and $\mod{N}$ a first order model over the language $\{R/2, \{E_p/1\colon p\text{ variable in }\varphi\}, c/0\}$ the following conditions hold:
\[\Gamma^l \Vdash_{\mod{M}^\times} \varphi^l \text{ if and only if }\Gamma \Vdash_{\mod{M}} \varphi \qquad \text{ and }\qquad
\Gamma \Vdash_{\mod{N}_\times} \varphi \text{ if and only if  } \Gamma^l \Vdash_{\mod{N}} \varphi^l.\]
\end{proposition}

Fixed a predicate language $\mathcal{P}$, we let $\forall	\Pi$ be the logic arising from all safe $\alg{A}$ predicate structures, for $\alg{A} \in \class{P}$; we let $\forall[0,1]_\Pi$ be the logic arising from all safe $[0,1]_\Pi$ predicate structures. Formally, for $\Gamma \cup \{\varphi\}$ sentences in $\mathcal{P}$, we write \footnote{We consider the logic between sentences only, since this is the one will use later on.}
\begin{align*}
\Gamma \Vdash_{\forall \Pi} \varphi &\text{ iff for any }\alg{A} \in \class{P}, \text{ and any }\mod{M}\ \alg{A}\text{-safe structure for }\mathcal{P}, \Gamma \Vdash_{\mod{M}} \varphi,\\
\Gamma \Vdash_{\forall [0,1]_\Pi} \varphi &\text{ iff for any }\mod{M}\ [0,1]_\Pi\text{-safe structure for }\mathcal{P}, \Gamma \Vdash_{\mod{M}} \varphi.\\
\end{align*}

\section{Decidability and standard completeness of crisp local modal product logic $\Vdash_{K\Pi}$}\label{sec:statements}

Decidability of modal product logics remains, in most cases, an open question in the literature. It is known \cite{Ha98} that first-order product logics are not complete with respect to either finite models or witnessed models. This incompleteness holds both for the logic arising from the so-called \textit{general semantics}—defined over first order models valued in all linearly ordered algebras—and for the \textit{standard} semantics—defined over first order models valued in the standard product algebra. Consequently, it is not possible, at least in any straightforward way, to use finite models for working with these logics. This stands in contrast, for example, to the \L ukasiewicz case \cite{Vi20}.

The lack of the finite model property also arises in first-order G\"odel logic. However, concerning the modal G\"odel logic\footnote{Modal G\"odel logic is defined analogously to $\Vdash_{M[0,1]\Pi}$ and $\Vdash_{K[0,1]_\Pi}$. In that case, it is immediate to see that the logic of Kripke models evaluated over all G\"odel algebras coincides with that of models over the standard G\"odel algebra.}, it is possible to find a semantics that enjoys the finite model property \cite{CaMe13, CaMeRo17}, relying on the very particular characteristics of G\"odel homomorphisms. However, this approach cannot be adapted to prove a finite model property even for the modal logic over the standard product algebra.\footnote{For the interested reader, we note that, while in the G\"odel case it is relatively clear how to reconstruct a (possibly infinite) $[0,1]$-valued G\"odel Kripke model-via order-preserving mappings, which are G\"odel homomorphisms-, and this procedure cannot be carried out for the product case.} 

Regarding the logic $\Vdash_{M[0,1]\Pi}$, we can draw inspiration from the work developed in \cite{CeEs22} within the context of Description Logics. The results presented there establish the decidability of the set of theorems of the standard product logic with valued accessibility relations. In other words, the question of whether $\Vdash_{M[0,1]\Pi} \varphi$ for an arbitrary formula $\varphi$ is decidable. As we will see in the final section, the ideas from that previous work can be generalized to prove the decidability of the logical entailment $\Vdash_{M[0,1]\Pi}$, and furthermore, to show that $\Vdash_{M[0,1]\Pi}$ coincides with $\Vdash_{M\Pi}$.

However, despite being inspiring, these previous ideas do not offer a solution, nor a direction in which to work to solve the analogous questions regarding logics of crisp models. In the earlier reference, the method used to reconstruct a (possibly infinite) $[0,1]$-valued Kripke model crucially relies on modifying the value of the accessibility relation and freely moving it within the interval $(0,1)$. This freedom is, of course, not possible when the accessibility relation is expected to be crisp (i.e., taking values precisely in $\{0,1\}$).

One of the main results we present in this paper is the following one, closing the above open question of decidability of $\Vdash_{K[0,1]_\Pi}$ in a positive way.

\begin{theorem}[Decidability of $\Vdash_{K[0,1]_\Pi}$] \label{th:dec}
For any $\Gamma \cup \{\varphi\}\subseteq_{\omega}Fm$, the problem of determining whether or not 
\[\Gamma \Vdash_{K[0,1]_\Pi} \varphi\]
is decidable
\end{theorem}
In particular, the set of theorems of $\Vdash_{K[0,1]_\Pi}$ is decidable too.

We will see that, in fact, there is a second problem which turns to be strongly related to the above one. Namely, whether the logics $\Vdash_{K[0,1]_\Pi}$ and $\Vdash_{K\Pi}$ coincide. The second main result we present in this paper answers again positively this question.
We refer to it by \textit{standard completeness of (local) modal product logic}, following the terminology of the propositional standard completeness.

\begin{theorem}[Standard Completeness of $\Vdash_{K\Pi}$] \label{th:standardComp}
For any $\Gamma \cup \{\varphi\}\subseteq_{\omega}Fm$,
\[\Gamma \Vdash_{K\Pi} \varphi \qquad \text{ if and only if }\qquad \Gamma \Vdash_{K[0,1]_\Pi} \varphi.\]
\end{theorem}
In particular, we get that their sets of theorems also coincide.
Of course, from the two previous theorems, we also get that the logic $\Vdash_{K\Pi}$ is decidable.

The two results above are immediate corollaries of the next theorem, which will be proven along the following sections. 

%
%

%
\begin{theorem}\label{th:all}Let $\Upsilon = \Gamma \cup \{\varphi\} \subseteq_\omega Fm$. There is a propositional formula $\Upsilon^\star$ that can be effectively constructed from $\Upsilon$ for which the following are equivalent:\footnote{Recall that $\psi_\alpha$ is a propositional formula defined according to Definition \ref{def:modalToPropFormulas}.}

\begin{enumerate}
\item $\Gamma \Vdash_{K\Pi}\varphi$, 
\item $\Gamma_{\langle 0 \rangle}, \Upsilon^\star \models_{\Pi_\Delta} \varphi_{\langle 0 \rangle}$,
\item $\Gamma \Vdash_{K[0,1]_\Pi}\varphi$.
\end{enumerate}
\end{theorem}


It is immediate  that \emph{(1)}$\Rightarrow $\emph{(3)}, since $K[0,1]_\Pi \subset K\Pi$. In Section \ref{sec:soundness} we will prove that \emph{(2)} $\Rightarrow $ \emph{(1)}, and in Section \ref{sec:completeness}, that \emph{(3)} $\Rightarrow $ \emph{(2)}.  

Theorem \ref{th:standardComp} is the equivalence of \emph{(1)} and  \emph{(3)} from Theorem \ref{th:all}. On the other hand, Theorem \ref{th:dec} follows from the equivalence of \emph{(1)} and \emph{(2)}, and the fact that $\models_{\Pi_\Delta}$, i.e.,  (propositional) product logic with $\Delta$, is decidable (see e.g. \cite{Ha11}).

Before proceeding with the proof of the previous theorem, we would like to emphasize the role of $\Pi_\Delta$ (rather than $\Pi$) as the propositional logic that serves as the central point of the reduction. First, note that there is no essential difference between using $\Pi_\Delta$ or $\Pi$ when it comes to proving the decidability and standard completeness of $\Vdash_{K\Pi}$.
That said, we introduce $\Delta$ for technical convenience. In particular, in the proof of \emph{(2)} $\Rightarrow$ \emph{(1)}, the presence of $\Delta$ allows for a simpler definition of the set of formulas $\Upsilon^\star$. Nevertheless, its use is not strictly necessary, and it could be avoided. However, omitting $\Delta$ significantly complicates the definition of $\Upsilon^\star$ in Definition \ref{def:definingFormulas} and makes Theorem \ref{th:standardComp} much harder to follow. For this reason, we choose to work with $\Pi_\Delta$.

\section{From $K\Pi$ to propositional product logic}\label{sec:soundness}

As we said before, first-order product logic is not complete with respect to witnessed models. Nevertheless, it is known to be complete with respect to so-called \textit{closed} or \textit{quasi-witnessed} models \cite{LaMa07,CeEs11}. These models may require infinitely-branching nodes, but only in cases where the evaluation of a formula beginning with a box at the parent node drops to the bottom element of the algebra through an infinite descending chain of values at the successor nodes.

More specifically, and more relevant for the ongoing work, \cite{LaMa07} proves that the logic is complete with respect to quasi-witnessed models evaluated over a particular product algebra. Let us now introduce this algebra, present the formal completeness result, and fix some notation concerning the algebra introduced.

%
%
%
%

\begin{definition}\label{def:lexsum}
The \termDef{lexicographic sum} $\alg{\mathds{R}^{\mathds{Q}}} = \langle \mathds{R}^{\mathds{Q}},  \oplus,  \leq \rangle$ is the ordered abelian group of functions $f \colon \mathds{Q}\rightarrow \mathds{R}$ whose support is well-ordered (namely,  such that $\{q \in \mathds{Q}\colon f(q) \neq 0\}$ is a well-ordered subset of $\mathds{Q}$).  The $\oplus$ operation is the component-wise sum and the order on $\alg{\mathds{R}^{\mathds{Q}}}$ is the lexicographic one.
\end{definition}

The transformation $\mathfrak{B}$ introduced in \cite{CiT00} can be applied to the previous l-group,  obtaining a product chain in the following way. 
Let us denote $(\alg{\mathds{R}^{\mathds{Q}}})^- \coloneqq \{a \in \alg{\mathds{R}^{\mathds{Q}}}\colon a \leq \mathbf{0}\}$, where $\mathbf{0}$ stands for the neutral element of the group (i.e., the constant function $0$). 

\begin{definition}
We denote by $\mathbf{R}$ the product algebra  $\mathfrak{B}(\alg{\mathds{R}^{\mathds{Q}}})$, namely
\[\mathbf{R} \coloneqq \langle (\alg{\mathds{R}^{\mathds{Q}}})^- \cup \{\perp\}, \wedge,\vee, +, \rightarrow, \perp, \mathbf{0} \rangle,\]
where the order is inherited from $\alg{\mathds{R}^{\mathds{Q}}}$ (for the definition of $\wedge,\vee$) with $\bot$ the bottom element and $\mathbf{0}$ the top one, and 
\begin{align*}
x + y &= \begin{cases} x \oplus y & \hbox{ if }x,y \in (\alg{\mathds{R}^{\mathds{Q}}})^-,\\ \perp &\hbox{ otherwise. } \end{cases}& x \rightarrow y &= \begin{cases} \mathbf{0} \wedge (y - x)  & \hbox{ if }x,y \in (\alg{\mathds{R}^{\mathds{Q}}})^-, \\
\mathbf{0} &\hbox{ if } x = \perp, \\
\perp &\hbox{ if } y = \perp \text{ and } x\in (\alg{\mathds{R}^{\mathds{Q}}})^-, \end{cases}
\end{align*}
where $-$ is defined component-wise.
\end{definition}
%
%
%
%

For an element $a \in \mathbf{R}$ and $q \in \mathds{Q}$, we will let $a_{\mapsfrom q} \in \mathbf{R}$ to be defined as $\perp$ if $a = \perp$ and, otherwise (namely, if $a \in (\alg{\mathds{R}^{\mathds{Q}}})^-$), the element such that \[a_{\mapsfrom q}(p) = \begin{cases} a(p) &\hbox{ if }p \leq q, \\ 0 & \hbox{ otherwise. }\end{cases}\]
Namely, $a_{\mapsfrom q} $ is the element of the algebra $\mathbf{R}$ which coincides with $a$ in all elements of the domain smaller or equal than $q$, and it is cropped to $0$ in the rest.

\begin{obs}\label{obs:ordermapsto}
Let $a,b \in \mathbf{R}$. Then the following conditions hold:
\begin{enumerate}
\item For any $q \in \mathds{Q}$ it holds that 
 $a + a_{\mapsfrom q} = \perp$ if and only if $a = \perp$, and $\mathbf{0}_{\mapsfrom q} = \mathbf{0}$,
\item If $a \leq b$, then $a_{\mapsfrom q} \leq b_{\mapsfrom q}$ for any $q \in \mathds{Q}$,
\item $a \leq b$ if and only if $a + a_{\mapsfrom q} \leq b + b_{\mapsfrom q}$,
\item For $a,b \in (\alg{\mathds{R}^{\mathds{Q}}})^-$ it holds that $a_{\mapsfrom q} + b_{\mapsfrom q} = (a+b)_{\mapsfrom q}$, and for $a > b$, $a_{\mapsfrom q} - b_{\mapsfrom q} = (a-b)_{\mapsfrom q}$.
\end{enumerate}
\end{obs}
\begin{proof} See Appendix 2.\end{proof}

Let us also check another fact, related to the computation of values in quasi-witnessed situations.
\begin{lemma}\label{lem:infSupmapsfrom}
Let $X \subseteq \mathbf{R}$ such that $\bigwedge X = \bot$ or $\bigwedge X \in X$, and $\bigvee X \in X$. For any $q \in \mathds{Q}$ the following hold:
\[
\bigwedge_{x \in X}(x + x_{\mapsfrom q}) =  \bigwedge_{x \in X}x + (\bigwedge_{x \in X} x)_{\mapsfrom q}\qquad \text{ and} \qquad 
\bigvee_{x \in X}(x + x_{\mapsfrom q}) =  \bigvee_{x \in X}x + (\bigvee_{x \in X} x)_{\mapsfrom q}.
\]
\end{lemma}
\begin{proof} See Appendix 2.\end{proof}

For any $a \in \mathbf{R}$, we let ${\vartriangleleft}a \in \{\perp'\} \cup \mathds{Q} \cup \{\top'\}$ (where the order is the one from $\mathds{Q}$ and $\bot' < q < \top'$ for any $q \in \mathds{Q}$)  be defined as 
\[{\vartriangleleft}a \coloneqq 	\begin{cases} \perp' &\hbox{ if } a = \perp, \\ \top' &\hbox{ if } a = \mathbf{0} \\
min\{q \in \mathds{Q}\colon a(q) \neq 0\} &\hbox{ otherwise.}\end{cases}\] 
Namely, ${\vartriangleleft}a$ is the mimimum element of the domain (i.e., of $\mathds{Q}$) for which its image is different from $0$ (or $\bot', \top'$ if $a$ was, respectively, the bottom or the top element of the algebra). 
Note that this element ${\vartriangleleft}a$ is well defined, since for any $a \in \mathbf{R}$ such that $\perp < a < \mathbf{0}$, the set $\{q \in \mathds{Q}\colon a(q) \neq 0\}$ is a non-empty well-ordered subset of $\mathds{Q}$.
The following are easy observations which follow immediately from the definitions.
\begin{obs}\label{obs:closetriangle}
\begin{enumerate}
\item 
For any $a \in \mathbf{R}, q \in \mathds{Q}$ with ${\vartriangleleft}a \leq q$ it holds that ${\vartriangleleft}a_{\mapsfrom q} = {\vartriangleleft}a$, and if $a > \bot$, $a_{\mapsfrom {\vartriangleleft}a}[{\vartriangleleft}a] = a[{\vartriangleleft}a]$.
\item For any $a,b \in \mathbf{R}$, if ${\vartriangleleft}a \leq {\vartriangleleft}b$ then ${\vartriangleleft}(a + b) =  {\vartriangleleft}a$. 
\item For any $a,b \in \mathbf{R}$, if ${\vartriangleleft}a < {\vartriangleleft}b$ then $b_{\mapsfrom {\vartriangleleft}a} = \mathbf{0}$, and ${\vartriangleleft}(a + a_{\mapsfrom{q}}) < {\vartriangleleft}(b + b_{\mapsfrom q})$.
\item For any $I\subseteq \omega$ and $a_i \in \alg{R}$ for each $i \in I$, \[\bigwedge_{i \in I}a_i = \bot \text{ if and only if } \forall q \in \mathds{Q}, \exists i_q \in I \text{ s.t. } {\vartriangleleft}a_{i_q} < q.\]
\end{enumerate}
\end{obs}

As we highlighted in the beginning of this section, first order logic $\forall\Pi$ is complete with respect to countable first order models valued over $\mathbf{R}$. Not only this, the authors show that it is complete with respect to quasi-witnessed models over that algebra.

\begin{theorem}[Theorem 2.9, \cite{LaMa07}]\label{Thm:LaMa07}
Let $\mathcal{P}$ be a countable predicate language,  and $\Gamma \cup \{\varphi\}$ a finite set of sentences in that language such that $\Gamma \not \vdash_{\forall\Pi} \varphi$.  Then there is a countable,  quasi-witnessed $\mathbf{R}$-structure $\mod{M}$ such that $\mod{M}\models \Gamma$ and $\mod{M} \not \models \varphi$. 
\end{theorem}

We can combine the previous result with the constructions and observations done in the Section \ref{subsec:prel-predLogic} relating modal and predicate logics, particularly with Proposition \ref{prop:modalToFo} and Observation \ref{obs:languageModalToFO}. Using also Corollary \ref{cor:LocalFiniteDepth} in order to obtain a tree of finite depth, we get the following corollary.
\begin{corollary}\label{cor:quasiWitComp}
Let $\Gamma \cup \{\varphi\}$ a set of modal formulas such that $\Gamma \not \vdash_{K\Pi} \varphi$.  Then there is a countable,  quasi-witnessed $\mathbf{R}$-Kripke tree  of finite depth $\mod{T}$ and root $v$ such that $\Gamma \not \Vdash_{\mod{T},v} \varphi$.
\end{corollary}
\begin{proof} See Appendix 2.\end{proof}

The study of $\mathbf{R}$-Kripke trees in a deeper detail will allow us to prove that, without loss of generality, we can refine completeness not only to $\mathbf{R}$-quasi-witnessed trees, but to those in which, for the quasi-witnessed nodes (i.e., those with a necessary infinite branching), the values taken in all the successors of all the formulas involved hold a regular relation with their respective value in a certain ``seed" successor node. This relation will turn to be expressible by means of a finite set of propositional formulas, allowing for the definition of the $\Upsilon^*$ formula from Theorem \ref{th:all}. In combination to the more usual encoding in the propositional language of the values of modal formulas in finitely-branching nodes as, respectively, minima or suprema over the values at the successors, this will allow us to conclude the implication \emph{(2)} $\Rightarrow$ \emph{(1)} from Theorem \ref{th:all} (by contraposition).

%
%
%
%
%

%

\begin{lemma}\label{lem:belowWorld}
Let $\mod{T}$ be a quasi-witnessed $\mathbf{R}$-Kripke tree, $v \in W$ and $\Upsilon$ a finite set of formulas. For every $\Box \varphi \in \Upsilon$ such that $e(v,\Box \varphi)=\bot$, there is some world\footnote{If $\Box \varphi \in uWit_{\mod{T}}(w, \Upsilon)$, i.e., $e(w, \varphi)>\bot$ for all $w \in W$ such that $Rvw$, there are in fact infinitely many.} $v_{\Box \varphi}$ such that $Rvv_{\Box \varphi}$ and 
\[{\vartriangleleft}e(v_{\Box \varphi}, \varphi) < \bigwedge \{{\vartriangleleft}e(v_{\Box \varphi}, \chi)\colon \Box \chi \in \Upsilon, \\ e(v, \Box \chi) > \bot\}.\]
\end{lemma}
\begin{proof}
 Since $\Upsilon$ is a finite set, and by the definition of the evaluation of the $\Box$ operator, we have 
\[\bot < \bigwedge \{e(v, \Box \chi)\colon \Box \chi \in \Upsilon, e(v, \Box \chi) > \bot\} \leq \bigwedge \{e(w, \chi)\colon \Box \chi \in \Upsilon, e(v, \Box \chi) > \bot \}\]
for all $w \in W$ s.t. $Rvw$.

From the definition of the $\vartriangleleft$ operation it follows that 
\[\bot < \bigwedge \{{\vartriangleleft}e(v, \Box \chi)\colon \Box \chi \in \Upsilon, e(v, \Box \chi) > \bot\} \leq \bigwedge \{{\vartriangleleft}e(w, \chi)\colon \Box \chi \in \Upsilon, e(v, \Box \chi) > \bot\}\]
for all $w \in W$ s.t. $Rvw$. 

Let us denote by $r$ the number $ \bigwedge \{{\vartriangleleft}e(v, \Box \chi)\colon \Box \chi \in \Upsilon, e(v, \Box \chi) > \bot\}$.

Since  by assumption $e(v,\Box \varphi) = \bot$, from Observation \ref{obs:closetriangle}(4) we get there is 
some $w_r \in W$ such that $Rvw_r$ and ${\vartriangleleft}e(w_r, \varphi) < r$. 
This proves the lemma, since
\[{\vartriangleleft}e(w_r, \varphi) < r \leq \bigwedge \{{\vartriangleleft}e(w_r, \chi)\colon \Box \chi \in \Upsilon, e(v, \Box \chi) > \bot\} .\qedhere \]

\end{proof}

Using the previous lemma, we can introduce the following convention. Consider a $\mathbf{R}$-Kripke tree $\mod{T}$ and a fixed finite set of formulas $\Upsilon$  (omitted from the notation). For every world $v \in W$ with at least one successor, and $\heartsuit \varphi \in \Upsilon$ a formula beginning with a modality, we let
\[v_{\heartsuit \varphi} \coloneqq \begin{cases} \text{ any world in }wits_{\mod{T}}(v, \heartsuit \varphi) &\hbox{ whenever this set is non-empty,}\\ \text{ any world identified in Lemma \ref{lem:belowWorld}} &\hbox{ otherwise.} \end{cases}\]
If $v$ has no successors, the assignment $v_{\heartsuit \varphi}$ is empty.

We will use the above labeling convention throughout a sequence of transformations. Roughly speaking, we first identify a distinguished subset of worlds in the tree $\mod{T}$, namely those labeled according to the previous convention. We then use this labeling to systematically rename the worlds of $\mod{T}$, thereby obtaining a canonical way of addressing them that will simplify the constructions developed later. Up to this stage, the process amounts merely to a relabeling of the original model.

In a final and more technical step, we expand the renamed model in a conservative way, obtaining certain regularity properties related to unwitnessed situations.

We will refer to sequences by letters $\sigma, \delta, \eta, \varrho$, and to single in them elements by letters $\alpha, \beta, \gamma$. The concatenation of two sequences $\sigma_1$ and $\sigma_2$ will be denoted by $\sigma_1^\smallfrown \sigma_2$. When $\sigma$ is a sequence and $\beta$ an element, we write $\sigma^\smallfrown \beta$ and $\beta^\smallfrown \sigma$ as a shorthand for $\sigma^\smallfrown \langle \beta \rangle$ and $\langle \beta \rangle^\smallfrown\sigma $. 

 

 Given three sequences of elements $\sigma_1, \sigma_2$ and $\delta$, the fact that $\sigma_1$ is an initial subsequence of $\delta$ by $\sigma_1 \sqsubseteq \delta$, and the sequence resulting of substituting the whole initial subsequence $\sigma_1$ by the sequence $\sigma_2$ in $\delta$ by $\delta[\sigma_1/\sigma_2]$.

\begin{definition}\label{def:relevanWorlds}
Let $\mod{T}$ be a quasi-witnessed $\mathbf{R}$-Kripke tree with root $r$, and $\Upsilon$ a finite set of formulas. We let the \termDef{relevant worlds} of $\mod{T}$ according to $\Upsilon$ to be the set $RW_{\Upsilon}(\mod{T}) \coloneqq \bigcup_{i \leq md(\Upsilon)} RW_{\Upsilon}^i(\mod{T})$ where

\begin{align*}
RW_{\Upsilon}^0(\mod{T}) \coloneqq& \{r\},\\
RW_{\Upsilon}^{i+1}(\mod{T}) \coloneqq& \{v_{\heartsuit \varphi} \in W\colon v \in RW_{\Upsilon}^{i}(\mod{T}), \heartsuit \varphi \in \Upsilon_{\depth(v)}  \}.
\end{align*}
\end{definition}

For a quasi-witnessed $\mathbf{R}$-Kripke tree $\mod{T}$ with root $r$, and $\Upsilon$ a finite set of formulas, let us recursively define the mapping $\eta$ assigning to each of the elements of its universe a set of particular sequences that will hold some information on the evaluation of the formulas of $\Upsilon$. In the following definition we make use of the characterization of the universe of a tree as $W=\bigcup_{i \in I}W_i$ following Definition \ref{def:tree}. Recall that $par(v)$ denotes the parent of $v$ in the tree.


\begin{align*}
\eta(r) \coloneqq & \{\langle 0 \rangle\},\\
\text{ for }v \in W_{i+1}, \eta(v) \coloneqq & \begin{cases} \{\sigma^\smallfrown\heartsuit \varphi\colon  \sigma \in \eta(par(v)), v = par(v)_{\heartsuit \varphi}\} &\hbox{ if } v \in RW_{\Upsilon}(\mod{T}),\\
\{\sigma^\smallfrown v \colon  \sigma \in \eta(par(v))\}&\hbox{ otherwise. }\end{cases}
\end{align*}

Intuitively, each relevant world is associated with a sequence indicating the ``path" (of witnesses or analogous labels) that lead to it from the root (and is replicated if some world is a witness of several formulas), while for the rest of the worlds, we associate to them the label of the last ancestor in the tree which was relevant, followed by the sequence of the original names of the worlds in $\mod{T}$. 
It is easy to see that $\eta(v) \cap \eta(w) = \emptyset$ for any $v \neq w$, so for any $\sigma \in \bigcup \{\eta(v) \colon v \in W_i, i \leq md(\Upsilon)\}$ we can write $\eta^{-1}(\sigma)$ to denote the (unique) element $v \in W$ such that $\sigma \in \eta(v)$.

We let $\eta_\Upsilon(\mod{T})$ be defined from $\mod{T}$ by using the previous tagging in the obvious way and cropping the model at depth $md(\Upsilon)$, namely, the model $\langle \eta_\Upsilon(W), \eta_\Upsilon(R), \eta_\Upsilon(e) \rangle$ where, 
\begin{align*}
\eta_\Upsilon(W) \coloneqq& \bigcup \{\eta(v) \colon v \in W_i, i \leq md(\Upsilon)\},\\
\eta_\Upsilon(R) \coloneqq& \{\langle \sigma, \sigma^\smallfrown \alpha\rangle\colon \sigma^\smallfrown \alpha \in \eta_\Upsilon(W) \}, \\
\eta_\Upsilon(e)(\sigma, p) \coloneqq& e(\eta^{-1}(\sigma),p) \text{ for all }p.
\end{align*}
For the sake of readability, we will remove the subindex $\Upsilon$ from the previous elements when it is clear from the context.

It is easy to see that for any $\sigma \in \eta(W)$, $\depth_{\eta(\mod{T})}(\sigma) = \depth_{\mod{T}}(\eta^{-1}(\sigma)) = \vert \sigma\vert -1$ (for $\vert \sigma \vert$ the number of elements of $\sigma$), and that 
$\eta^{-1}(\sigma ^\smallfrown \heartsuit \varphi) = \eta^{-1}(\sigma)_{\heartsuit \varphi}$. 
It is also clear that if $\sigma$ is a successor (not necessarily immediate) of $\delta$ in $\eta_\Upsilon(\mod{T})$, then $\delta \sqsubseteq \sigma$.

\begin{exa}\label{example:running1}
Let us introduce a detailed example which will be running through this section, on how the modifications we define affect a particular tree over a particular set of formulas. 

Instead of fully defining the model (which would turn to be rather cumbersome and possibly confusing), let us specify only its determining characteristics for what concerns the transformations we will apply to it.  
Consider the frame of a tree $\mod{T}$ to be the following one: 
\[\begin{tikzcd}[sep=small]
	{u_1} &&& r \\
	& {v_1} &&&& {} \\
	{u_2} &&& {v_k} \\
	& {w_1} &&&& {} \\
	&&& {w_m}
	\arrow[curve={height=6pt}, from=1-4, to=2-2]
	\arrow[from=1-4, to=3-4]
	\arrow[""{name=0, anchor=center, inner sep=0}, curve={height=12pt}, dotted, no head, from=3-4, to=2-6]
	\arrow[from=2-2, to=1-1]
	\arrow[from=2-2, to=3-1]
	\arrow[""{name=1, anchor=center, inner sep=0}, curve={height=6pt}, dotted, no head, from=4-2, to=5-4]
	\arrow[""{name=2, anchor=center, inner sep=0}, curve={height=6pt}, dotted, no head, from=5-4, to=4-6]
	\arrow[from=3-4, to=5-4]
	\arrow[curve={height=6pt}, from=3-4, to=4-2]
	\arrow[""{name=3, anchor=center, inner sep=0}, curve={height=12pt}, dotted, no head, from=2-2, to=3-4]
	\arrow[shorten >=4pt, Rightarrow, dotted, from=1-4, to=0]
	\arrow[shorten >=4pt, Rightarrow, dotted, from=3-4, to=2]
	\arrow[shorten >=8pt, dotted, from=3-4, to=1]
	\arrow[shorten >=4pt, dotted, from=1-4, to=3]
\end{tikzcd}\]

 By a dotted arrow we denote access to a finite number of accessible worlds, and by a double dotted arrow, to an infinite one.

 Assume that we are interested in the values taken at the root $r$ by the set of formulas $\Upsilon = \{\Box (y \vee \Box x), \Diamond \Diamond y\}$. Henceforth, $\Upsilon_0 = \Upsilon$, $\Upsilon_1 = \{y, \Box x, \Diamond y\}$ and $\Upsilon_2 = \{x, y\}$. Suppose $e(r, \Box (y \vee \Box x))= \bot$ unwitnessed, and $e(v_k, \Box x) = \bot$ unwitnessed, and such that $r_{\Box (y \vee \Box x)} = v_k$ and $(v_k)_{\Box x} = w_m$ (the worlds identified in Lemma \ref{lem:belowWorld}). For the other formulas, assume witnessed evaluations, and let 
 \begin{align*}
 r_{\Diamond \Diamond y} &= v_1, &
 (v_1)_{\Diamond y} &= u_1,\\
 (v_1)_{\Box x} &= u_1, &
 (v_k)_{\Diamond y} &= w_1.
 \end{align*}

 We can first clearly identify the set of relevant worlds, which are the ones we previously introduced: 
 \[RW_{\Upsilon}(\mod{T}) = \{r, v_1, v_k, u_1, w_1, w_m\}.\]

 This leads to obtain the following structure of the model $\eta(\mod{T})$, where we have depicted in gray the worlds coming from Lemma \ref{lem:belowWorld}:
\[\begin{tikzcd}[sep=small]
	{\langle 0, \Diamond \Diamond y, \Diamond y\rangle} &&& {\langle 0 \rangle} \\
	{\langle 0, \Diamond \Diamond y, \Box x\rangle} & {\langle 0, \Diamond \Diamond y\rangle} &&&& {} \\
	{\langle 0, \Diamond \Diamond y, u_2\rangle} &&& \textcolor{rgb,255:red,120;green,120;blue,120} {\langle 0, \Box(y \vee \Box x)\rangle} \\
	& {\langle 0, \Box(y \vee \Box x), \Diamond y\rangle} &&&& {} \\
	&&&  \textcolor{rgb,255:red,120;green,120;blue,120}{\langle 0, \Box(y \vee \Box x), \Box x\rangle}
	\arrow[curve={height=6pt}, from=1-4, to=2-2]
	\arrow[from=1-4, to=3-4]
	\arrow[from=2-2, to=1-1]
	\arrow[from=2-2, to=2-1]
	\arrow[from=2-2, to=3-1]
	\arrow[""{name=0, anchor=center, inner sep=0}, curve={height=12pt}, dotted, no head, from=2-2, to=3-4]
	\arrow[""{name=1, anchor=center, inner sep=0}, curve={height=12pt}, dotted, no head, from=3-4, to=2-6]
	\arrow[curve={height=6pt}, from=3-4, to=4-2]
	\arrow[from=3-4, to=5-4]
	\arrow[""{name=2, anchor=center, inner sep=0}, curve={height=6pt}, dotted, no head, from=4-2, to=5-4]
	\arrow[""{name=3, anchor=center, inner sep=0}, curve={height=6pt}, dotted, no head, from=5-4, to=4-6]
	\arrow[shorten >=5pt, Rightarrow, dotted, from=1-4, to=1]
	\arrow[shorten >=6pt, dotted, from=1-4, to=0]
	\arrow[shorten >=5pt, Rightarrow, dotted, from=3-4, to=3]
	\arrow[shorten >=11pt, dotted, from=3-4, to=2]
\end{tikzcd}\] 
 
 \qed

\end{exa}

\begin{lemma}\label{lem:inverseEtaPreservation}
Let $\mod{T}$ be a quasi-witnessed $\alg{R}$-Kripke tree. For any $\sigma \in \eta(W)$ and any formula $\varphi$ it holds that,
$\eta(e)(\sigma, \varphi) = e(\eta^{-1}(\sigma), \varphi).$
In particular, for each $\sigma'$ such that $\eta^{-1}(\sigma') = \eta^{-1}(\sigma)$ it holds that $\eta(e)(\sigma, \varphi) = \eta(e)(\sigma', \varphi)$. 
\end{lemma}
\begin{proof} See Appendix 2.\end{proof}

\begin{corollary}\label{cor:inverseEtaPreservation}
Let $\mod{T}$ be a quasi-witnessed $\alg{R}$-Kripke tree. For any $\sigma \in \eta(W)$ and any formula $\Box \varphi \in \Upsilon_{d(\eta^{-1}(\sigma))}$ with $\Box \varphi \in uWit_{\mod{T}}(\eta^{-1}(\sigma), \Upsilon_{d(\eta^{-1}(\sigma))})$, the world $\sigma^\smallfrown\Box \varphi$ satisfies (in $\eta(\mod{T})$) the condition from Lemma \ref{lem:belowWorld}, namely, 
\[{\vartriangleleft}\eta(e)(\sigma^\smallfrown\Box \varphi, \varphi) < \bigwedge \{{\vartriangleleft}\eta(e)(\sigma^\smallfrown\Box \varphi, \chi) \colon \Box \chi \in \Upsilon_{d(\eta^{-1}(\sigma))}, \eta(e)(\sigma, \Box \chi) > \bot\}.\]
\end{corollary}

As we anticipated, we will now extend the tree $\eta(\mod{T})$ with additional sequences in which some modal formulas (of those starting with a $\Box$, as we will see) appear primed, e.g., $\delta = \langle 0, \Diamond \psi_1, \Box \psi_2', \Box \psi_3, \Box \psi_4'\rangle$. Intuitively, what we are doing is replicating subtrees obtained from worlds of the form $v_{\Box \chi}$ where $\Box \chi$ is unwitnessed in $v$, modifying the evaluation there in a particular way (depending precisely on the value taken by $\chi$ in $v_{\Box \chi}$) and gluing these subtrees in the natural position withing the original tree.

Define the set \[Gens_i^\Upsilon(\mod{T}) \coloneqq \{\langle \sigma, \Box \varphi\rangle \colon \sigma \in \eta(RW(\mod{T})), \vert \sigma \vert = i +1, \Box \varphi \in uWit_{\eta(\mod{T})}(\sigma, \Upsilon_i)\}.\] These are, for $i$ indicating the depth of the model (i.e., distance from the root), the pairs formed by a world and the unwitnessed formulas at that world, at depth $i$ (of the formulas of $\Upsilon$ of the corresponding modal depth). 
Observe that, according to the previous notation, for any $\langle \sigma, \Box \varphi\rangle \in Gens_i$, we have that the world $\sigma^\smallfrown\Box \varphi$ belongs to $\eta(RW(\mod{T}))$.

\begin{definition}\label{def:mod+}
Let $\mod{T}$ be a quasi-witnessed $\mathbf{R}$-Kripke tree with root $r$ and $\Upsilon$ a finite set of formulas. 
We will define $\mod{T}^0, \ldots, \mod{T}^{md(\Upsilon)}$ as follows. First, let $\mod{T}^{md(\Upsilon)} = \langle W^{md(\Upsilon)}, R^{md(\Upsilon)}, e^{md(\Upsilon)}\rangle \coloneqq \eta_\Upsilon(\mod{T})$. Now, for $\mod{T}^i$ with $md(\Upsilon) \geq i > 0$ already defined, let\footnote{For the sake of readability, the recursion is given in inverse order, namely from the largest -i.e., $md(\Upsilon)$- to the lowest index ($0$), because the worlds considered at each step go from the deepest level to the root of the model.}
$\mod{T}^{i-1} \coloneqq \langle W^{i-1},R^{i-1},e^{i-1}\rangle$ for

 \begin{align*}
 W^{i-1} \coloneqq& W^i \cup \bigcup\{W^i_{\langle \sigma, \Box \varphi\rangle}\colon \langle \sigma, \Box \varphi\rangle \in Gens_i^\Upsilon(\mod{T})\}, \\
 R^{i-1} \coloneqq& R^i \cup \bigcup\{R^i_{\langle \sigma, \Box \varphi\rangle} \cup \{ \langle \sigma, \sigma ^\smallfrown\Box \varphi'\rangle \colon \langle \sigma, \Box \varphi\rangle \in Gens_i^\Upsilon(\mod{T})\}\} ,\\
   e^{i-1}(\sigma,p) \coloneqq & \begin{cases} 
   e^i(\sigma,p) &\hbox{ if } \sigma \in W^i,\\ 
   e^i_{\langle \sigma, \Box \varphi\rangle}(\sigma,p) &\hbox{ if } \sigma \in W^i_{\langle \sigma, \Box \varphi\rangle} \text{ for some } \langle \sigma, \Box \varphi\rangle \in Gens_i^\Upsilon(\mod{T}).\end{cases}
 \end{align*}

with 
\begin{itemize}
\item $W^i_{\langle \sigma, \Box \varphi\rangle} \coloneqq \{\delta[\sigma^\smallfrown\Box \varphi/\sigma^\smallfrown\Box \varphi']\colon \delta \in W[\mod{T^i}_{\sigma^\smallfrown\Box \varphi}]\}$,\footnote{Recall that $\mod{T^i}_{\sigma^\smallfrown\Box \varphi}$ denotes the subtree of $\mod{T^i}$ generated from the world $\sigma^\smallfrown\Box \varphi$, and that $\delta[\sigma^\smallfrown\Box \varphi/\sigma^\smallfrown\Box \varphi']$ is the sequence resulting from substituting, in $\delta$, the initial subsequence $\sigma^\smallfrown\Box \varphi$ by $\sigma^\smallfrown\Box \varphi'$ -namely, changing the corresponding apparition of $\Box \varphi$ by its primed variant. It is clear that, by construction, $\sigma^\smallfrown\Box \varphi$ will be an initial sequence of each world in $\mod{T^i}_{\sigma^\smallfrown\Box \varphi}$.}

\item $R^i_{\langle \sigma, \Box \varphi\rangle} \coloneqq \{\langle \delta_1[\sigma^\smallfrown\Box \varphi/\sigma^\smallfrown\Box \varphi'], \delta_2[\sigma^\smallfrown\Box \varphi/\sigma^\smallfrown\Box \varphi']\rangle\colon \langle \delta_1, \delta_2 \rangle \in R^i\}$,
\item $e^i_{\langle \sigma, \Box \varphi\rangle}(\delta[\sigma^\smallfrown\Box \varphi/\sigma^\smallfrown\Box \varphi'], p) \coloneqq e^i(\delta, p) + e^i(\delta, p)_{\mapsfrom{{\vartriangleleft}e^i(\sigma^\smallfrown\Box \varphi, \varphi)}}.$

\end{itemize}

We define \[\mod{T}^+ \coloneqq \mod{T}^0.\]

For a model $\mod{T}^+$, we call \termDef{Hereditarily  Relevant}, and denote them by $HRW(\mod{T}^+)$, to the worlds in $\mod{T}^+$ that are named with a sequence of the form $\langle 0, \alpha_1, \ldots, \alpha_n\rangle$ with $\alpha_i$ a modal formula or a primed modal formula. Namely, those that are in $\eta(RW(\mod{T}))$, or copies of these where some formula is primed.

\end{definition}

In the following, in order to lighten the notation, we assume a fixed quasi-witnessed $\alg{R}$ tree $\mod{T}$, and will enunciate results about its extension $\mod{T}^+$ using the usual notation $\mod{T}^+ = \langle W^+, R^+, e^+\rangle$ and the notation (for the construction steps) used in the previous definition.

We are going to extend the notation about sequences to handle specific questions about the sequences with primed elements introduced in the previous extension. 
For any sequence $\sigma$, we will denote by $\sigma_-$ the sequence in which the \textit{first} prime (if any) has been removed, and by $\underline{\delta}$ the sequence in which \textit{all} primes (if any) have been removed (so it is, in fact an element of $\eta(W)$). For instance, for $\delta = \langle 0, \Diamond \psi_1, \Box \psi_2', \Box \psi_3, \Box \psi_4'\rangle$, the previous lets $\delta_- = \langle 0, \Diamond \psi_1, \Box \psi_2, \Box \psi_3, \Box \psi_4'\rangle$ and $\underline{\delta} = \langle 0, \Diamond \psi_1, \Box \psi_2, \Box \psi_3, \Box \psi_4\rangle$.

Furthermore, for $\sigma = \sigma_-$ (namely, without primed elements) and any sequence $\delta$ we will define:
\[c(\sigma^\smallfrown\Box \varphi'^\smallfrown \delta) = \sigma^\smallfrown\Box \varphi \qquad \text{ and } \qquad
{pr}(\sigma^\smallfrown\Box \varphi'^\smallfrown \delta) = \varphi.\]
Following the previous example, we would get that $c(\delta) = \langle 0, \Diamond \psi_1, \Box \psi_2\rangle$ and ${pr}(\delta) = \psi_2$.

It is clear that if $\sigma \in W^{i-1}\setminus W^i$ for some $1 \leq i \leq md(\Upsilon)$, both $\sigma_-$ and $c(\sigma)$ belong to $W^i$.\footnote{Indeed, observe that in this case, $\sigma = \delta[\xi ^\smallfrown\Box \varphi/\xi ^\smallfrown\Box \varphi']$ for some corresponding $\langle \xi, \Box\varphi \rangle $ and $\delta \in W[\mod{T}^i_{\sigma ^\smallfrown\Box \varphi}] \subseteq W^i$. Hence $\sigma_- = \delta$, and $c(\sigma) = \xi ^\smallfrown\Box \varphi \sqsubseteq \delta$.}
  Therefore, in Definition \ref{def:mod+}, observe that for $\sigma, \delta \in W^{i-1}\setminus W^i$ (i.e., new worlds added at step $i-1$), we have that $R^{i-1} \sigma \delta$ if and only if $R^i \sigma_-\delta_-$, and $e^{i-1}(\sigma, p) = e^i(\sigma_-, p) + e^i(\sigma_-, p)_{\mapsfrom{{\vartriangleleft}e^i(c(\sigma), {pr}(\sigma))}}$.

\begin{exa}\label{example:running2}
Let us continue Example \ref{example:running1}, and see how that model is extended to $\mod{T}^+$. The whole model is $\mod{T}^0$, and we have depicted in black the original model $\mod{T}$, hence $\mod{T}^2$. To go from $\mod{T}^2$ to $\mod{T}^1$ we add the worlds and relations colored in green, and to go from $\mod{T}^1$ to $\mod{T}^0$ we add the worlds and relations colored in red.  

\[\begin{tikzcd}[cramped, scale cd=0.73]
	&&& {\langle 0 \rangle} && \textcolor{rgb,255:red,190;green,14;blue,17}{\langle 0, \Box(y \vee \Box x)', \Diamond y\rangle} \\
	{\langle 0, \Diamond \Diamond y, \Diamond y\rangle} \\
	{\langle 0, \Diamond \Diamond y, \Box x\rangle} & {\langle 0, \Diamond \Diamond y\rangle} &&&& \textcolor{rgb,255:red,190;green,14;blue,17}{\langle 0, \Box(y \vee \Box x)', \Box x\rangle} \\
	{\langle 0, \Diamond \Diamond y, u_2\rangle} && {\langle 0, \Box(y \vee \Box x)\rangle} && \textcolor{rgb,255:red,190;green,14;blue,17}{\langle 0, \Box(y \vee \Box x)'\rangle} \\
	& {\langle 0, \Box(y \vee \Box x), \Diamond y\rangle} && {} && {} \\
	&&&&& \textcolor{rgb,255:red,190;green,14;blue,17}{\langle 0, \Box(y \vee \Box x)', \Box x'\rangle} \\
	& {\langle 0, \Box(y \vee \Box x), \Box x\rangle} &&& \textcolor{rgb,255:red,38;green,107;blue,30}{\langle 0, \Box(y \vee \Box x), \Box x'\rangle} \\
	&& {}
	\arrow[curve={height=6pt}, from=1-4, to=3-2]
	\arrow[from=1-4, to=4-3]
	\arrow[draw={rgb,255:red,190;green,14;blue,17}, from=1-4, to=4-5]
	\arrow[""{name=0, anchor=center, inner sep=0}, draw={rgb,255:red,190;green,14;blue,17}, dotted, no head, from=1-6, to=3-6]
	\arrow[from=3-2, to=2-1]
	\arrow[from=3-2, to=3-1]
	\arrow[from=3-2, to=4-1]
	\arrow[""{name=1, anchor=center, inner sep=0}, curve={height=12pt}, dotted, no head, from=3-2, to=4-3]
	\arrow[""{name=2, anchor=center, inner sep=0}, draw={rgb,255:red,190;green,14;blue,17}, dotted, no head, from=3-6, to=5-6]
	\arrow[curve={height=6pt}, from=4-3, to=5-2]
	\arrow[""{name=3, anchor=center, inner sep=0}, curve={height=12pt}, dotted, no head, from=4-3, to=5-4]
	\arrow[from=4-3, to=7-2]
	\arrow[draw={rgb,255:red,38;green,107;blue,30}, curve={height=30pt}, from=4-3, to=7-5]
	\arrow[draw={rgb,255:red,190;green,14;blue,17}, curve={height=-12pt}, from=4-5, to=1-6]
	\arrow[draw={rgb,255:red,190;green,14;blue,17}, from=4-5, to=3-6]
	\arrow[draw={rgb,255:red,190;green,14;blue,17}, curve={height=12pt}, from=4-5, to=6-6]
	\arrow[""{name=4, anchor=center, inner sep=0}, curve={height=6pt}, dotted, no head, from=5-2, to=7-2]
	\arrow[""{name=5, anchor=center, inner sep=0}, curve={height=6pt}, dotted, no head, from=7-2, to=8-3]
	\arrow[shorten >=10pt, dotted, from=1-4, to=1]
	\arrow[shorten >=8pt, Rightarrow, dotted, from=1-4, to=3]
	\arrow[shorten >=14pt, dotted, from=4-3, to=4]
	\arrow[shorten >=8pt, Rightarrow, dotted, from=4-3, to=5]
	\arrow[draw={rgb,255:red,190;green,14;blue,17}, shorten >=15pt, from=4-5, to=0]
	\arrow[draw={rgb,255:red,190;green,14;blue,17}, shorten >=11pt, Rightarrow, dotted, from=4-5, to=2]
\end{tikzcd}\]
  \qed
\end{exa}

Let us now prove some facts about $\mod{T}^+$ that will allow us to understand better how it is built and how it behaves.

It is first immediate that at each step $i-1$ in the construction of $\mod{T}^+$, at any point $\delta \in \mod{T}^{i}$, we add (if any) a finite number of immediate successors to it in the construction of $\mod{T}^{i-1}$. Therefore, we know that each $\mod{T}^i$ (and so, also $\mod{T}^+$) is a quasi-witnessed tree too.

\begin{obs}\label{obs:freshSuccessorsOnlyOriginalWorlds}
\begin{enumerate}
\item $\delta \in W^{i-1}\setminus W^i$ if and only if $\delta = \delta_1^\smallfrown\Box \varphi'^\smallfrown\delta_2$ for some $\langle \delta_1, \Box \varphi \rangle \in Gens_i^\Upsilon(\mod{T})$.
\item For $\delta \in W^i, \delta^\smallfrown \alpha \in W^{i-1}\setminus W^i$ if and only if $\alpha = \Box \varphi'$ for some  $\langle \delta, \Box \varphi \rangle \in Gens_i^\Upsilon(\mod{T})$. 
\item For $\delta \in W^{i-1}\setminus W^i$ for some $1 \leq i \leq md(\Upsilon)$, it holds that 
${\delta_-}^\smallfrown \alpha \in W^{i-1}$ if and only if ${\delta_-}^\smallfrown \alpha \in W^{i}$ already. In particular for $\delta \in W^+$, $\delta ^\smallfrown \alpha \in W^+$ if and only if ${\delta_-}^\smallfrown \alpha \in W^+$, and if and only if $\underline{\delta}^\smallfrown \alpha \in W^+$.
\end{enumerate}
\end{obs}
\begin{proof} See Appendix 2.\end{proof}

The previous implies that only worlds belonging to the original model receive fresh immediate successors in the construction of $\mod{T}^+$.

It is now fairly easy to check a sort of truth-lemma extending the definition of the evaluation of the propositional letters in the new worlds of the model $\mod{T}^+$ to the whole model and all formulas. 

\begin{lemma}\label{lem:truthLemmaM+}
Let $0 < i \leq md(\Upsilon)$. For any $\sigma \in W^{i-1}\setminus W^i$ and any formula $\varphi \in \Upsilon_{\depth(\sigma)}$,

\[e^{i-1}(\sigma, \varphi) = e^i(\sigma_-, \varphi) + e^i(\sigma_-, \varphi)_{\mapsfrom{{\vartriangleleft}e^i(c(\sigma), {pr}(\sigma))}}.\]
\end{lemma} 
\begin{proof}
We will prove it for arbitrary $i$ by induction on the complexity of the formula. For propositional variables the lemma is exactly the definition of $e^{i-1}$. 

For propositional connectives, the proof follows easily by applying the conditions from Observation \ref{obs:ordermapsto}. 

For formulas starting with the modal operator $\Box$, by definition and applying I.H. we first have that 
\[e^{i-1}(\sigma, \Box \varphi) = \bigwedge_{\sigma^\smallfrown\alpha \in W^{i-1}} e^{i-1}(\sigma^\smallfrown\alpha, \varphi) = \bigwedge_{\sigma^\smallfrown\alpha \in W^{i-1}} (e^{i}((\sigma^\smallfrown\alpha)_-, \varphi) +  e^i((\sigma^\smallfrown\alpha)_-, \varphi)_{\mapsfrom{{\vartriangleleft}e^i(c(\sigma^\smallfrown\alpha), {pr}(\sigma^\smallfrown\alpha))}}).\]
Hence, since $\sigma^\smallfrown\alpha \in W^{i-1}$ implies that $(\sigma^\smallfrown\alpha)_- \in W^{i-1}$, we get that 
\[e^{i-1}(\sigma, \Box \varphi) = \bigwedge_{(\sigma^\smallfrown\alpha)_- \in W^{i-1}} (e^{i}((\sigma^\smallfrown\alpha)_-, \varphi) +  e^i((\sigma^\smallfrown\alpha)_-, \varphi)_{\mapsfrom{{\vartriangleleft}e^i(c(\sigma^\smallfrown\alpha), {pr}(\sigma^\smallfrown\alpha))}}).\]

Given that $\sigma \in W^{i-1}\setminus W^i$ (hence $\sigma \neq \sigma_-$), clearly $(\sigma^\smallfrown\alpha)_- = {\sigma_-}^\smallfrown\alpha$, and both $c(\sigma^\smallfrown\alpha) = c(\sigma)$ and ${pr}(\sigma^\smallfrown\alpha) = {pr}(\sigma)$. Moreover, as we checked in Observation \ref{obs:freshSuccessorsOnlyOriginalWorlds} (2), ${\sigma_-}^\smallfrown\alpha \in W^{i-1}$ if and only if ${\sigma_-}^\smallfrown\alpha \in W^{i}$.
Therefore, the above equality implies that 
 \[e^{i-1}(\sigma, \Box \varphi) = \bigwedge_{{\sigma_-}^\smallfrown\alpha \in W^{i}} (e^{i}({\sigma_-}^\smallfrown\alpha, \varphi) +  e^i({\sigma_-}^\smallfrown\alpha, \varphi)_{\mapsfrom{{\vartriangleleft}e^i(c(\sigma), {pr}(\sigma))}}).\]

Recall that $\mod{T}^i$ is quasi-witnessed, so 
$\bigwedge_{{\sigma_-}^\smallfrown\alpha \in W^{i}} (e^{i}({\sigma_-}^\smallfrown\alpha, \varphi)$ equals either $\bot$ or $e^{i}({\sigma_-}^\smallfrown\alpha_0, \varphi)$ for some ${\sigma_-}^\smallfrown\alpha_0 \in W^i$. Therefore, we can apply Lemma \ref{lem:infSupmapsfrom} to get that 
\[e^{i-1}(\sigma, \Box \varphi) = 
(\bigwedge_{{\sigma_-}^\smallfrown\alpha \in W^{i}} (e^{i}({\sigma_-}^\smallfrown\alpha, \varphi)) +  (\bigwedge_{{\sigma_-}^\smallfrown\alpha \in W^{i}} e^i({\sigma_-}^\smallfrown\alpha, \varphi))_{\mapsfrom{{\vartriangleleft}e^i(c(\sigma), {pr}(\sigma))}}).\]
Since by definition $e^{i}(\sigma_-, \Box \varphi) = \bigwedge_{{\sigma_-}^\smallfrown\alpha \in W^{i}} e^{i}({\sigma_-}^\smallfrown\alpha, \varphi)$, this concludes the proof of the $\Box$ case. 

The proof for formulas starting with $\Diamond$ is done analogously. 
\end{proof}

Let us also prove that at each step, the extended model $\mod{T}^{i-1}$ is a conservative extension of $\mod{T}^i$. This will imply that $\mod{T}^+$ is a conservative extension of $\mod{T}$. 
\begin{lemma}\label{lem:conservativeM+}
Let $0 < i \leq md(\Upsilon)$. For any $\sigma \in W^i$ and any formula $\varphi \in \Upsilon_{\depth(\sigma)}$ it holds that

\[e^{i-1}(\sigma, \varphi) = e^i(\sigma, \varphi).\]
In particular, if $\sigma \in W^{md(\Upsilon)}$ (i.e., in $\eta(RW(\mod{T}))$), \[e^{i-1}(\sigma, \varphi) = e(\eta^{-1}(\sigma), \varphi).\]
\end{lemma}
\begin{proof} See Appendix 2.\end{proof}

Together, the last two lemmas allow us to conclude the following result.

\begin{corollary}\label{cor:truthLemmaT+}
For any world $\sigma \in W^+$ such that $\sigma \neq \sigma_-$, and $\varphi \in \Upsilon_{\depth(\sigma)}$, it holds that 
\[e^+(\sigma, \varphi) = e^+(\sigma_-, \varphi) + e^+(\sigma_-, \varphi)_{\mapsfrom {\vartriangleleft}e^+(c(\sigma), {pr}(\sigma))}\]
\end{corollary}
\begin{proof}
Assume $\sigma \in W^{i-1}\setminus W^i$ for some $1 \leq i \leq md(\Upsilon)$. Thus, $\sigma_-$ and $c(\sigma)$ belong to $W^i$, and so we know from Lemma \ref{lem:conservativeM+} that 
$e^+(\sigma, \varphi) = e^{i-1}(\sigma, \varphi)$, $e^+(\sigma_-, \varphi) = e^i(\sigma_-, \varphi)$ and $e^+(c(\sigma), {pr}(\sigma)) = e^i(c(\sigma), {pr}(\sigma))$. 

Lemma \ref{lem:truthLemmaM+} implies that 
$e^{i-1}(\sigma, \varphi) = e^i(\sigma_-, \varphi) + e^i(\sigma_-, p)_{\mapsfrom{{\vartriangleleft}e^i(c(\sigma), {pr}(\sigma))}}$, which concludes the proof. 
 
\end{proof}

All the previous information can be summarized in a result that will allow us to propositionally express a certain regularity condition between two successors of a world $\sigma$ with $\Box \varphi \in uWit_{\mod{T}^+}(\sigma,\Upsilon_{\depth(\sigma)})$. In particular, this result establishes a relation between the evaluations of the formulas in $\Upsilon_{\depth(\sigma)}$ at the worlds $\sigma$ and $\sigma_-$.

\begin{definition}\label{def:avalue}
Let $\sigma \in W^+$ with $\sigma_-\neq \sigma$,  and $\varphi \in \Upsilon_{\depth(\sigma)}$. We let 
\[a_{\sigma, \varphi} \coloneqq e^+(\sigma_-, \varphi)_{\mapsfrom {\vartriangleleft}e^+(c(\sigma), {pr}(\sigma))}.\]
\end{definition}

Roughly speaking, the previous value will be the difference (in terms of the product operation of the algebra) between the value of a formula $\varphi$ in a world $\sigma$ and in its ``seed" $\sigma_-$.

\begin{proposition}\label{prop:valueProperties}
Let $\sigma \in W^+$ with $\sigma \neq \sigma_-$, and $\varphi \in \Upsilon_{\depth(\sigma)}$. Then:
%
\begin{enumerate}
\item $e^+(\sigma, \varphi) = e^+(\sigma_-, \varphi) + a_{\sigma, \varphi}$,
\item If $\varphi, \chi \in \Upsilon_{\depth(\sigma)}$ and $e^+(\sigma_-, \chi) \leq e^+(\sigma_-, \varphi)$ then $a_{\sigma, \chi} \leq a_{\sigma, \varphi}$, 
\item $a_{\sigma, \varphi} = \bot$ if and only if $e^+(\sigma_-, \varphi) = \bot$, 
\item If $\sigma = \delta^\smallfrown\Box \psi'$ for some $\langle \delta, \Box \psi\rangle \in Gens_{\depth(\sigma)}^\Upsilon(\mod{T})$, then

\begin{enumerate}
\item For any $\Box \varphi \in \Upsilon_{\depth(\delta)}$ such that $e^+(\delta, \Box \varphi) > \bot$, then $a_{\sigma, \varphi} = \mathbf{0}$;
\item $\bot < a_{\sigma, \psi} < \mathbf{0}$.
\end{enumerate}
\end{enumerate}

\end{proposition}
\begin{proof}
Let us see that the value from Definition \ref{def:avalue} satisfies the conditions of the proposition.
\begin{enumerate}
\item Is exactly Corollary \ref{cor:truthLemmaT+} (and the definition of $a_{\sigma, \varphi}$).
\item Follows directly from Observation \ref{obs:ordermapsto} (2).
\item Follows from  Observation \ref{obs:ordermapsto} (1).
\item Assume $\sigma = \delta^\smallfrown\Box \psi'$ for some $\langle \delta, \Box \psi\rangle \in Gens_{\depth(\delta)}^\Upsilon(\mod{T})$, and for convenience let us denote $v = \eta^{-1}(\delta)$. By definition, $c(\sigma) = \sigma_- = \delta ^\smallfrown\Box \psi$ and ${pr}(\sigma) = \psi$.   Therefore, $\eta^{-1}(\sigma_-) = v_{\Box \psi}$, and since $\Box \psi \in uWit_{\mod{T}}(v, \Upsilon_{\depth(\delta)})$ we know $v_{\Box \psi}$ satisfies conditions of Lemma \ref{lem:belowWorld}. Using this we can prove the two subcases:
\begin{enumerate}
\item For any $\Box \varphi \in \Upsilon_{\depth(v)}$ such that $e(v, \Box \varphi) > \bot$ we have that 
${\vartriangleleft}e(v_{\Box \psi}, \psi) < {\vartriangleleft}e(v_{\Box \psi}, \varphi)$. Therefore, from Lemmas  \ref{lem:conservativeM+}
 and \ref{lem:inverseEtaPreservation} (recall that $\sigma_-$ and $\delta$ belong to $\eta(RW(\mod{T}))$), and lastly applying Observation \ref{obs:closetriangle} (3)  we get that 
\[ a_{\sigma, \varphi} = e^+(\sigma_-, \varphi)_{\mapsfrom {\vartriangleleft}e^+(c(\sigma), {pr}(\sigma))} =
 e(\eta^{-1}(\sigma_-), \varphi)_{\mapsfrom {\vartriangleleft}e^+(\delta^\smallfrown\Box \psi, \psi)} =
 e(v_{\Box \psi}, \varphi)_{\mapsfrom {\vartriangleleft}e(v_{\Box \psi}, \psi)} = \mathbf{0}. \] 
%
%

\item We first know that $e(v_{\Box \psi}, \psi) > \bot$ because  $\Box \psi \in uWit_{\mod{T}}(v, \Upsilon_{v})$. Using Observation \ref{obs:ordermapsto} (1) and Lemmas  \ref{lem:inverseEtaPreservation} and \ref{lem:conservativeM+}, the previous implies that $a_{\sigma, \psi} > \bot$.

On the other hand, by definition (see also Observation \ref{obs:closetriangle} (1)) and Lemma  \ref{lem:inverseEtaPreservation},
 \[e(v_{\Box \psi}, \psi)_{\mapsfrom {\vartriangleleft}e(v_{\Box \psi}, \psi)}[q] = \begin{cases} e(v_{\Box \psi}, \psi)[{\vartriangleleft}e(v_{\Box \psi}, \psi)] &\hbox{ if }q = {\vartriangleleft}e(v_{\Box \psi}, \psi) \\ 0 &\hbox{ otherwise. }\end{cases}\]

By the choice of $v_{\Box \psi}$, necessarily  $e(v_{\Box \psi}, \psi)[{\vartriangleleft}e(v_{\Box \psi}, \psi)] < 0$. Going back to $\mod{T}^+$ using Lemma \ref{lem:conservativeM+} allows us to conclude that $a_{\sigma, \psi} < \mathbf{0}$. 
\end{enumerate}
\end{enumerate}
\end{proof}

\begin{lemma}\label{lem:awitnesses}
 Let $\sigma \in HRW(\mod{T}^+)$ with $\sigma_- \neq \sigma$, $\Diamond \varphi, \Box \psi \in \Upsilon_{\depth(\sigma)}$ and such that
 $\sigma ^\smallfrown\Box \psi' \not \in W^+$. The following properties hold:
\begin{enumerate}
\item $e^+(\sigma, \Box \psi) = e^+(\sigma^\smallfrown\Box \psi, \psi)$,  $a_{\sigma, \Box \psi} = a_{\sigma^\smallfrown\Box \psi, \varphi}$ and  $a_{\sigma^\smallfrown\Box \psi, \psi} \leq a_{\delta, \psi } $ for any $\delta \in W^+$ with $R^+\sigma\delta$;
\item $e^+(\sigma, \Diamond \varphi) = e^+(\sigma^\smallfrown\Diamond \varphi, \varphi)$, $a_{\sigma, \Diamond \varphi} = a_{\sigma^\smallfrown\Diamond \varphi, \varphi}$ and  $a_{\sigma^\smallfrown\Diamond \varphi, \varphi} \geq a_{\delta, \varphi} $ for any $\delta \in W^+$ with $R^+\sigma\delta$.
\end{enumerate}
\end{lemma}
\begin{proof} See Appendix 2.\end{proof}

We can abstract the labeling of the hereditarily relevant worlds in $\mod{T}^+$ away from the model \textit{per se}: all the possible sequences that might be needed for this (for any arbitrary model) are fully determined by the formulas in $\Upsilon$, as we will see. Stating this in full generality will later allow us to follow an idea similar to those in \cite{CeEs10}, \cite{Vi20}, and others: namely, moving from a modal entailment to a propositional one by splitting each modal formula into a corresponding set of propositional formulas. In this process, we assign to each propositional variable—and to each formula beginning with a modality—a subscript indicating ``the world" at which that formula is being considered. The crucial point is that the number of labels resulting from this process is finite, which makes it possible to reason within propositional logic (which is finitarily complete).

\begin{definition}\label{def:possibleWorlds}
Let $\Upsilon$ be a finite set of formulas, and let $\Upsilon_i^{\Diamondbox} \coloneqq \{\Box \varphi\colon \Box \varphi \in \Upsilon_i\} \cup \{\Diamond \varphi\colon \Diamond \varphi \in \Upsilon_i\}$, for each $0 \leq i \leq md(\Upsilon)$. We define
\[\Sigma^\Upsilon \coloneqq \{\langle 0 \rangle\} \cup \{\langle 0, \beta_1,\ldots, \beta_k \rangle \colon 1 \leq k \leq md(\Upsilon), \beta_i \in \Upsilon_{i-1}^{\Diamondbox} \cup \{\Box \psi' \colon \Box \psi \in \Upsilon_{i-1}\}\}\]
\end{definition}

It is clear that $\vert \Sigma^\Upsilon \vert = \Sigma_{i = 0}^{md(\Upsilon)-1} lv_i$, for 
$lv_0 \coloneqq 1$ (the single world at the root) and 
\[lv_{i+1} \coloneqq lv_i \cdot \vert \{\Diamond \psi \in \Upsilon_i\} \vert + 2 \cdot \vert \{\Box \psi \in \Upsilon_i\} \vert\] (the amount of possible sequences at each modal level). Therefore, $\Sigma^\Upsilon$ is a finite set. 

As before, unless stated otherwise we will assume $\Upsilon$ is a fixed set of formulas and drop the super-index.
It is also clear that now we can produce an equivalent definition of the hereditarily relevant worlds of a model $\mod{T}^+$: it turns out  that $HRW(\mod{T}^+) = W^+ \cap \Sigma$.

\begin{exa}
Continuing Example \ref{example:running1}, the set $\Sigma^\Upsilon$ is given by the following set of sequences:
\begin{align*}
\{\langle 0 \rangle, \langle 0, \Box(y \vee \Box x)\rangle, \langle 0, \Box(y \vee \Box x)'\rangle,  \langle 0, \Diamond \Diamond y\rangle, \\
\langle 0, \Box(y \vee \Box x), \Box x\rangle, \langle 0, \Box(y \vee \Box x), \Box x'\rangle,\\
\langle 0, \Box(y \vee \Box x)', \Box x\rangle, \langle 0, \Box(y \vee \Box x)', \Box x'\rangle,\\
\langle 0, \Box(y \vee \Box x), \Diamond y\rangle, \langle 0, \Box(y \vee \Box x)', \Diamond y\rangle, \\
\langle 0, \Diamond \Diamond y, \Box x\rangle, \langle 0, \Diamond \Diamond y, \Box x'\rangle\}
\end{align*}
 \qed
\end{exa}

Let us define a class of subsets of $\Sigma$ such that the hereditarily relevant worlds of each possible model $\mod{T}^+$ are in one-to-one correspondence with the elements of these subsets. This impacts both the existence (or absence) of primed elements (associated with the possible distribution of unwitnessed worlds) and the very existence of longer sequences (related to the endpoints of the model at depths strictly less than the modal depth of the formulas).

\begin{definition}\label{def:coherentSet}
We say that a set $\Omega \subseteq \Sigma$ is \termDef{coherent} (for $\Upsilon$), and write $\Omega \underset{\text{coh}}{\subseteq} \Sigma$, whenever $\langle 0 \rangle \in \Omega$ and, for any $\sigma^\smallfrown\alpha\in \Omega$ the following conditions hold:
\begin{enumerate}
\item $\sigma \in \Omega$, 
\item $\sigma^\smallfrown\heartsuit \chi \in \Omega$ for any $\heartsuit \chi \in \Upsilon_{\vert \sigma\vert -1}$,  
\item $(\sigma^\smallfrown\alpha)_- \in \Omega$, 
\item For any $\beta$, $\sigma^\smallfrown\beta \in \Omega$ if and only if ${\sigma_-}^\smallfrown\beta \in \Omega$.
\end{enumerate}

\end{definition}

As a way to speak, one can think of coherent sets as those arising from any possible choice regarding which formulas remain unwitnessed in the worlds of a (finite-depth) structure, as well as the depth assigned to each of its branches.
It is immediate that, for $\Omega$ coherent and $\sigma \in \Omega$, we have that $\sigma^\smallfrown\beta \in \Omega$ if and only if $\underline{\sigma}^\smallfrown\beta \in \Omega$ (by iterated application of condition 4).

\begin{lemma}\label{lemma:HRWcoherent}
The hereditarily relevant worlds of $\mod{T}^+$ are a coherent subset of $\Sigma$.
\end{lemma}

\begin{proof}
$\langle 0 \rangle \in W^+$ is by definition a hereditarily relevant world of $\mod{T}^+$. 
Assume $\sigma ^\smallfrown\alpha\in HRW(\mod{T}^+)$, and let us check the conditions of the previous definition:

\begin{enumerate}
\item Clearly, $\sigma \in W^+$ too, hence in $HRW(\mod{T}^+)$,
\item Since $\sigma ^\smallfrown\alpha\in HRW(\mod{T}^+)$ it means $\sigma$ has successors in $\mod{T}^+$. The proof of the case can be done easily by induction on the number of primed elements in $\sigma$. If there is none, we know that $\sigma \in  \tau(RW(\mod{T}))$ and by definition, $\tau^{-1}(\sigma)_{\heartsuit \chi} \in RW(\mod{T})$ for each $\heartsuit \chi \in \Upsilon_{\vert \sigma\vert -1}$. Therefore, $\sigma^\smallfrown \heartsuit \chi \in \tau(\tau^{-1}(\sigma)_{\heartsuit \chi})$ and so, in $W^{md(\Upsilon)} \subseteq W^+$. Otherwise, we know that ${\sigma_-}^\smallfrown \alpha \in HRW(\mod{T}^+)$ by definition. By I.H., ${\sigma_-}^\smallfrown \heartsuit \chi \in HRW(\mod{T}^+)$ for each $\heartsuit \chi \in \Upsilon_{\vert \sigma\vert -1}$. Following Observation \ref{obs:freshSuccessorsOnlyOriginalWorlds} (2), we get that $\sigma ^\smallfrown \heartsuit \chi \in HRW(\mod{T}^+)$ too.
\item It is clear by definition of $\mod{T}^+$ that if $\sigma ^\smallfrown\alpha\in HRW(\mod{T}^+)$, namely, if it is a sequence of only formulas and primed formulas, then the sequence resulting from removing any of the intermediate primes is still a world in $HRW(\mod{T}^+)$.
\item Follows from Observation \ref{obs:freshSuccessorsOnlyOriginalWorlds} (2).\qedhere
\end{enumerate}

\end{proof}

Let us now introduce, for each coherent subset $\Omega$,\footnote{Recall that $\Omega$ depends on the fixed set of formulas $\Upsilon$.} a set of symbols, $\mathcal{V}^\Omega$, which will be used as (propositional) variables in propositional product logic. 
For an arbitrary set of modal formulas $\Theta$, let $pV(\Theta)$ be the propositional variables appearing outside any modality, in the formulas from $\Theta$. 
For $0 \leq i \leq md(\Upsilon)$, consider the set $vars_i \coloneqq pV(\Upsilon_i) \cup \{\heartsuit \psi \colon \heartsuit \psi \in \Upsilon_i\}$. Namely, the propositional variables appearing in $\Upsilon_i$ outside a modality, and the formulas beginning with a modality of the corresponding level. 
Furthermore, we denote by $\Omega_i \coloneqq \{\sigma \in \Omega \colon \vert \sigma \vert = i\}$, the elements in $\Omega$ of length $i$. 

\begin{definition}\label{def:setOfVars}
Let $\Omega$ be a coherent subset of $\Sigma$. 
We define 
\[\mathcal{V}^\Omega \coloneqq \bigcup_{0 \leq i \leq md(\Upsilon)}\{\delta_{\sigma}\colon \delta \in vars_i, \sigma \in \Omega_{i+1}\} \cup \{\alpha_{\sigma, \varphi} \colon \sigma \in \Omega_{i+1}, \sigma \neq \sigma_-, \varphi \in \Upsilon_i\}.\]
\end{definition}

\begin{exa}
For $\Upsilon = \{\Box(y \vee \Box x), \Diamond \Diamond y\}$ as in Example \ref{example:running1}, we have that 
\[vars_0 = \{\Box(y \vee \Box x), \Diamond \Diamond y\},\quad vars_1 = \{y, \Box x, \Diamond y\},\quad vars_2 = \{x,y\}.\]
If we let $\Omega$ to be the set $HRW(\mod{T}^+)$ from Example \ref{example:running2}, we have that 
\begin{align*}
\mathcal{V}^\Omega = \{ \Box(y \vee \Box x)_{\langle 0 \rangle},\ (\Diamond \Diamond y)_{\langle 0 \rangle}\\
y_{\langle 0, \Diamond \Diamond y\rangle},\ y_{\langle 0, \Box(y \vee \Box x) \rangle },\ y_{\langle 0, \Box(y \vee \Box x)' \rangle }\\
(\Box x)_{\langle 0, \Diamond \Diamond y\rangle},\ (\Box x)_{\langle 0, \Box(y \vee \Box x) \rangle },\ (\Box x)_{\langle 0, \Box(y \vee \Box x)' \rangle }\\
(\Diamond y)_{\langle 0, \Diamond \Diamond y\rangle},\ (\Diamond y)_{\langle 0, \Box(y \vee \Box x) \rangle },\ (\Diamond y)_{\langle 0, \Box(y \vee \Box x)' \rangle }\\
x_{\langle 0, \Diamond \Diamond y, \Diamond y\rangle},\ x_{\langle 0, \Diamond \Diamond y, \Box x\rangle}, \\
x_{\langle 0, \Box(y \vee \Box x), \Diamond y \rangle },\ x_{\langle 0, \Box(y \vee \Box x), \Box x \rangle },\ x_{\langle 0, \Box(y \vee \Box x), \Box x' \rangle },\\
x_{\langle 0, \Box(y \vee \Box x)', \Diamond y \rangle },\ x_{\langle 0, \Box(y \vee \Box x)', \Box x \rangle },\ x_{\langle 0, \Box(y \vee \Box x)', \Box x' \rangle }, \\
y_{\langle 0, \Diamond \Diamond y, \Diamond y\rangle},\ y_{\langle 0, \Diamond \Diamond y, \Box x\rangle},\ y_{\langle 0, \Diamond \Diamond y, \Box x'\rangle}, \\
y_{\langle 0, \Box(y \vee \Box x), \Diamond y \rangle },\ y_{\langle 0, \Box(y \vee \Box x), \Box x \rangle },  \\
y_{\langle 0, \Box(y \vee \Box x)', \Diamond y \rangle },\ y_{\langle 0, \Box(y \vee \Box x)', \Box x \rangle },\ y_{\langle 0, \Box(y \vee \Box x)', \Box x' \rangle }, \\
\alpha_{\langle 0, \Box(y \vee \Box x)'\rangle, y},\ \alpha_{\langle 0, \Box(y \vee \Box x)'\rangle, \Box x},\ \alpha_{\langle 0, \Box(y \vee \Box x)'\rangle,  \Diamond y},\\
\alpha_{\langle 0, \Box(y \vee \Box x)', \Box x\rangle, x},\ 
\alpha_{\langle 0, \Box(y \vee \Box x)', \Box x\rangle, y }, \\
\alpha_{\langle 0, \Box(y \vee \Box x)', \Box x'\rangle, x },\ 
\alpha_{\langle 0, \Box(y \vee \Box x)', \Box x'\rangle, y }, \\
\alpha_{\langle 0, \Box(y \vee \Box x)', \Diamond y\rangle, x},\ 
\alpha_{\langle 0, \Box(y \vee \Box x)', \Diamond y\rangle, y}\}
\end{align*}
While a bit cumbersome, the above is simply the list of all possible ``values" that fully determine the behavior of the model: namely all the pairs $\langle$ formula-world$\rangle$ whose values in the model determine the value of $\Upsilon$ at the root, plus variables $\alpha$ that will encode the difference between the value of the formula at a world $\sigma$ and at is ``seed" $\sigma_-$.  \qed
\end{exa}

%
%
%

Observe that, according to this definition, for each propositional variable $p \in pV(\Upsilon_i)$ we have a new variable $p_\sigma$ for each $\sigma \in \Omega$ with $\vert\sigma\vert = i+1$. In the same spirit, for each formula $\heartsuit \psi \in \Upsilon_i$ beginning with a modality, we also have a fresh variable $(\heartsuit \psi)_\sigma$ for the corresponding $\sigma$. Any propositional combination of elements with the same sequence as subindex will be denoted accordingly. Namely, $(\varphi \star \psi)_\sigma \coloneqq \varphi_\sigma \star \psi_\sigma$ for any propositional connective $\star$, and $(\neg \varphi)_\sigma \coloneqq \neg \varphi_\sigma$.
Henceforth, for any (modal) formula $\varphi \in \Upsilon_i$ and $\sigma \in \Omega$ with $\vert \sigma\vert = i+1$, we denote by $\varphi_{\sigma}$ the corresponding propositional formula over the language $\mathcal{V}^{\Omega}$, moving the subindex $\sigma$ inwards through the propositional connectives, until the level of formulas beginning with a modality (that turn into the new $(\heartsuit \chi)_\sigma$ propositional variables) and to propositional variables.

Fixed $\Omega$ a coherent subset of $\Sigma$, let us now define the following sets of propositional formulas in the language of $\Pi_\Delta$ over the variables $\mathcal{V}^\Omega$:
\[\mathcal{EP}^{\Omega}, \mathcal{W}it_\Diamond^{\Omega}, \mathcal{W}it_\Box^{\Omega}, u\mathcal{W}it^{\Omega}, \mathcal{T}op_\alpha^\Omega, 2\mathcal{V}^{\Omega},\mathcal{O}rd_\alpha \text{ and } \mathcal{N}eg^{\Omega}.\]
In order to help with the reading, after each of the definitions we include in italics its intuitive meaning, for a model satisfying it. In those explanations, whenever we talk about a formula taking a value in a certain ``world" $\sigma$, we will mean a formula whose value at that world is relevant for the evaluation of $\Upsilon$ at $\langle 0 \rangle$ (namely, formulas belonging to $\Upsilon_{\vert \sigma \vert -1}$). 
\begin{enumerate}[\label={}]
\item 
\begin{fleqn}\begin{align*}
\bullet\ \mathcal{EP}^{\Omega} \coloneqq  \bigcup_{\sigma \in \Omega, \not \exists  \sigma^\smallfrown \beta \in \Omega }\ &
 \{ (\Box \psi)_\sigma \colon (\Box \psi)_\sigma \in \mathcal{V}^\Omega\} \cup \{\alpha_{\sigma, \Box \psi} \colon \alpha_{\sigma, \Box \psi} \in \mathcal{V}^\Omega\} \cup\\ 
 &\{ \neg (\Diamond \psi)_\sigma \colon (\Diamond \psi)_\sigma \in \mathcal{V}^\Omega\} \cup \{\neg \alpha_{\sigma, \Diamond \psi} \colon \alpha_{\sigma, \Diamond \psi} \in \mathcal{V}^\Omega \}.
 \end{align*}
\end{fleqn}


\textit{If $\sigma$ is understood as end-point of $\Omega$ -i.e., no $\sigma^\smallfrown \beta$ belongs to $\Omega$-, the modal formulas (and their respective $\alpha$ values) are evaluated to the corresponding value ($\top$ for $\Box$ formulas and $\bot$ for $\Diamond$ ones).}

In our running example, taking $\Omega$ as the set $HRW(\mod{T}^+)$, the previous set would be empty (there are no modal formulas evaluated at any ending point, because all ending points arising from $\Omega$ are of depth 2, as $\Upsilon$).
%
\item 
\begin{fleqn}\[\bullet\  \mathcal{W}it_\Diamond^{\Omega} \coloneqq \{W_\Diamond((\Diamond \psi)_\sigma) \colon \sigma^\smallfrown \Diamond \psi \in \Omega \}\text{, with } \]\end{fleqn}
\begin{align*}
W_\Diamond((\Diamond \psi)_\sigma) \coloneqq &
((\Diamond \psi)_\sigma \leftrightarrow \psi_{\sigma^\smallfrown\Diamond \psi}) \wedge  (\underset{\sigma^\smallfrown \beta \in \Omega}{\bigvee} \psi_{\sigma^\smallfrown \beta} \rightarrow (\Diamond \psi)_\sigma)
\end{align*}
if $\sigma = \sigma_-$ (i.e., if $\sigma$ does not have primed elements), and 
\begin{align*}
W_\Diamond((\Diamond \psi)_\sigma) \coloneqq &
((\Diamond \psi)_\sigma \leftrightarrow \psi_{\sigma^\smallfrown\Diamond \psi}) \wedge  (\underset{\sigma^\smallfrown \beta \in \Omega}{\bigvee} \psi_{\sigma^\smallfrown \beta} \rightarrow (\Diamond \psi)_\sigma) \wedge\\
&(\alpha_{\sigma, \Diamond \psi} \leftrightarrow  \alpha_{\sigma^\smallfrown\Diamond \psi, \psi}) \wedge (\bigvee_{\sigma^\smallfrown \beta \in \Omega} \alpha_{\sigma^\smallfrown\beta, \psi} \rightarrow  \alpha_{\sigma^\smallfrown\Diamond \psi, \psi}) 
\end{align*} otherwise.

\textit{If $\sigma$ is not an end-point, the formulas beginning with a $\Diamond$ are assigned the corresponding maximum value among the successors of $\sigma$, which also coincides with the value of that formula in the corresponding witnessing world. The analogous is done for the associated $\alpha$ value, whenever it exists, following the behaviour identified in Lemma \ref{lem:awitnesses}.}

In our running example, a formula appearing (among others) in the previous set is 
\begin{align*}
((\Diamond \Diamond y)_{\langle 0\rangle} \leftrightarrow (\Diamond  y)_{\langle 0, \Diamond \Diamond y\rangle}) \wedge (((\Diamond y)_{\langle 0, \Box(y \vee \Box x)\rangle}\vee (\Diamond y)_{\langle 0, \Box(y \vee \Box x)'  \rangle}) \rightarrow (\Diamond \Diamond y)_{\langle 0\rangle}).
\end{align*}

\item 
\begin{fleqn}\[
\bullet\ \mathcal{W}it_\Box^{\Omega} \coloneqq \{W_\Box((\Box \psi)_\sigma) \colon  \sigma^\smallfrown \Box \psi  \in \Omega, \sigma^\smallfrown \Box \psi'  \not \in \Omega \} \text{, with }\]\end{fleqn}
\begin{align*}
W_\Box((\Box \psi)_\sigma) \coloneqq &
((\Box \psi)_\sigma \leftrightarrow \psi_{\sigma^\smallfrown\Box \psi}) \wedge  ((\Box  \psi)_\sigma \rightarrow \underset{\sigma^\smallfrown \beta \in \Omega}{\bigwedge} \psi_{\sigma^\smallfrown \beta})
\end{align*}
if $\sigma = \sigma_-$ (i.e., if $\sigma$ does not have primed elements), and 
\begin{align*}
W_\Box((\Box \psi)_\sigma) \coloneqq &
((\Box \psi)_\sigma \leftrightarrow \psi_{\sigma^\smallfrown\Box \psi}) \wedge  ((\Box  \psi)_\sigma \rightarrow \underset{\sigma^\smallfrown \beta \in \Omega}{\bigwedge} \psi_{\sigma^\smallfrown \beta}) \wedge\\
&(\alpha_{\sigma, \Box  \psi} \leftrightarrow  \alpha_{\sigma^\smallfrown\Box  \psi, \psi}) \wedge (\alpha_{\sigma^\smallfrown\Box  \psi, \psi} \rightarrow 
\bigwedge_{\sigma^\smallfrown \beta \in \Omega} \alpha_{\sigma^\smallfrown\beta, \psi}) 
\end{align*} otherwise.

\textit{If $\sigma$ is not an end-point, the formulas $\Box \psi$ that ``are witnessed in $\sigma$ in $\Omega$" (namely, $\sigma^\smallfrown \Box \psi'\not \in \Omega$) are assigned the corresponding minimum value among the successors of $\sigma$, which also coincides with the value of that formula in the corresponding witness. The analogous is done for the associated $\alpha$ value, whenever it exists, following Lemma \ref{lem:awitnesses}}.

In our running example, the previous set is given by the singleton
\[
\{((\Box x)_{\langle 0, \Diamond \Diamond y\rangle} \leftrightarrow x_{\langle 0, \Diamond \Diamond y, \Box x \rangle}) \wedge 
((\Box x)_{\langle 0, \Diamond \Diamond y\rangle} \rightarrow x_{\langle 0, \Diamond \Diamond y, \Diamond y\rangle})\}.
\]

\item 
\begin{fleqn}\[
\bullet\ \mathcal{G}ens^{\Omega}_\Box \coloneqq \{\neg (\Box \psi)_\sigma \wedge \neg \Delta \alpha_{\sigma^\smallfrown \Box \psi', \psi} \colon \sigma^\smallfrown \Box \psi' \in \Omega,  \sigma = \sigma_-\}.\]\end{fleqn} 


 \textit{For any pair $\langle \sigma, \Box \psi\rangle$ read as a source of unwitnessing (in the sense of $Gens_i$ in model $\mod{T}$), the value of that formula is assigned to $\bot$ and the $\alpha$-value of $\psi$ in the world $\sigma^\smallfrown \Box \psi'$ is forced to be below $\top$. This follows point (4).b) of Proposition \ref{prop:valueProperties}.}

In our running example, the previous set is given by
\begin{align*}
\{\neg ((\Box (y \vee \Box x))_{\langle 0\rangle} \wedge \neg \Delta \alpha_{\langle 0,\Box(y \vee \Box x)' \rangle, (y \vee \Box x)},\\
\neg (\Box x)_{\langle 0,\Box(y \vee \Box x) \rangle} \wedge \neg \Delta \alpha_{\langle 0,\Box(y \vee \Box x), \Box x' \rangle, x}\}.
\end{align*}

\item 

\begin{fleqn}\[
\bullet\ \mathcal{T}op_\alpha^\Omega \coloneqq \{\neg \neg (\Box \psi)_\sigma \rightarrow \alpha_{\sigma^\smallfrown \Box \chi', \psi}\colon \sigma^\smallfrown \Box \chi' \in \Omega, \sigma = \sigma_-, (\Box \psi)_\sigma \in \mathcal{V}^\Omega, \psi \neq \chi \}.\]\end{fleqn}

\textit{For formulas $\Box\psi$ not evaluated to $\bot$ in a world $\sigma$ which is source of unwitnessing as in the previous point, the associated value of $\alpha$ to $\psi$ in $\sigma^\smallfrown \Box \chi'$ (for any such sequence belonging to $\Omega$) is $\top$, following point (4).a) of Proposition \ref{prop:valueProperties}.  }

In our running example, the previous set is empty (there are no other $\Box$-formulas to be evaluated in the selected worlds).

\item 
\begin{fleqn}\[
\bullet\ 2\mathcal{V}^{\Omega}_\alpha \coloneqq \{\psi_\sigma \leftrightarrow (\psi_{\sigma_-} \with \alpha_{\sigma, \psi}) \colon \sigma \in \Omega, \sigma \neq \sigma_-, \psi \in \Upsilon_{\vert \sigma \vert -1}\}.\]\end{fleqn}


\footnote{Here and in the following formulas we are exploiting the first and third conditions of the definition of coherent set to use that also $\sigma_- \in \Omega$.} 
\textit{For any $\sigma \in \Omega$ with at least some primed element, the value of any formula in $\sigma$ is assigned to be in relation to that of $\sigma_-$ through the corresponding $\alpha$-value, following point (1) of Proposition \ref{prop:valueProperties}.}

In our running example, some formulas included in the previous set are the following ones:
\begin{align*}
y_{\langle 0, \Box(y \vee \Box x), \Box x'\rangle} \leftrightarrow&\  (y_{\langle 0, \Box(y \vee \Box x), \Box x\rangle} \with \alpha_{\langle 0, \Box(y \vee \Box x), \Box x'\rangle, y}), \\
y_{\langle 0, \Box(y \vee \Box x)', \Box x\rangle} \leftrightarrow&\  (y_{\langle 0, \Box(y \vee \Box x), \Box x\rangle} \with \alpha_{\langle 0, \Box(y \vee \Box x)', \Box x\rangle, y}), \\
y_{\langle 0, \Box(y \vee \Box x)', \Box x'\rangle} \leftrightarrow&\  (y_{\langle 0, \Box(y \vee \Box x), \Box x'\rangle} \with \alpha_{\langle 0, \Box(y \vee \Box x)', \Box x'\rangle, y}). \\
\end{align*}

\item 
%

\begin{fleqn}\[
\bullet\ \mathcal{O}rd^{\Omega}_\alpha \coloneqq \{\Delta(\psi_{\sigma_-} \rightarrow \chi_{\sigma_-}) \rightarrow (\alpha_{\sigma, \psi} \rightarrow \alpha_{\sigma, \chi}) \colon \sigma \in \Omega, \sigma \neq \sigma_-,  \psi,\chi \in \Upsilon_{\vert \sigma \vert -1}\}.\]\end{fleqn}
 

\textit{For any $\sigma \in \Omega$ with at least some primed element, and any pair of formulas, the corresponding $\alpha$-values are ordered according to the point (2) of Proposition \ref{prop:valueProperties}.}

In our running example, some formulas included in the previous set are the following ones:
\begin{align*}
\Delta(y_{\langle 0, \Box(y \vee \Box x)\rangle} \rightarrow (\Box x)_{\langle 0, \Box(y \vee \Box x)\rangle}) \rightarrow (\alpha_{\langle 0, \Box(y \vee \Box x)'\rangle, y} \rightarrow \alpha_{\langle 0, \Box(y \vee \Box x)'\rangle, \Box x}),\\
\Delta(y_{\langle 0, \Box(y \vee \Box x), \Box x'\rangle} \rightarrow x_{\langle 0, \Box(y \vee \Box x), \Box x'\rangle}) \rightarrow (\alpha_{\langle 0, \Box(y \vee \Box x)', \Box x'\rangle, y} \rightarrow \alpha_{\langle 0, \Box(y \vee \Box x)', \Box x'\rangle,  x}).\\
\end{align*}

\item 
%

\begin{fleqn}\[
\bullet\ \mathcal{N}eg^{\Omega} \coloneqq \{\neg \alpha_{\sigma,\psi} \leftrightarrow \neg \psi_{\sigma_-} \colon \alpha_{\sigma,\psi} \in \mathcal{V}^\Omega\}.\]\end{fleqn}



\textit{For any $\alpha_{\sigma,\psi} \in \mathcal{V}^\Omega$(namely, $\sigma \in \Omega$ with at least some primed element, and $\psi$ in the corresponding depth of $\Upsilon$) its value is assigned $\bot$ if and only if the corresponding formula in $\sigma_-$ is also $\bot$ (namely, following point (3) of Proposition \ref{prop:valueProperties}).}

In our running example, some formulas included in the previous set are the following ones:
\begin{align*}
\neg \alpha_{\langle 0, \Box(y \vee \Box x)'\rangle, y} \leftrightarrow \neg y_{\langle 0, \Box(y \vee \Box x)\rangle},\quad
\neg \alpha_{\langle 0, \Box(y \vee \Box x)'\rangle, \Box x} \leftrightarrow \neg (\Box x)_{\langle 0, \Box(y \vee \Box x)\rangle},\\
\neg \alpha_{\langle 0, \Box(y \vee \Box x)', \Box x\rangle, x} \leftrightarrow \neg x_{\langle 0, \Box(y \vee \Box x), \Box x},\quad
\neg \alpha_{\langle 0, \Box(y \vee \Box x)', \Box x'\rangle, x } \leftrightarrow \neg x_{\langle 0, \Box(y \vee \Box x), \Box x'}.
\end{align*}

\end{enumerate}

\begin{definition}\label{def:definingFormulas}
Let $\Upsilon$ be a finite set of modal formulas. We let 

\[\Upsilon^\star \coloneqq \bigvee_{\Omega \underset{\text{coh}}{\subseteq} \Sigma}\bigwedge \Upsilon^\Omega, \]
for \[\Upsilon^\Omega \coloneqq \mathcal{EP}^{\Omega} \cup \mathcal{W}it_\Diamond^{\Omega} \cup \mathcal{W}it_\Box^{\Omega} \cup u\mathcal{W}it^{\Omega} \cup \mathcal{T}op_\alpha^\Omega \cup 2\mathcal{V}^{\Omega} \cup \mathcal{O}rd_\alpha \cup \mathcal{N}eg^{\Omega}.\]

\end{definition}

Observe that the definition of $\Upsilon^\star$ is constructive and depends only on the set of formulas $\Upsilon$. It is completely independent of the initial model $\mod{T}$, which was obtained non-constructively. This distinction is crucial for proving the decidability of $\vdash_{K\Pi}$. Indeed, if $\mod{T}$ were used in the definition of $\Upsilon^\star$ (for example, by defining $\Upsilon^\star$ directly as $\Upsilon^{HRW(\mod{T}^+)}$, which would give a shorter encoding) then the resulting set would not necessarily be computable and therefore could not be used in a decidability proof.

It is worth noting, however, that this simpler alternative definition would still suffice to establish the standard completeness of the logic, which in any case follows from the general approach we are adopting.

As we anticipated in the statement of Theorem \ref{th:all}, formula $\Upsilon^\star$ can now provide us a way to move from modal entailments to propositional ones. The next result formalizes this idea and proves the direction \emph{(2)} $\Rightarrow$ \emph{(1)} of the main Theorem \ref{th:all}.

\begin{theorem}\label{prop:ModToProp}
Let $\Upsilon= \Gamma \cup \{\varphi\}$ be a finite set of formulas such that $\Gamma \not \Vdash_{K\Pi} \varphi$. Then 
\[\Gamma_{\langle 0 \rangle}, \Upsilon^\star
\not \models_{\Pi_\Delta} \varphi_{\langle 0 \rangle}.\]
\end{theorem}
\begin{proof}
From Corollary \ref{cor:quasiWitComp} we know that if $\Gamma \not \Vdash_{K\Pi} \varphi$ then there is a quasi-witnessed $\mathbf{R}$-Kripke tree $\mod{T}$ with root $r$ such that $\Gamma \not \Vdash_{\mod{T},r} \varphi$. From Lemma \ref{lem:conservativeM+}, we also know that $\Gamma \not \models_{\mod{T}^+, \langle 0 \rangle} \varphi$.

Consider the set $ \Omega = HRW(\mod{T}^+)$, which is a  coherent subset of $\Sigma$, following Lemma  \ref{lemma:HRWcoherent}. Moreover, let $h$ be the propositional $\Pi_\Delta$-homomorphism from $\mathcal{V}_p(\Upsilon^\star)$ into the $\Delta$-product chain $\mathbf{R}_\Delta$ defined by:
\begin{align*}
h(p_\sigma) \coloneqq& e^+(\sigma, p) \text{ for }p \in \mathcal{V}(\Upsilon), \sigma \in \Omega,\\
h((\heartsuit \psi)_\sigma) \coloneqq& e^+(\sigma, \heartsuit \psi) \text{ for }\heartsuit \in \{\Box, \Diamond\}, \sigma \in \Omega,\\
h(\alpha_{\sigma,\phi}) \coloneqq& a_{\sigma, \phi} \text{ (from Definition \ref{def:avalue}).}
\end{align*}

Since the modal evaluations are, world-wise, propositional homomorphisms, it is immediate that 
for each $\sigma \in \Omega$ and each $\psi \in \Upsilon_{\vert \sigma \vert-1}$ it holds that $h(\psi_\sigma) = e^+(\sigma, \psi) = e(\tau^{-1}(\sigma), \psi)$. In particular, $h(\gamma_{\langle 0 \rangle}) = e(r,\gamma) =\top$ for each $\gamma \in \Gamma$ and $h(\varphi_{\langle 0 \rangle}) =e(r,\varphi) < \top$. Let us prove below that the formula $\Upsilon^\star$ is evaluated to $\top$, by checking that, for $\Omega$ as above, every formula in $\Upsilon^\Omega$ is evaluated to $\top$ under $h$. This will conclude the proof of the theorem.

\subsection*{The case of $\mathcal{EP}^{\Omega}$} Recall that 
\begin{align*}
\mathcal{EP}^{\Omega} \coloneqq  \bigcup_{\sigma \in \Omega, \not \exists  \sigma^\smallfrown \beta \in \Omega }\ &
 \{ (\Box \psi)_\sigma \colon (\Box \psi)_\sigma \in \mathcal{V}^\Omega\} \cup \{\alpha_{\sigma, \Box \psi} \colon \alpha_{\sigma, \Box \psi} \in \mathcal{V}^\Omega\} \cup\\ 
 &\{ \neg (\Diamond \psi)_\sigma \colon (\Diamond \psi)_\sigma \in \mathcal{V}^\Omega\} \cup \{\neg \alpha_{\sigma, \Diamond \psi} \colon \alpha_{\sigma, \Diamond \psi} \in \mathcal{V}^\Omega\}.
 \end{align*}
 
 For each $\sigma \in \Omega$ such that there is not any  $\sigma^\smallfrown \beta \in \Omega$, by definition of $\Omega$ (i.e., $HRW(\mod{T}^+)$) there are no successors of $\sigma$ in $HRW(\mod{T}^+)$. Therefore, either there are no successors of $\sigma$ in $\mod{T}^+$ at all, or there are no formulas in $\Upsilon_{\vert \sigma \vert -1}$. In the second case, there is no $(\Box \psi)_\sigma$, or $(\Diamond \psi)_\sigma$ in $\mathcal{V}^\Omega$. In the first case, indeed $\sigma$ is an endpoint of the tree $\mod{T}^+$ so indeed $e^+(\sigma, \Box \chi) = \top$ and $e^+(\sigma, \Diamond \chi) = \bot$ for any $\chi$, as expected. From Proposition \ref{prop:valueProperties} (1) we get that $h(\alpha_{\sigma, \Box \psi})= \top$ (since $\top = e^+(\sigma, \Box \psi) = e^+(\sigma_-, \Box \psi) + a_{\sigma, \Box \psi}$, necessarily both elements are equal to $\top$) , and from (3) that $h(\alpha_{\sigma, \Diamond \psi})= \bot$, whenever they exist.
 
 \subsection*{The case of $\mathcal{W}it_\Diamond^{\Omega}$} 
  Recall that 
\[\mathcal{W}it_\Diamond^{\Omega} \coloneqq \{W_\Diamond((\Diamond \psi)_\sigma) \colon \sigma^\smallfrown \Diamond \psi \in \Omega \}\text{, with } \]
\begin{align*}
W_\Diamond((\Diamond \psi)_\sigma) \coloneqq &
((\Diamond \psi)_\sigma \leftrightarrow \psi_{\sigma^\smallfrown\Diamond \psi}) \wedge  (\underset{\sigma^\smallfrown \beta \in \Omega}{\bigvee} \psi_{\sigma^\smallfrown \beta} \rightarrow (\Diamond \psi)_\sigma),
\end{align*}
if $\sigma = \sigma_-$ (i.e., if $\sigma$ does not have primed elements), and otherwise
\begin{align*}
W_\Diamond((\Diamond \psi)_\sigma) \coloneqq &
((\Diamond \psi)_\sigma \leftrightarrow \psi_{\sigma^\smallfrown\Diamond \psi}) \wedge  (\underset{\sigma^\smallfrown \beta \in \Omega}{\bigvee} \psi_{\sigma^\smallfrown \beta} \rightarrow (\Diamond \psi)_\sigma) \wedge\\
&(\alpha_{\sigma, \Diamond \psi} \leftrightarrow  \alpha_{\sigma^\smallfrown\Diamond \psi, \psi}) \wedge (\bigvee_{\sigma^\smallfrown \beta \in \Omega} \alpha_{\sigma^\smallfrown\beta, \psi} \rightarrow  \alpha_{\sigma^\smallfrown\Diamond \psi, \psi}).
\end{align*}   
Follows directly from Lemma \ref{lem:awitnesses} (2).

 \subsection*{The case of $\mathcal{W}it_\Box^{\Omega}$} 
The definitions are analogous to those of the previous case, and the proof again follows from Lemma \ref{lem:awitnesses} (1).
 
 \subsection*{The case of $\mathcal{G}ens^{\Omega}_\Box$}
  Recall that 
\[
  \mathcal{G}ens^{\Omega}_\Box \coloneqq \{\neg (\Box \psi)_\sigma \wedge \neg \Delta \alpha_{\sigma^\smallfrown \Box \psi', \psi} \colon \sigma^\smallfrown \Box \psi' \in \Omega,  \sigma = \sigma_-\}.\]

If $\sigma^\smallfrown \Box \psi' \in \Omega$ and $ \sigma = \sigma_-$, necessarily $\langle \sigma, \Box \psi \rangle \in Gens_{\vert \sigma \vert -1}$ and so $\Box \psi$ is unwitnessed in $\sigma$ in $\mod{T}^+$. Therefore, $e^+(\sigma, \Box \psi) = \bot$, and $a_{\sigma^\smallfrown \Box \psi', \psi} < \mathbf{0}$ (so $h(\neg \Delta \alpha_{\sigma^\smallfrown \Box \psi', \psi}) = \neg \Delta a_{\sigma^\smallfrown \Box \psi', \psi}= \top$), from Proposition \ref{prop:valueProperties}, (4), (b).

  \subsection*{The case of $\mathcal{T}op_\alpha^\Omega$}
    Recall that 
\[\mathcal{T}op_\alpha^\Omega \coloneqq \{\neg \neg (\Box \psi)_\sigma \rightarrow \alpha_{\sigma^\smallfrown \Box \chi', \psi}\colon \sigma^\smallfrown \Box \chi' \in \Omega, \sigma = \sigma_-, (\Box \psi)_\sigma \in \mathcal{V}^\Omega \}.\]
  The proof is done as that of the previous case, relying in Proposition \ref{prop:valueProperties}, (4), (a). 

  \subsection*{The case of $2\mathcal{V}^{\Omega}_\alpha$}
  
    Recall that 
\[2\mathcal{V}^{\Omega}_\alpha \coloneqq \{\psi_\sigma \leftrightarrow (\psi_{\sigma_-} \with \alpha_{\sigma, \psi}) \colon \sigma \in \Omega, \sigma \neq \sigma_-, \psi \in \Upsilon_{\vert \sigma \vert -1}\}.\]
Follows from Proposition \ref{prop:valueProperties}, (1). 

  \subsection*{The case of $2\mathcal{O}rd^{\Omega}_\alpha$}
  
    Recall that 
\[\mathcal{O}rd^{\Omega}_\alpha \coloneqq \{\Delta(\psi_{\sigma_-} \rightarrow \chi_{\sigma_-}) \rightarrow (\alpha_{\sigma, \psi} \rightarrow \alpha_{\sigma, \chi}) \colon \sigma \in \Omega, \sigma \neq \sigma_-,  \psi,\chi \in \Upsilon_{\vert \sigma \vert -1}\}.\]

Follows from Proposition \ref{prop:valueProperties}, (2).

  \subsection*{The case of $\mathcal{N}eg^{\Omega}$}
  
    Recall that 
\[\mathcal{N}eg^{\Omega} \coloneqq \{\neg \alpha_{\sigma,\psi} \leftrightarrow \neg \psi_{\sigma_-} \colon \sigma \in \Omega, \sigma \neq \sigma_-,  \psi \in \Upsilon_{\vert \sigma\vert - 1}\}.\]

Follows from Proposition \ref{prop:valueProperties}, (3). Indeed, $h(\psi_{\sigma_-}) = e^+(\sigma_-, \psi)$ and, from Proposition \ref{prop:valueProperties}, (3), this value equals $\bot$ if and only if also $a_{\sigma, \psi} = \bot$.
\end{proof}

\section{From Propositional product logic to $K[0,1]_\Pi$}\label{sec:completeness}

In this section we will see how to go from the propositional condition characterized in the last theorem of the previous section to a (possibly infinite) standard Kripke model that is a counter model for the original entailment.

%

Let us begin by pointing out certain simple conditions on the behavior of the variables $\alpha$ under homomorphisms into $[0,1]_{\Pi_\Delta}$ that satisfy$\Upsilon^\Omega$, for some $\Omega$ coherent subset of $\Sigma$ (following the definitions of $\Sigma$ -Def.\ref{def:possibleWorlds}, of coherent set -Def. \ref{def:coherentSet} and of $\Upsilon^\star$ -Def.\ref{def:definingFormulas} from the previous section) for a fixed set of modal formulas $\Upsilon$.

Before doing so, let us define for convenience, given any $\sigma \in \Sigma$,

\begin{align*}
init(\sigma) \coloneqq \begin{cases} \emptyset &\hbox{ if }\sigma = \sigma_-, \\ init(\sigma_-) \cup \{\sigma\} &\hbox{ otherwise. } \end{cases}
\end{align*}

The next observation follows from the iterated application of the third condition of Definition \ref{def:coherentSet}.
\begin{obs}\label{obs:initPreservation}
For $\Omega$ coherent subset of $\Sigma$, if $\sigma \in \Omega$ then $init(\sigma) \subseteq \Omega$. 
\end{obs}

The following are immediate from the use of the formulas in the sets $2\mathcal{V}^\Omega$ and in $\mathcal{O}rd_\alpha^\Omega$.

\begin{obs}\label{obs:XXXXXX}
Let $\Omega \underset{coh}{\subseteq}\Sigma$, 
 $\sigma \in \Omega$, and $h$ an homomorphism 
 into $[0,1]_{\Pi_\Delta}$ such that $h(\Upsilon^\Omega) = 1$. Then, for any $\sigma \in \Omega$, and formulas $\varphi_\sigma, \alpha_{\sigma, \varphi_1 \with \varphi_2}$ and $\alpha_{\sigma, \psi_1 \rightarrow \psi_2} \in \mathcal{V}^\Omega$ the following hold:
 \setcounter{equation}{0}
\begin{align}
h(\varphi_\sigma) &= h(\varphi_{\underline{\sigma}}) \cdot \prod_{\eta \in init(\sigma)}h(\alpha_{\eta, \varphi}),\\
h(\alpha_{\sigma, \varphi_1 \with \varphi_2}) &=  h(\alpha_{\sigma,\varphi_1}) \cdot h(\alpha_{\sigma,\varphi_2})\qquad\text{ and }\qquad
h(\alpha_{\sigma, \psi_1 \rightarrow \psi_2}) = h(\alpha_{\sigma, \psi_1}) \rightarrow h(\alpha_{\sigma, \psi_2}).
\end{align}
\end{obs}
\begin{proof}
For the first observation, recall that from formulas in $2\mathcal{V}^\Omega_\alpha$, we have that 
$h(\varphi_\sigma) = h(\varphi_{\sigma_-}) \cdot h(\alpha_{\sigma, \varphi})$, and applying $2\mathcal{V}^\Omega_\alpha$ again, we get 
$h(\varphi_\sigma) = h(\varphi_{(\sigma_-)_-}) \cdot h(\alpha_{\sigma_-, \varphi}) \cdot h(\alpha_{\sigma, \varphi})$. Iterating this procedure until reaching $\underline{\sigma}$ in the first element results in the claimed equality. For the second, recall that if $\alpha_{\sigma, \varphi_1 \with \varphi_2} \in \mathcal{V}^\Omega$, by definition also $\alpha_{\sigma, \varphi_1}, \alpha_{\sigma, \varphi_2 }\in \mathcal{V}^\Omega$, since $\varphi_1 \with \varphi_2 \in \Upsilon_{\vert \sigma \vert -1}$ implies that $\varphi_1, \varphi_2 \in \Upsilon_{\vert \sigma \vert -1}$ too -and the same holds for the implication operation. The claim follows again from $2\mathcal{V}^\Omega_\alpha$.
\end{proof}

The above properties are stated for one single homomorphim $h$, but we can gain inspiration in them to, roughly speaking, iterate the first behavior relying in the conditions from the second point along a family of homomorphisms. In order to do so, we will adapt those notions to a relation between two homomorphisms, for which we first  introduce a simple idea concerning deformation of homomorphisms.

\begin{definition}\label{def:orderPresDeformation}
Let $\Theta$ be a closed set of (propositional) formulas. Given two homomorphisms $f,g$ into  $[0,1]_{\Pi_\Delta}$, we say that \termDef{$g$ is a $\Theta$ order-preserving deformation of $f$} whenever for each $\varphi, \psi \in \Theta$, 
\[g(\varphi) = 0 \Longleftrightarrow f(\varphi) = 0\qquad \text{ and }\qquad f(\varphi) \leq f(\psi) \Longrightarrow g(\varphi) \leq g(\psi).\] 
%
\end{definition}

While the first condition could equivalently be stated referring to variables only, the second one cannot be derived from its analogous version over the involved propositional variables, which is a strictly weaker condition. Clearly, any homomorphism is an order-preserving deformation of itself, for any set of formulas.

For two homomorphisms $f,g$ into $[0,1]_{\Pi_\Delta}$, we denote by $fg$ the product of the two (and by $f^n$ the product of $f$ by itself $n$ times), i.e., the homomorphism defined by
\[fg(p) \coloneqq f(p) \cdot g(p).\] 

An interesting feature is that order-preserving deformations of the same homomorphism inherit the structure of the definition, in the following sense.

\begin{proposition}\label{prop:truthLemmaOP}
Let $f,g,h$ be homomorphisms into  $[0,1]_{\Pi_\Delta}$ such that $f$ and $g$ are $\Theta$ order-preserving deformations of $h$.For any $\psi \in \Theta$ it holds that 
\[fg(\psi) = f(\psi) \cdot g(\psi).\]
\end{proposition}
\begin{proof} See Appendix 3.\end{proof}

The following observations are immediate, since the product in $[0,1]$ is order-preserving and $a_1 \cdot \ldots \cdot a_n = 0$ if and only if $a_i = 0$ for some $1\leq i \leq n$, for any finite family $a_1,\ldots, a_n \in [0,1]$.

\begin{obs}\label{obs:ordPres}
Let $f$ be an homomorphism into $[0,1]_{\Pi_\Delta}$. 
\begin{enumerate}
\item If $\{f_i\colon i \in I\}$ is a finite family of $\Theta$ order-preserving deformations of $f$, then $\Pi_{i \in I}f_i$ is also a $\Theta$ order-preserving deformation of $f$.
\item If $g$ is a $\Theta$ order-preserving deformation of $f$, then for any $n \in \omega$, $g^n$ is also a $\Theta$ order-preserving deformation of $f$.
\item If $\{f_i\colon i \in I\}$ is a finite family of $\Theta$ order-preserving deformations of $f$, and $mult\colon I \rightarrow \omega$ is a multiplicity function for $I$, $\Pi_{i \in I}f_i^{mult(i)}$ is also a $\Theta$ order-preserving deformation of $f$.
\end{enumerate}
\end{obs}

Exploiting the previous notions we can construct a (possibly infinite) standard Kripke model from a propositional homomorphism arising from Proposition \ref{prop:ModToProp}. To do so in the simplest possible way, rather than working with a single homomorphism that evaluates the full set of variables $\mathcal{V}^\Omega$ for some $\Omega$ coherent subset of $\Sigma$, we will work with two families of homomorphisms, indexed by the elements in $\Omega$. One family will address  the values taken by the $\alpha$ variables (reading these values as ``offsets" or deformation degrees), and the other evaluate the formulas without any $\alpha_s$ (i.e., the original formulas). 
%
%
The set of variables which is now important for each of the homomorphisms in the two families also varies depending on the element of the coherent set $\Omega$ that is indexing them. Specifically, for $\sigma \in \Omega$, we let 
\[\mathcal{V}_\sigma \coloneqq \mathcal{V}ars(\Upsilon_{\vert \sigma \vert -1}) \cup \{\heartsuit \psi\colon \heartsuit \psi \in  \Upsilon_{\vert \sigma \vert -1}\}.\]
Namely, as in the previous section, we treat each formula in $\Upsilon_{\vert \sigma \vert -1}$ that begins with a modality as a fresh propositional variable —removing its subindex in the process. In fact, for each $\sigma \in \Omega$, we have:

\[\{p_\sigma\colon p \in \mathcal{V}_\sigma\}  = \{p_\sigma \colon p_\sigma  \in \mathcal{V}^\Omega\}.\]

In practice, each of the homomorphisms we define below will evaluate the corresponding set of formulas $\Upsilon_{\vert \sigma \vert -1}$, now viewed as propositional formulas over the extended variable set  $\mathcal{V}_\sigma$. 

\begin{definition}
Let $\Omega \underset{coh}{\subseteq}\Sigma$. Given $\sigma \in \Omega$, and a homomorphism $h$ from $\mathcal{V}^\Omega$ into $[0,1]_{\Pi_\Delta}$, we define the pair of homomorphisms $\varrho_{h,\sigma}$ and $ \delta_{h, \sigma}$ 
from $\mathcal{V}_\sigma$ into $[0,1]_{\Pi_\Delta}$ defined by
\[\varrho_{h,\sigma}(p) \coloneqq h(p_\sigma) \qquad \text{ and } \qquad \delta_{h, \sigma}(p) \coloneqq h(\alpha_{\sigma,p}).\]
\end{definition}

If $h$ is clear from the context, we will simply write $\varrho_\sigma$ and $\delta_\sigma$.

Recall that the set of variables $\mathcal{V}_\sigma$ includes all formulas beginning by a modality in $\Upsilon_{\vert \sigma \vert -1}$. Since $\varrho_{h,\sigma}$ and $\delta_{h, \sigma}$ are propositional homomorphisms, it follows that
for each $\varphi \in \Upsilon_{\vert \sigma \vert -1}$ (seen as set of propositional formulas over the set of variables $\mathcal{V}_\sigma$), it holds that
\[\varrho_{h,\sigma}(\varphi) = h(\varphi_\sigma) \qquad \text{ and } \qquad \delta_{h, \sigma}(\varphi) = h(\alpha_{\sigma,\varphi}).\]
The first equality is immediate by the definition of $\varphi_\sigma$, while the second one follows from Observation \ref{obs:XXXXXX} (2). 

%

\begin{lemma}\label{lem:PropToHoms}
Let $\Omega \underset{coh}{\subseteq}\Sigma$, and $h$ a homomorphism from $\mathcal{V}^\Omega$ into $[0,1]_{\Pi_\Delta}$ such that $h(\Upsilon^\Omega) = \{1\}$. Then, for each $\sigma \in \Omega$, and each
 $\eta \in init(\sigma)$, 
 $\delta_{h, \eta}$ is a $\Upsilon_{\vert \sigma \vert -1}$ order-preserving deformation of $\varrho_{h,\underline{\sigma}}$. 

\end{lemma}
\begin{proof}
Consider any $\varphi \in  \Upsilon_{\vert \sigma \vert -1}$. Clearly $\vert \sigma \vert = \vert \eta\vert = \vert\underline{\sigma}\vert$ for any $\eta \in init(\sigma)$. Hence, since $h(\mathcal{N}eg^\Omega) = h(2\mathcal{V}^\Omega) = \{1\}$, it follows that $h(\alpha_{\eta,\varphi}) = 0$ if and only if $h(\varphi_{\underline{\sigma}}) = 0$
 for each $\eta \in init(\sigma)$. By definition, $\delta_{h, \eta}(\varphi) = h(\alpha_{\eta,\varphi})$ and $\varrho_{h,\underline{\sigma}}(\varphi) = h(\varphi_{\underline{\sigma}})$. This proves the first condition of the definition of an order-preserving deformation (Def. \ref{def:orderPresDeformation}).

On the other hand, let $\varphi,\psi \in  \Upsilon_{\vert \sigma \vert -1}$ be such that $h(\varphi_{\underline{\sigma}}) \leq h(\psi_{\underline{\sigma}})$. Iterating that $h(\mathcal{O}rd^\Omega) = h(2\mathcal{V}^\Omega) = \{1\}$, and since the product in $[0,1]$ is order preserving, we get that, for any $\eta \in init(\sigma)$, also $h(\varphi_{\eta_-}) \leq h(\psi_{\eta_-})$. Applying $\mathcal{O}rd^\Omega$ a last time, it follows that $h(\alpha_{\eta, \varphi}) \leq h(\alpha_{\eta, \psi})$. Again, by the definition of $\delta_{h, \eta}$ and of $\varrho_{h,\underline{\sigma}}$, we conclude the second condition of Definition \ref{def:orderPresDeformation}.
\end{proof}

To use a finite coherent set $\Omega$ to produce the universe of a possibly infinite standard Kripke model, we are going to turn the elements of $\Omega$ that have some primed element into infinitely many replicas, that will behave, in a certain sense, as powers of the basic one. 

We let
%
%

\begin{align*}
pows(\sigma) \coloneqq& \{\sigma\}, && \text{ for } \sigma \in \Omega \text{ with } \sigma = \sigma_-,\\
pows(\sigma_1^\smallfrown\Box \varphi'^\smallfrown\sigma_2) \coloneqq&  \{\varrho^\smallfrown(\Box \varphi)_k^\smallfrown\sigma_2\colon k \in \omega, \varrho \in pows(\sigma_1)\}, && \text{ for }\sigma_1^\smallfrown\Box \varphi'^\smallfrown\sigma_2 \in \Omega, \sigma_2 \text{ prime-free.}
\end{align*}
Now, given $\eta = \eta_1^\smallfrown(\Box \varphi)_k^\smallfrown\eta_2 \in pows(\sigma)$ for some $\sigma \in \Omega$ and for $\eta_1$ subindexes-free, and similarly to the notation used for elements in $\Sigma$, we  let\footnote{For $\eta = \eta_-$, $init(\eta) = \emptyset$.}:
\begin{align*}
\eta_- \coloneqq  \eta_1^\smallfrown \Box \varphi^\smallfrown\eta_2, &  & \underline{\eta} \coloneqq \underline{\sigma}, & & \muwave{\eta} \coloneqq \sigma,\\
init(\eta) \coloneqq init(\eta_-) \cup \{\eta\}, & & mult(\eta) \coloneqq k. &&
\end{align*}

\begin{obs}\label{obs:onPowProperties}
If $\eta \in init(\nu)$ then $\muwave{\eta} \in init(\muwave{\nu})$. Moreover,
$\underline{\nu} = \underline{\muwave{\nu}}$.

\end{obs}

\begin{definition}\label{def:canonicalModel}
Let $\Omega$ be a coherent subset of $\Sigma$. We define the crisp frame $\mathfrak{F}_\Omega = \langle W_\Omega, R_\Omega\rangle $ by
\begin{align*}
W_\Omega \coloneqq&  \bigcup_{\sigma \in \Omega} pows(\sigma),\\
R_\Omega \coloneqq& \{\langle \nu, \nu^\smallfrown \beta\rangle\colon \nu^\smallfrown \beta \in W_\Omega\}\ (\beta \text{ may be a formula beginning by a modality or a subindexed } \Box\text{-formula)}.
\end{align*}
Moreover, let $h$ be any homomorphism  into $[0,1]_{\Pi_\Delta}$ such that $h(\Upsilon^\Omega) \subseteq\{1\}$. We let $\mod{T}_{\langle h, \Omega\rangle} \coloneqq \langle \mathfrak{F}_\Omega, e_h\rangle$ be the 
$[0,1]_\Pi$ crisp Kripke model given by extending the frame above with the evaluation

\[e_h(\nu, p) \coloneqq \varrho_{h,\underline{\nu}}(p) \cdot \prod\limits_{\eta \in init(\nu)}\delta_{h,\muwave{\eta}}(p)^{mult(\eta)}, \text{ for } p \in \mathcal{V}ars(\Upsilon).\]
\end{definition}

Namely, according to the notation on the product of homomorphisms we introduced before, $e_h(\nu, p) = (\varrho_{h,\underline{\nu}} \prod\limits_{\eta \in init(\nu)}\delta_{h,\muwave{\eta}}^{mult(\eta)}) (p)$. An example of a frame $\mod{F}_\Omega$ can be seen in Example \ref{example:4}.

\begin{exa}\label{example:4}
Let us use Examples \ref{example:running1} and \ref{example:running2} as a source for a coherent set. The following is the frame $\mathfrak{F}_\Omega$ for $\Omega = HRW(\mod{T}^+)$ from Example \ref{example:running2}.

\[\begin{tikzcd}[sep=tiny,scale cd=0.78]
	&&&& {\langle 0, \Diamond \Diamond y, \Box x\rangle} \\
	&& {\langle 0, \Diamond \Diamond y \rangle} && {\langle 0, \Diamond \Diamond y, \Diamond y\rangle} \\
	&&&& {\langle 0, \Box (y \vee \Box x), \Diamond y\rangle} \\
	&&&& {\langle 0, \Box (y \vee \Box x), \Box x\rangle} \\
	&& {\langle 0, \Box (y \vee \Box x)\rangle} && {\langle 0, \Box (y \vee \Box x), (\Box x)_1\rangle} \\
	&&&&& {\langle 0, \Box (y \vee \Box x), (\Box x)_k\rangle} \\
	\\
	{\langle 0 \rangle} &&&& {\langle 0, \Box (y \vee \Box x)_1, \Diamond y\rangle} && {} \\
	&&&& {\langle 0, \Box (y \vee \Box x)_1, \Box x\rangle} \\
	&& {\langle 0, \Box (y \vee \Box x)_1\rangle} && {\langle 0, \Box (y \vee \Box x)_1, (\Box x)_1\rangle} \\
	\\
	&&&&& {\langle 0, \Box (y \vee \Box x)_1, (\Box x)_k\rangle} \\
	\\
	&&&& {\langle 0, \Box (y \vee \Box x)_n, \Diamond y\rangle} && {} &&& {} \\
	&&&& {\langle 0, \Box (y \vee \Box x)_n, \Box x\rangle} \\
	&& {\langle 0, \Box (y \vee \Box x)_n\rangle} && {\langle 0, \Box (y \vee \Box x)_n, (\Box x)_1\rangle} \\
	\\
	&&&&& {\langle 0, \Box (y \vee \Box x)_n, (\Box x)_k\rangle} \\
	&& {} \\
	&&&&&& {}
	\arrow[curve={height=-24pt}, from=8-1, to=2-3]
	\arrow[curve={height=-6pt}, from=8-1, to=5-3]
	\arrow[from=2-3, to=1-5]
	\arrow[from=2-3, to=2-5]
	\arrow[from=5-3, to=5-5]
	\arrow[from=5-3, to=4-5]
	\arrow[from=5-3, to=3-5]
	\arrow[curve={height=-6pt}, from=8-1, to=10-3]
	\arrow[from=10-3, to=9-5]
	\arrow[from=10-3, to=10-5]
	\arrow[curve={height=12pt}, from=10-3, to=12-6]
	\arrow[""{name=0, anchor=center, inner sep=0}, curve={height=12pt}, Rightarrow, dotted, no head, from=12-6, to=14-7]
	\arrow[""{name=1, anchor=center, inner sep=0}, dotted, no head, from=10-3, to=16-3]
	\arrow[""{name=2, anchor=center, inner sep=0}, curve={height=12pt}, dotted, no head, from=5-5, to=6-6]
	\arrow[""{name=3, anchor=center, inner sep=0}, Rightarrow, dotted, no head, from=16-3, to=19-3]
	\arrow[""{name=4, anchor=center, inner sep=0}, curve={height=12pt}, dotted, no head, from=10-5, to=12-6]
	\arrow[from=16-3, to=15-5]
	\arrow[from=16-3, to=16-5]
	\arrow[""{name=5, anchor=center, inner sep=0}, curve={height=12pt}, Rightarrow, dotted, no head, from=18-6, to=20-7]
	\arrow[""{name=6, anchor=center, inner sep=0}, curve={height=12pt}, Rightarrow, dotted, no head, from=6-6, to=8-7]
	\arrow[from=10-3, to=8-5]
	\arrow[curve={height=18pt}, from=5-3, to=6-6]
	\arrow[from=16-3, to=14-5]
	\arrow[curve={height=18pt}, from=16-3, to=18-6]
	\arrow[""{name=7, anchor=center, inner sep=0}, curve={height=6pt}, dotted, no head, from=16-5, to=18-6]
	\arrow[curve={height=18pt}, from=8-1, to=16-3]
	\arrow[shorten >=12pt, dotted, from=8-1, to=1]
	\arrow[curve={height=30pt}, Rightarrow, dotted, from=8-1, to=3]
	\arrow[curve={height=18pt}, shorten >=22pt, Rightarrow, dotted, from=10-3, to=0]
	\arrow[curve={height=30pt}, shorten >=23pt, Rightarrow, dotted, from=16-3, to=5]
	\arrow[curve={height=12pt}, shorten >=13pt, dotted, from=5-3, to=2]
	\arrow[curve={height=18pt}, shorten >=22pt, Rightarrow, dotted, from=5-3, to=6]
	\arrow[curve={height=6pt}, dotted, from=10-3, to=4]
	\arrow[curve={height=6pt}, shorten >=13pt, dotted, from=16-3, to=7]
\end{tikzcd}\]
 \qed
\end{exa}

\begin{lemma}\label{lem:truthLemma}
Let $\Omega \underset{coh}{\subseteq}\Sigma$ and $h$ be a homomorphism from $\mathcal{V}^\Omega$ into $[0,1]_{\Pi_\Delta}$ such that $h(\Upsilon^\Omega) \subseteq\{1\}$. For any $\nu \in W_\Omega$ and $\varphi \in \Upsilon_{\vert \nu \vert -1}$,\footnote{Recall that, following to the definition of the (propositional) homomorphisms $h_{\underline{\nu}}$ and $\delta^h_\muwave{\eta}$, in the right side of the equation below the formula $\varphi$ is understood as a purely propositional formula over the extended set of variables $\mathcal{V}_{\muwave{\nu}}$.}
\[
e_h(\nu, \varphi) = \varrho_{h,\underline{\nu}}(\varphi) \cdot \prod\limits_{\eta \in init(\nu)}\delta_{h,\muwave{\eta}}(\varphi)^{mult(\eta)}.
\]
\end{lemma}
\begin{proof}
To lighten the notation, and since $h$ is fixed thorough the lemma, we will drop it from the subindexes above. Namely, in the proof below we will speak about $\mod{T}_\Omega$, and we will write (and prove that)
\begin{equation}\label{eq:truthLemmaeqToProve}
e(\nu, \varphi) = \varrho_{\underline{\nu}}(\varphi) \cdot \prod\limits_{\eta \in init(\nu)}\delta_{\muwave{\eta}}(\varphi)^{mult(\eta)}.
\end{equation}

We will prove the lemma by induction on the complexity of the formula. 
Before proceeding, recall that, as 
we proved in Lemma \ref{lem:PropToHoms}, $\delta_{\muwave{\eta}}$ is a 
$\Upsilon_{\vert \nu \vert -1}$ order-preserving deformation of $\varrho_{\underline{\nu}}$, for each $\eta \in init(\nu)$ (seeing $\Upsilon_{\vert \nu \vert -1}$
 as a set of propositional formulas over the variables $\mathcal{V}_{\muwave{\nu}}$). We then know that $\prod_{\eta \in init(\nu)}\delta_{\muwave{\eta}}^{mult(\eta)}$ is also a $\Upsilon_{\vert \nu \vert -1}$ order-preserving deformation of $\varrho_{\underline{\nu}}$, following Observation \ref{obs:ordPres}. Therefore, from Proposition \ref{prop:truthLemmaOP}, we get that for any formula $\psi$ in $\Upsilon_{\vert \nu \vert -1}$ (seen as a propositional formula over the variables $\mathcal{V}_{\muwave{\nu}}$) it holds that  
 \begin{equation}\label{eq:prodInherited}
(\varrho_{\underline{\nu}}\prod_{\eta \in init(\nu)}\delta_{\muwave{\eta}}^{mult(\eta)})(\psi) = \varrho_{\underline{\nu}}(\psi) \cdot\prod_{\eta \in init(\nu)} \delta_{\muwave{\eta}}(\psi)^{mult(\eta)}.
 \end{equation}
 We will also resort, without further notice, to the fact that $\eta \in init(\nu)$ implies that $\muwave{\eta} \in init(\muwave{\nu})$, and that 
$\underline{\nu} = \underline{\muwave{\nu}}$ (Observation \ref{obs:onPowProperties}).

We begin by considering the case where the outermost connective of $\varphi$ is propositional.
Assume that $\varphi = \varphi_1 \star \varphi_2 \in \Upsilon_{\vert \nu \vert -1}$ with $\star$ propositional connective (the case for unary connective $\Delta$ is proven analogously). By I.H.

\[e(\nu, \varphi_1 \star \varphi_2) =  e(\nu, \varphi_1) \star e(\nu,  \varphi_2) = (\varrho_{\underline{\nu}}(\varphi_1) \cdot \prod_{\eta \in init(\nu)}\delta_{\muwave{\eta}}(\varphi_1)^{mult(\eta)}) \star (\varrho_{\underline{\nu}}(\varphi_2) \cdot \prod_{\eta \in init(\nu)}\delta_{\muwave{\eta}}(\varphi_2)^{mult(\eta)}).\]
Following equation (\ref{eq:prodInherited}) (which is applicable to $\varphi_1, \varphi_2$ since, by definition, $\Upsilon_{\vert \nu \vert -1}$ is closed under propositional subformulas) we get that 

\[e(\nu, \varphi_1 \star \varphi_2) =  (\varrho_{\underline{\nu}} \prod_{\eta \in init(\nu)}\delta_{\muwave{\eta}}^{mult(\eta)})(\varphi_1) \star (\varrho_{\underline{\nu}} \prod_{\eta \in init(\nu)}\delta_{\muwave{\eta}}^{mult(\eta)})(\varphi_2).\]
Since $\varrho_{\underline{\nu}} \cdot \prod_{\eta \in init(\nu)}\delta_{\muwave{\eta}}^{mult(\eta)}$ is a propositional homomorphism it preserves propositional operations. Therefore, using equation (\ref{eq:prodInherited}) again (since $\varphi = \varphi_1 \star \varphi_2 \in \Upsilon_{\vert \nu \vert -1})$, it follows that 

\[e(\nu, \varphi_1 \star \varphi_2) = (\varrho_{\underline{\nu}} \prod_{\eta \in init(\nu)}\delta_{\muwave{\eta}}^{mult(\eta)})(\varphi) = \varrho_{\underline{\nu}}(\varphi) \cdot \prod_{\eta \in init(\nu)}\delta_{\muwave{\eta}}(\varphi)^{mult(\eta)}.\]


 For formulas beginning by a modality, consider first the case when there is no $\muwave{\nu}^\smallfrown\alpha \in \Omega$. It is clear that there is no 
 $\underline{\nu}^\smallfrown\alpha \in \Omega$, nor any $\eta^\smallfrown\alpha \in \Omega $ for $\eta \in init(\nu)$, following the Definition of coherent set (\ref{def:coherentSet}, see also observation below it). 
  Therefore, the worlds $\nu, \underline{\nu}$ and $\eta$ for $\eta \in init(\nu)$ don't have successors in the model $\mod{T}_\Omega$.  
Henceforth, by definition, \[e(\nu, \Box \psi) = 1 \qquad \text{ and } \qquad e(\nu, \Diamond \psi) = 0, \qquad \text{ for any }\psi.\]
On the other hand, since $h$ evaluates to $1$ the set $\mathcal{EP}^\Omega$, we get that for any $\Box \psi, \Diamond \psi \in \Upsilon_{\vert \nu \vert -1}$,
\begin{align*}
\varrho_{\muwave{\nu}}(\Box \psi) &=  h((\Box \psi)_{\muwave{\nu}}) = \delta_{\muwave{\nu}}(\Box \psi) = h(\alpha_{\muwave{\nu}, \Box \psi}) = 1, \text{ and }\\
\varrho_{\muwave{\nu}}(\Diamond \psi) &= h((\Diamond \psi)_{\muwave{\nu}}) = \delta_{\muwave{\nu}}(\Diamond \psi) = h(\alpha_{\muwave{\nu}, \Diamond \psi}) = 0.
\end{align*}

This immediately concludes that $e(\nu, \Diamond \psi) = \varrho_{\muwave{\nu}}(\Diamond \psi) = \varrho_{\muwave{\nu}}(\Diamond \psi) \cdot \prod_{\eta \in init(\nu)}\delta_{\muwave{\eta}}(\Diamond \psi)^{mult(\eta)} = 0$.
 
For the $\Box$ case, also
  $h((\Box \psi)_{\muwave{\eta}}) = 1$ for any $\eta \in init(\nu)$. Since $h$ evaluates to $1$ the formulas in $2\mathcal{V}^\Omega$, we get that 
$\delta_{\muwave{\eta}}(\Box \psi)^{mult(\eta)} = 1$. It follows that 
$e(\nu, \Box \psi) =  \varrho_{\muwave{\nu}}(\Box \psi) \cdot \prod_{\eta \in init(\nu)}\delta_{\muwave{\eta}}(\Box \psi)^{mult(\eta)} = 1$, proving Equation (\ref{eq:truthLemmaeqToProve}) for the case of ``ending points" according to $\Omega$.

We will consider now the case when there is some $\muwave{\nu}^\smallfrown\alpha \in \Omega$.

 Let us begin with the more complex case of formulas beginning with a $\Box$ operation. By definition of the Kripke model $\mod{T}_\Omega$, 
\begin{align*}
e(\nu, \Box \varphi) &= \bigwedge \{e(\nu^\smallfrown\heartsuit \chi , \varphi) \colon \nu^\smallfrown\heartsuit \chi \in W_\Omega\} \wedge \bigwedge \{e(\nu^\smallfrown(\Box \chi)_k, \varphi)\colon \nu^\smallfrown(\Box \chi)_k \in W_\Omega\} =\\
&  
\bigwedge \{e(\nu^\smallfrown\heartsuit \chi , \varphi)\colon \nu^\smallfrown\heartsuit \chi \in W_\Omega\} \wedge \bigwedge \{e(\nu^\smallfrown(\Box \chi)_k, \varphi)\colon k \in \omega, \muwave{\nu}^\smallfrown\Box \chi' \in \Omega\}.
\end{align*}
By Induction Hypothesis, the previous equals to $A \wedge B$ for 
\begin{align*}
A \coloneqq& \bigwedge \{(\varrho_{\underline{\nu^\smallfrown\heartsuit \chi}}(\varphi) \cdot \prod_{\eta \in init(\nu^\smallfrown\heartsuit\chi)}\delta_{\muwave{\eta}}(\varphi)^{mult(\eta)})\colon \nu^\smallfrown\heartsuit \chi \in W_\Omega\}, \text{ and }\\
B \coloneqq& \bigwedge \{(\varrho_{\underline{\nu^\smallfrown(\Box \chi)_k}} (\varphi) \cdot \prod_{\eta \in init(\nu^\smallfrown(\Box \chi)_k)}\delta_{\muwave{\eta}}(\varphi)^{mult(\eta)}) \colon k \in \omega, \muwave{\nu}^\smallfrown\Box \chi' \in \Omega\}.
\end{align*}

We split the proof in two cases.
\begin{enumerate}
\item Assume $\muwave{\nu}^\smallfrown\Box \varphi' \in \Omega$. Therefore, $\underline{\nu}^\smallfrown\Box \varphi' \in \Omega$ too, and from  $h(\mathcal{G}ens_\Box^\Omega) = 1$ we get that, for any $k$, $\delta_{\muwave{\underline{\nu}^\smallfrown(\Box \varphi)_k}}(\varphi) = \delta_{\underline{\nu}^\smallfrown \Box \varphi'}(\varphi) = h(\alpha_{\underline{\nu}^\smallfrown \Box \varphi', \varphi}) < 1$. Since $\underline{\nu}^\smallfrown \Box \varphi_k \in init(\nu^\smallfrown \Box \varphi_k)$, in particular
\[B \leq \bigwedge_{\underline{\nu}^\smallfrown(\Box \varphi)_k\colon k \in \omega} (h(\varphi_{\underline{\nu^\smallfrown(\Box \varphi)_k}}) \cdot \delta_{\muwave{\underline{\nu}^\smallfrown(\Box \varphi)_k}}(\varphi)^k).\] 

This implies that $B = 0$, since it is below all the powers of an element strictly below $1$, and so, 
$e(\nu, \Box \varphi) =0$ too.
Moreover, from $h(\mathcal{G}ens_\Box^\Omega) = 1$ we also get that $\varrho_{\underline{\nu}}(\Box \varphi) = h((\Box \varphi)_{\underline{\nu}}) = 0$, which proves Equation (\ref{eq:truthLemmaeqToProve}) since 
\[e(\nu, \Box \varphi) = 0 = \varrho_{\underline{\nu}}(\Box \varphi) = \varrho_{\underline{\nu}}(\Box \varphi) \cdot \prod_{\eta \in init(\nu)}\delta_{\muwave{\eta}}(\Box \varphi)^{mult(\eta)}.\]

\item Otherwise, assume $\muwave{\nu}^\smallfrown\Box \varphi' \not \in \Omega$ (and hence, since $\Omega$ is coherent, $\underline{\nu}^\smallfrown\Box \varphi' \not \in \Omega$ either, as before). 
From the fact that $h(\mathcal{W}it^\Omega_\Box) = \{1\}$ we can infer the following properties\footnote{We are making use of the definition of coherent set to use that certain elements belong to $\Omega$. In particular, for $C1$, $C2$ and $C4$ we are using the first and second conditions, and for $C3$ we rely in Observation \ref{obs:initPreservation}, which is essentially the third condition of the definition.} :
\begin{align*}
(C1)& \varrho_{\underline{\nu}}(\Box \varphi) = \varrho_{\underline{\nu}^\smallfrown\Box \varphi}(\varphi),\\
(C2)& \varrho_{\underline{\nu}^\smallfrown\Box \varphi}(\varphi) \leq \varrho_{\underline{\nu^\smallfrown\beta}}(\varphi) \text{ for any } \nu^\smallfrown\beta \in W,\\
(C3)& \delta_{\muwave{\nu}}(\Box \varphi) = \delta_{\muwave{\nu}^\smallfrown\Box \varphi}(\varphi),\\
(C4)& \delta_{\muwave{\eta}^\smallfrown \Box \varphi}(\varphi) \leq \delta_{\muwave{\eta ^\smallfrown \beta}}(\varphi) \text{ for any }
\eta \in init(\nu), \eta^\smallfrown \beta \in W.
\end{align*}
In what follows we will resort to these conditions, and to the fact that $\muwave{\nu}^\smallfrown\heartsuit \varphi = \muwave{\nu^\smallfrown\heartsuit \varphi}$, and $\underline{\nu}^\smallfrown\heartsuit \varphi = \underline{\nu^\smallfrown\heartsuit \varphi}$ for any $\nu, \heartsuit \varphi$.

Since $\eta \in  init(\nu^\smallfrown \heartsuit \chi)$ if and only if $\eta = \eta_1^\smallfrown \heartsuit \chi$ for some $\eta_1 \in init(\nu)$, we have that 
\begin{eqnarray*}
\varrho_{\underline{\nu^\smallfrown \Box \varphi}}(\varphi) \leq&  \varrho_{\underline{\nu^\smallfrown \beta}}(\varphi)\text{ (due to (C2)), and }\\
\prod_{\eta \in init(\nu^\smallfrown \Box \varphi) }\delta_{\muwave{\eta}}(\varphi)^{mult(\eta)} \leq& \prod_{\eta \in init(\nu^\smallfrown \heartsuit \chi) }\delta_{\muwave{\eta}}(\varphi)^{mult(\eta)} \text{ (due to (C4)).}
\end{eqnarray*}
%
%

Applying the previous inequalities in the definition of $A$ we get that


\[A = \varrho_{\underline{\nu^\smallfrown \Box \varphi}}(\varphi) \cdot \prod_{\eta \in init(\nu^\smallfrown \Box \varphi) }\delta_{\muwave{\eta}}(\varphi)^{mult(\eta)}.\]

Lastly, using $(C1)$ and $(C3)$ in the previous equation  we get that 
\[A = 
\varrho_{\underline{\nu}}(\Box \varphi) \cdot \prod_{\eta \in init(\nu)}\delta_{\muwave{\eta}}(\Box \varphi)^{mult(\eta)}.\]



Let us now consider two subcases.

Assume that $h((\Box \varphi)_{\muwave{\nu}}) = 0$. From $h(2\mathcal{V}^\Omega_\alpha) = \{1\}$ and from $h(\mathcal{N}eg^\Omega_\alpha) = \{1\}$ we get that 
$h((\Box \varphi)_{\muwave{\nu}}) = h((\Box \varphi)_{\muwave{\nu}_-}) \cdot h(\alpha_{\muwave{\nu}, \Box \varphi})$ and that 
$\neg h(\alpha_{\muwave{\nu}, \Box \varphi}) \leftrightarrow \neg h((\Box \varphi)_{\muwave{\nu}_-})$. Together, the previous imply that
 $h((\Box \varphi)_{\muwave{\nu}}) = 0$ iff $h((\Box \varphi)_{\muwave{\nu}_-})=0$, and iterating this, that $h((\Box \varphi)_{\muwave{\nu}}) = 0$
  iff $h((\Box \varphi)_{\underline{\nu}}) = 0$. Since we assumed that $h(\Box \varphi_{\muwave{\nu}}) = 0$, also $h(\Box \varphi_{\underline{\nu}}) = 0$ and so, by definition, $\varrho_{\underline{\nu}}(\Box \varphi) = 0$ too. Using all the previous in the above equality for $A$, we obtain that $A = 0$ and so, that $e(\nu, \Box \psi ) = 0$. This immediately proves Equation (\ref{eq:truthLemmaeqToProve}) in this subcase, since $\varrho_{\underline{\nu}}(\Box \varphi)  = 0$ implies that also $\varrho_{\underline{\nu}}(\Box \varphi) \cdot \prod\limits_{\eta \in init(\nu)}\delta_{\muwave{\eta}}(\varphi)^{mult(\eta)} = 0$ (the right side of the Equation).

Otherwise, assume that $h((\Box \varphi)_{\muwave{\nu}}) > 0 $. From $h(\mathcal{T}op^\Omega_\alpha) = \{1\}$ we know that $\delta_{\muwave{\nu}^\smallfrown\Box \chi'}( \varphi) = 1$ for any $\muwave{\nu}^\smallfrown\Box \chi' \in \Omega$. Hence, since $\muwave{\nu^\smallfrown(\Box \chi)_k} = \muwave{\nu}^\smallfrown\Box \chi'$ for any $k$, and using the definitions of $init$ and $_{\sim}$, we get that for any $\nu^\smallfrown(\Box \chi)_k \in \Omega$,

\[\prod_{\eta \in init(\nu^\smallfrown(\Box \chi)_k)}\delta_{\muwave{\eta}}(\varphi)^{mult(\eta)} = \delta_{\muwave{\nu}^\smallfrown \Box \chi'}(\varphi)^k \cdot
(\prod_{\eta \in init(\nu^\smallfrown \Box \chi)}\delta_{\muwave{\eta}}(\varphi)^{mult(\eta)}) =  \prod_{\eta \in init(\nu^\smallfrown \Box \chi)}\delta_{\muwave{\eta}}(\varphi)^{mult(\eta)}\] 
and therefore, 

\[B =  \bigwedge_{\muwave{\nu}^\smallfrown\Box \chi' \in \Omega} (\varrho_{\underline{\nu^\smallfrown(\Box \chi)}} (\varphi) \cdot   \prod_{\eta \in init(\nu^\smallfrown \Box \chi)}\delta_{\muwave{\eta}}(\varphi)^{mult(\eta)}).\]
Since $\muwave{\nu}^\smallfrown\Box \chi' \in \Omega$ implies $\muwave{\nu}^\smallfrown\Box \chi \in \Omega$, we get that $B \geq A$. 


This allows us to conclude Equation (\ref{eq:truthLemmaeqToProve}) in this subcase, since
\[e(\nu, \Box \varphi) = A = \varrho_{\underline{\nu}}(\Box \varphi) \cdot \prod_{\eta \in init(\nu)}\delta_{\muwave{\eta}}(\Box \varphi)^{mult(\eta)}.\]

\end{enumerate}

Let us now prove the case of formulas beginning with a $\Diamond$ operation. The idea of the proof is similar to the one for the $\Box$ case, but a bit simpler since the new worlds only diminish the value of the formulas. By definition of the Kripke model $\mod{T}_\Omega$, 
\begin{align*}
e(\nu, \Diamond \varphi) &= \bigvee_{\nu^\smallfrown\heartsuit \chi \in W} e(\nu^\smallfrown\heartsuit \chi , \varphi) \vee \bigvee_{\nu^\smallfrown(\Box \chi)_k \in W} e(\nu^\smallfrown(\Box \chi)_k, \varphi) =\\
&  
\bigvee_{\nu^\smallfrown\heartsuit \chi \in W} e(\nu^\smallfrown\heartsuit \chi , \varphi) \vee \bigvee_{\substack{\nu^\smallfrown(\Box \chi)_k\colon\\ k \in \omega, \muwave{\nu}^\smallfrown\Box \chi' \in \Omega}} e(\nu^\smallfrown(\Box \chi)_k, \varphi).
\end{align*}
By Induction Hypothesis, the previous equals to  $A' \vee B'$ for 
\begin{align*}
A' \coloneqq& \bigvee_{\nu^\smallfrown\heartsuit \chi \in W} (\varrho_{\underline{\nu^\smallfrown\heartsuit \chi}}(\varphi) \cdot \prod_{\eta \in init(\nu^\smallfrown\heartsuit\chi)}\delta_{\muwave{\eta}}(\varphi)^{mult(\eta)}), \\
B' \coloneqq& \bigvee_{\substack{\nu^\smallfrown(\Box \chi)_k\colon\\ k \in \omega, \muwave{\nu}^\smallfrown\Box \chi' \in \Omega}} (\varrho_{\underline{\nu^\smallfrown(\Box \chi)_k}}(\varphi) \cdot \prod_{\eta \in init(\nu^\smallfrown(\Box \chi)_k)}\delta_{\muwave{\eta}}(\varphi)^{mult(\eta)}).
\end{align*}

As in the $\Box$-case, from $h(\mathcal{W}it_\Diamond^\Omega) = \{1\}$ we can obtain some immediate properties:
\begin{align*}
(C1_\Diamond)& h_{\underline{\nu}}(\Diamond \varphi) = h_{\underline{\nu}^\smallfrown\Diamond \varphi}(\varphi) = h_{\underline{\nu^\smallfrown\Diamond \varphi}}(\varphi),\\
(C2_\Diamond)& h_{\underline{\nu}^\smallfrown\Diamond \varphi}(\varphi) \geq h_{\underline{\nu^\smallfrown\beta}}(\varphi) \text{ for any } \nu^\smallfrown\beta \in W,\\
(C3_\Diamond)& \delta_{\muwave{\nu}}(\Diamond \varphi) = \delta_{\muwave{\nu}^\smallfrown\Diamond \varphi}(\varphi),\\
(C4_\Diamond)& \delta_{\muwave{\eta}^\smallfrown \Diamond \varphi} (\varphi) \geq \delta_{\muwave{\eta^\smallfrown \beta}}(\varphi) \text{ for any }\eta \in init(\nu), \eta^\smallfrown\beta \in W.
\end{align*}

As in the $\Box$ case, $\eta \in init(\nu^\smallfrown \Diamond \varphi) $ if and only if $\eta = \eta_1^\smallfrown \Diamond \varphi$ for some $\eta_1 \in init(\nu)$. Hence
\begin{eqnarray*}
\prod_{\eta \in init(\nu^\smallfrown \Diamond \varphi) }\delta_{\muwave{\eta}}(\varphi)^{mult(\eta)} \geq &\prod_{\eta \in init(\nu^\smallfrown \heartsuit \chi) }\delta_{\muwave{\eta}}(\varphi)^{mult(\eta)}, \text{ (due to }(C4_\Diamond)) \text{, and }\\
h_{\underline{\nu^\smallfrown \Diamond \varphi}}(\varphi) \geq & h_{\underline{\nu^\smallfrown \heartsuit \chi}}(\varphi) \text{ (due to }(C2_\Diamond))
\end{eqnarray*}

In combination with $(C1_\Diamond)$ and $(C3_\Diamond)$, this implies that 
%
%
%
%
%

\[A' = \varrho_{\underline{\nu}}(\Diamond \varphi) \cdot \prod_{\eta \in init(\nu)}\delta_{\muwave{\eta}}(\Diamond \varphi)^{mult(\eta)}.\]

On the other hand, it is immediate that $B' \leq A'$, since
\[\prod_{\eta \in init(\nu^\smallfrown(\Box \chi)_k)}\delta_{\muwave{\eta}}(\varphi)^{mult(\eta)} =  \prod_{\eta \in init(\nu^\smallfrown \Box \chi)}\delta_{\muwave{\eta}}(\varphi)^{mult(\eta)} \cdot \delta_{\muwave{\nu^\smallfrown(\Box \chi)_k}}(\varphi)^k \leq 
\prod_{\eta \in init(\nu^\smallfrown \Box \chi)}\delta_{\muwave{\eta}}(\varphi)^{mult(\eta)}.\]

This concludes the proof of Equation (\ref{eq:truthLemmaeqToProve}) for the $\Diamond$-step, since \[e(\nu, \Diamond \varphi) = A' = h_{\underline{\nu}}(\Diamond \varphi) \cdot \prod_{\eta \in init(\nu)}\delta_{\muwave{\eta}}(\Diamond \varphi)^{mult(\eta)}. \qedhere\]
\end{proof}

\begin{proposition}\label{prop:PropToMod}
Let $\Upsilon = \Gamma \cup \{\varphi\} \subseteq_{\omega} Fm$ and $\Omega$ a coherent subset of $\Sigma$ such that 

\[\Gamma_{\langle 0 \rangle}, \Upsilon^\Omega
\not \models_{\Pi_\Delta} \varphi_{\langle 0 \rangle}.\]

Then, $\Gamma \not \Vdash_{K[0,1]_\Pi} \varphi$. 
\end{proposition}

\begin{proof}

From the premises of the proposition (recall that $\models_{\Pi_\Delta}$ is standard complete)  follows that there is a propositional homomorphism $h$ from $\mathcal{V}^\Omega$ into $[0,1]_{\Pi_\Delta}$ such that $h(\Gamma_{\langle 0 \rangle}) \subseteq \{1\}$, $h(\Upsilon^\Omega) \subseteq \{1\}$ and for which $h(\varphi_{\langle 0 \rangle}) < 1$. Consider the model $\mod{T}_{\langle h, \Omega \rangle}$ from Definition \ref{def:canonicalModel}.

It is clear that $init(\langle 0 \rangle) = \emptyset$ and $\underline{\langle 0 \rangle} = \langle 0 \rangle$, and given that $\Gamma \cup \{\varphi\} = \Upsilon = \Upsilon_0 = \Upsilon_{\vert \langle 0\rangle \vert -1}$, from Lemma \ref{lem:truthLemma} it follows that
\begin{align*}
e_h(\langle 0 \rangle, \gamma) &= h_{\underline{\langle 0 \rangle}}(\gamma) \cdot \prod_{\eta \in init(\langle 0 \rangle)}\delta_{\muwave{\eta}}(\gamma)^{mult(\eta)} = h(\gamma_{\langle 0 \rangle}) = 1,\text{ for each } \gamma \in \Gamma, \text{ and } \\
e_h(\langle 0 \rangle, \varphi) &= h_{\underline{\langle 0 \rangle}}(\varphi) \cdot \prod_{\eta \in init(\langle 0 \rangle)}\delta_{\muwave{\eta}}(\varphi)^{mult(\eta)} = h(\varphi_{\langle 0 \rangle}) < 1.
\end{align*}

\end{proof}

\section{On the logics over valued frames: $\Vdash_{M[0,1]_\Pi}$ and $\Vdash_{M\Pi}$}\label{sec:valued}

Logics arising from valued frames were partially studied in \cite{CeEs22} within the framework of description logics. There, it is shown (among other results, more directly related to description logics and not addressed here) that the set of theorems of $\Vdash_{\mathtt{M}[0,1]\Pi}$ is decidable. However, it is not established whether the logic itself (i.e., the consequence relation) is decidable\footnote{For the interested reader, note that Theorem 24 from \cite{CeEs22}, which proves the decidability of the set \textit{Subs}, does not imply the decidability of the consequence relation. For example, $(x, x^2) \not \in Subs$, even though it is clearly a valid entailment. The absence of a Deduction Theorem makes this issue nontrivial.}, nor is the relation between $\Vdash{M\Pi}$ and $\Vdash_{M[0,1]_\Pi}$ addressed.

In this work, we present results concerning modal logics over valued frames. As one of the consequences, we recover the decidability of the set of theorems of $\Vdash_{M[0,1]\Pi}$, a result already established in \cite{CeEs22}. A reasoning similar—though significantly simpler—to the one developed in previous sections for logics with crisp accessibility relations can be applied here to obtain the analogue of Theorem \ref{th:all}.

We want to point out that although our proof differs significantly from the one given in \cite{CeEs22}, some underlying ideas from that paper have inspired our approach (for instance, the formulas in \cite[Definition 10]{CeEs22} play a role similar, though not identical, to those in Definition \ref{def:definingFormulasValued}). 
We provide the complete proof both to extend the decidability results so as to address the standard completeness question and the full logic understood as an entailment, and to present these results in a purely modal-logic setting (whereas the proof in \cite{CeEs22} is framed in fuzzy description logics). Thanks to the notation and ideas developed in the previous sections, our presentation is more concise and streamlined.

In what follows, we establish the next result, analogous to Theorem \ref{th:all}. As in the crisp case, this will immediately imply the standard completeness of $\Vdash_{M\Pi}$ and its decidability (and, in particular, the decidability of its set of theorems).

\begin{theorem}\label{th:allValued}
Let $\Upsilon = \Gamma \cup \{\varphi\} \subseteq_\omega Fm$. There is a propositional formula $\Upsilon^\ast$ that can be effectively constructed from $\Upsilon$ for which the following are equivalent:\footnote{Recall that $\psi_\alpha$ is a propositional formula defined according to Definition \ref{def:modalToPropFormulas}.}

\begin{enumerate}
\item $\Gamma \Vdash_{M\Pi}\varphi$, 
\item $\Gamma_{\langle 0 \rangle}, \Upsilon^\ast \models_{\Pi_\Delta} \varphi_{\langle 0 \rangle}$,
\item $\Gamma \Vdash_{M[0,1]_\Pi}\varphi$.
\end{enumerate}
\end{theorem}

As for the case for crisp frames, the previous theorem immediately implies the following two corollaries.

\begin{corollary}[Decidability of $\Vdash_{M[0,1]_\Pi}$] \label{th:decVal}
For any $\Gamma \cup \{\varphi\}\subseteq_{\omega}Fm$, the problem of determining whether or not 
\[\Gamma \Vdash_{M[0,1]_\Pi} \varphi\]
is decidable. 

In particular, the set of theorems of $\Vdash_{M[0,1]_\Pi}$ is decidable too. 
\end{corollary}
This follows from the equivalence between \emph{(2)} and \emph{(3)} in Theorem \ref{th:allValued}, and the decidability of the logical entailment over $\models_{\Pi_\Delta}$.

\begin{corollary}[Standard Completeness of $\Vdash_{M\Pi}$] \label{th:standardCompVal}
For any $\Gamma \cup \{\varphi\}\subseteq_{\omega}Fm$,
\[\Gamma \Vdash_{M\Pi} \varphi \qquad \text{ if and only if }\qquad \Gamma \Vdash_{M[0,1]_\Pi} \varphi\]
In particular,  their sets of theorems also coincide.
\end{corollary}
This is precisely the equivalence between \emph{(1)} and \emph{(3)} in Theorem \ref{th:allValued}.

From the two previous corollaries, we also get that the logic $\Vdash_{M\Pi}$ is decidable.

The structure of the proof of Theorem \ref{th:allValued} is similar to the one for crisp logics, but several steps can be strongly simplified or directly eliminated. This is due to the fact that moving the value of the accessibility relation in the interval $(0,1)$ allows us to use a single world to capture the information of the unwitnessed situations, in opposition to the crisp case where we needed two worlds.

Theorem 2.9 from \cite{LaMa07} yields, as before, completeness of $\Vdash_{M\Pi}$ with respect to quasi-witnessed valued trees (now with a valued accessibility) over $\mathbf{R}$.

A simplified version of Lemma \ref{lem:belowWorld} which will be enough to work in the new setting is phrased as follows:

\begin{lemma}\label{lem:belowWorldValued}[c.f. Lemma \ref{lem:belowWorld}]
Let $\mod{T}$ be a quasi-witnessed valued $\mathbf{R}$-Kripke tree, $v \in W$ and $\Upsilon$ a finite set of formulas. For every $\Box \varphi \in \Upsilon$ such that $e(v, \Box \varphi) = \bot$ there is some world $v_{\Box \varphi}$ such that both
\[
e(v_{\Box \varphi}, \varphi) < \bigwedge\{ e(v_{\Box \varphi}, \chi) \colon \Box \chi \in \Upsilon, e(v, \Box \chi) > \bot\}\qquad \text{ and }\qquad
e(v_{\Box \varphi}, \varphi) < R(v,v_{\Box \varphi}). 
\]
\end{lemma}
\begin{proof}
If $\Box \varphi$ is witnessed in $v$, the witnessing world clearly satisfies both conditions. Otherwise, let $a = \bigwedge \{e(v_{\Box \varphi}, \chi) \colon \Box \chi \in \Upsilon, e(v, \Box \chi) > \bot\}$, which is strictly greater than $\bot$ since $\Upsilon$ is a finite set. Since we assume  $e(v, \Box \varphi) = \bot$, necessarily there is some $w$ such that $R(v,w) \rightarrow e(w, \varphi) < a$. The first inequality of the lemma follows from the fact that in any product algebra $x \leq y \rightarrow x$ (following from residuation law and integrality). The second one is immediate (otherwise, the previous implication would equal $\mathbf{0}$). 
\end{proof}

Hence, for the previous premises, we can pick a value $r_{v,\Box \varphi} \in \mathbf{R} $ such that \[e(v_{\Box \varphi}, \varphi) < r_{v, \Box \varphi} < R(v,v_{\Box \varphi}) \wedge \bigwedge\{ e(v_{\Box \varphi}, \chi) \colon \Box \chi \in \Upsilon, e(v, \Box \chi) > \bot\}.\]

This provides an alternative definition of relevant worlds (see Definition \ref{def:relevanWorlds}), which we can name $RW_v(\mod{T})$, exactly as the one for the classical case for witnessed formulas and instead using the world $v_{\Box \varphi}$ identified in Lemma \ref{lem:belowWorldValued} for the unwitnessed cases. Also
the mapping $\eta$ (which is defined right after Definition \ref{def:relevanWorlds}) works as before, assigning sequences of formulas to the worlds of the model, while the new definition of the analogous of model $\eta_{\Upsilon}(\mod{T})$ now simply takes into account the valuedness of the accessibility relation. Namely, we define the model $\eta^v_{\Upsilon}(\mod{T}) = \langle \eta^v_\Upsilon(W), \eta^v_\Upsilon(R), \eta^v_\Upsilon(e) \rangle$ by 
\begin{align*}
\eta^v_\Upsilon(W) \coloneqq& \bigcup \{\eta(v) \colon v \in W_i, i \leq md(\Upsilon)\},\\
\eta^v_\Upsilon(R) (\sigma, \delta) \coloneqq& \begin{cases} R(\eta^{-1}(\sigma), \eta^{-1}(\delta)) &\hbox{ if }\delta = \sigma^\smallfrown \alpha \text{ for some symbol }\alpha,\\ 0 &\hbox{ otherwise,}\end{cases} \\
\eta^v_\Upsilon(e)(\sigma, p) \coloneqq& e(\eta^{-1}(\sigma),p) \text{ for all }p.
\end{align*}
For the sake of readability, we will remove the subindex $\Upsilon$ from the previous elements when it is clear from the context. Observe the map $\eta$ is the same as the one from the previous sections, so we will still use $\eta^{-1}$.

The new definition of the extended model (the valued version of Definition \ref{def:mod+}) will be now much simpler. It will not be necessary to modify the values of the formulas, but only to add a copy of the subtree generated by the previously identified world $v_{\Box \varphi}$ (for $\Box \varphi$ unwitnessed in $v$), and slightly modify the value of the accessibility from $v$ to the root of the new subtree. Observe that, as we pointed out in the introduction, this value is crucially settled in the open interval $(0,1)$.

Let the set $(Gens_i^\Upsilon(\mod{T}))_v \coloneqq \{\langle \sigma, \Box \varphi\rangle \colon \sigma \in \eta(RW_v(\mod{T})), \vert \sigma \vert = i +1, \Box \varphi \in uWit_{\eta(\mod{T})}(\sigma, \Upsilon_i)\}$, analogous to the corresponding set used in Definition \ref{def:mod+} but taking into account the new ``valued" relevant worlds. Again, these are simply the pairs of worlds and of unwitnessed formulas at that world, at corresponding depth of the model and of the set of formulas.


\begin{definition}\label{def:mod+Valued}
Let $\mod{T}$ be a quasi-witnessed valued $\mathbf{R}$-Kripke tree with root $r$, and $\Upsilon$ a finite set of formulas. We define $\mod{T}^0_v, \ldots, \mod{T}_v^{md{\Upsilon}}$ as follows. First, let $\mod{T}^{md(\Upsilon)}_v = \langle W^{md(\Upsilon)}_v, R^{md(\Upsilon)}_v, e^{md(\Upsilon)}_v\rangle \coloneqq \eta^v_\Upsilon(\mod{T})$. Now, for $\mod{T}_v^i$ with $md(\Upsilon) \geq i > 0$ already defined, let
$\mod{T}^{i-1} \coloneqq \langle W^{i-1},R^{i-1},e^{i-1}\rangle$ for \[ W^{i-1}_v \coloneqq W^i_v \cup \bigcup\{W^i_{\langle \sigma, \Box \varphi\rangle}\colon \langle \sigma, \Box \varphi\rangle \in (Gens_i^\Upsilon(\mod{T}))_v\} \text{  with: }\] 

\begin{align*}
 W^i_{\langle \sigma, \Box \varphi\rangle} \coloneqq&  \{\delta[\sigma^\smallfrown\Box \varphi/\sigma^\smallfrown\Box \varphi']\colon \delta \in W[\mod{T^i}_{\sigma^\smallfrown\Box \varphi}]\},\\
R^{i-1}_v(\sigma,\delta) \coloneqq& \begin{cases} 
              R^{i}_v(\sigma,\delta) &\hbox{ if }\sigma, \delta \in W^i_v,\\ 
              R^i_v(\sigma_-,\delta_-) &\hbox{ if } \sigma, \delta \in W^i_{\langle \sigma, \Box \varphi\rangle} \text{ for some } \langle \sigma, \Box \varphi\rangle \in (Gens_i^\Upsilon(\mod{T}))_v, \\
              r_{\eta^-1(\sigma),\Box \varphi} &\hbox{ if } \delta = \sigma^\smallfrown \Box \varphi' \text{ and } \langle \sigma, \Box \varphi\rangle \in (Gens_i^\Upsilon(\mod{T}))_v,\\
              \bot &\hbox{ otherwise, }\end{cases}\\
e^{i-1}_v(\sigma, p) \coloneqq& \begin{cases} e^i_v(\sigma, p) &\hbox{ if }\sigma \in W^i_v, \\ e^i_v(\sigma_-, p) &\hbox{ otherwise.} \end{cases}
 \end{align*}

 We define \[\mod{T}^+_v \coloneqq \mod{T}_v^0.\]
\end{definition}

As before, the Hereditarily Relevant worlds of $\mod{T}^+_v$ are those named with a sequence of formulas or primed formulas (hence, not including the names of any of the original worlds from $\mod{T}$).

Proving properties of this valued model is now simpler than checking the analogous results of the crisp case. As before, the only worlds that have fresh successors are those in the original model. Since the evaluation is now preserved fully (only the new accessibility relations change), we can prove easily in the following lemma that in the new worlds, the evaluation of all formulas is preserved.
\begin{lemma}\label{lem:newWorldsValues}
For any $1 \leq i \leq md(\Upsilon)$ and any $\sigma \in W^{i-1}_v\setminus W^i_v$ it holds that \[e^{i-1}_v(\sigma, \varphi) = e^i_v(\sigma_-, \varphi).\]
\end{lemma}
\begin{proof}
The lemma follows by induction on the complexity of the formula, with the propositional steps immediate from the definition of the evaluation $e^{i-1}_v$. 

For the modal steps, observe that if $\sigma \in W^{i-1}_v\setminus W^i_v$, necessarily $\sigma \not \in W^{md(\Upsilon)}$, and hence, $\langle \sigma, \psi\rangle \not \in (Gens_i^{\Upsilon}(\mod{T}))_v$ for any $\psi$. Also, by construction, $\sigma \in W^i_{\langle \delta, \Box \varphi\rangle}$ for some $\langle \delta, \Box \varphi\rangle \in (Gens_i^{\Upsilon}(\mod{T}))_v$. By definition, for $R^{i-1}_v(\sigma,\delta) > \bot$, necessarily $\delta \in W^i_{\langle \delta, \Box \varphi\rangle}$ too, and therefore, it holds that $R^{i-1}_v(\sigma,\delta) = R^i_v(\sigma_-,\delta_-)$. This means that the worlds in each $W^i_{\langle \delta, \Box \varphi\rangle}$ behave as an exact copy of the tree $\mod{T^i}_{\sigma^\smallfrown\Box \varphi}$. This immediately implies that the lemma holds also for modal formulas.
\end{proof}

It is now easy to conclude that the new model is conservative.

\begin{lemma}\label{lem:conservativeValuedModel}
Let $0 < i \leq md(\Upsilon)$. For any $\sigma \in W^i_v$ and any formula $\varphi \in \Upsilon_{\depth(\sigma)}$,

\[e^{i-1}_v(\sigma, \varphi) = e^i(\sigma, \varphi).\]
In particular, if $\sigma \in W^{md(\Upsilon)}_v$ (i.e., in $\eta^v(\mod{T})$), \[e^{i-1}_v(\sigma, \varphi) = e(\eta^{-1}(\sigma), \varphi).\]
\end{lemma}
\begin{proof} See Appendix 4.\end{proof}

The following are the (now much simpler) properties we will use to codify in a finite formula the infinitary characteristics of the unwitnessed situations.
\begin{obs}\label{obs:propsValued}
For any $\langle \sigma, \Box \varphi \rangle \in (Gens_i^{\Upsilon}(\mod{T}))_v$, for some $0 \leq i \leq md(\Upsilon)$ the following properties hold:
\begin{align}\label{eq:(P1)}\tag{P1}
R^+_v(\sigma, \sigma^\smallfrown \Box \varphi') \rightarrow e^+_v(\sigma^\smallfrown \Box \varphi', \varphi) <\ &  \mathbf{0} \text{ and }\\
\label{eq:(P2)}\tag{P2}
R^+_v(\sigma, \sigma^\smallfrown \Box \varphi') \rightarrow e^+_v(\sigma^\smallfrown \Box \varphi', \chi)   =\ & \mathbf{0} \text{ for each } \Box \chi \in \Upsilon_{\depth(\sigma)} \text{ s.t. } e^+_v(\sigma,\Box \chi) > \bot.
\end{align}
\end{obs}
\begin{proof}
Suppose $\langle \sigma, \Box \varphi \rangle \in (Gens_i^{\Upsilon}(\mod{T}))_v$ for some $i$. Then, $\sigma \in W^{md(\Upsilon)}_v$ and $R^{i-1}_v(\sigma, \sigma^\smallfrown \Box \varphi') = r_{\eta^{-1}(\sigma), \Box \varphi}$. The corollary follows from the  properties of this latter value (see Lemma \ref{lem:belowWorldValued} and below). 
\end{proof}

We will keep using the set of sequences $\Sigma$ from Definition \ref{def:possibleWorlds} and the coherent sets from Definition  \ref{def:coherentSet}, even if now, as we will see, for the unwitnessed situations (namely, where both $\sigma^\smallfrown\Box \chi \in \Omega$ and $\sigma^\smallfrown\Box \chi' \in \Omega$) we will only care about what happens in the world $\sigma^\smallfrown\Box \chi'$.
The fact that the hereditarily relevant worlds of a model $\mod{T}^+_v$ form a coherent set is proven as in Lemma \ref{lemma:HRWcoherent}.

The set of propositional formulas that will encode the previous information is now based on a slightly different set of propositional variables than those from Definition \ref{def:setOfVars}. Given a coherent set $\Omega$, we need the basic variables obtained from the formulas in $\Upsilon$ paired with their corresponding elements from $\Omega$, but instead of the $\alpha$-values (no longer used), we need to add variables that will encode the value of the accessibility relation. Namely, recall (from before Definition \ref{def:setOfVars})  that $vars_i \coloneqq pV(\Upsilon_i) \cup \{\heartsuit \psi\colon \heartsuit \in \Upsilon_i\}$ and $\Omega_i \coloneqq \{\sigma \in \Omega\colon \vert \sigma \vert = i\}$. Then we let

\[\mathcal{V}^{\Omega}_v \coloneqq \bigcup_{0 \leq i \leq md(\Upsilon)} \{\delta_\sigma\colon \delta\in vars_i, \sigma\in \Omega_{i+1}\} \cup \bigcup \{r_{\sigma, \sigma^\smallfrown\alpha }\colon \sigma, \sigma^\smallfrown\alpha \in \Omega\}.\]

Now, for a coheret subset of $\Omega$, we define the following sets of formulas (over the language $\mathcal{V}^\Omega_v$):
\[\mathcal{EP}^{\Omega}_v, (\mathcal{W}it^{\Omega}_\Diamond)_v, (\mathcal{W}it^{\Omega}_\Box)_v, (u\mathcal{W}it^{\Omega}_\Diamond)_v, (\mathcal{C}trl^\Omega_\Box)_v\]

where: 
\begin{enumerate}[\label={}]
\item 

\begin{fleqn}\begin{align*}
\bullet\ \mathcal{EP}^{\Omega}_v \coloneqq  \bigcup_{\sigma \in \Omega, \not \exists  \sigma^\smallfrown \beta \in \Omega } 
 \{ (\Box \psi)_\sigma \colon (\Box \psi)_\sigma \in \mathcal{V}^\Omega\}\cup \{ \neg (\Diamond \psi)_\sigma \colon (\Diamond \psi)_\sigma \in \mathcal{V}^\Omega\} .
 \end{align*}
\end{fleqn}


\textit{If $\sigma$ is understood as an end-point of $\Omega$ -i.e., no $\sigma^\smallfrown \beta$ belongs to $\Omega$, so it is understood that the accessibility to any world is $\bot$-, the modal formulas are evaluated to the corresponding value ($\mathbf{0}$ for $\Box$ formulas or $\bot$ for $\Diamond$ ones).}

%
\item 

\begin{fleqn}\begin{align*}
\bullet\ (\mathcal{W}it^{\Omega}_\Diamond)_v \coloneqq  \{W_\Diamond((\Diamond \psi)_\sigma) \colon \sigma^\smallfrown\Diamond \psi \in \Omega\}
\end{align*}\end{fleqn}
 for 
\[
W_\Diamond((\Diamond \psi)_\sigma) \coloneqq
((\Diamond \psi)_\sigma \leftrightarrow (r_{\sigma, \sigma^\smallfrown\Diamond \psi} \with \psi_{\sigma^\smallfrown\Diamond \psi})) \wedge  (\bigvee_{\sigma^\smallfrown \beta \in \Omega}(r_{\sigma, \sigma^\smallfrown \beta} \with \psi_{\sigma^\smallfrown \beta}) \rightarrow (\Diamond \psi)_\sigma).\]

\textit{The value of the $\Diamond$ formulas is given following the corresponding definition of the Kripke structure.}

\item 

\begin{fleqn}\begin{align*}
\bullet\ (\mathcal{W}it^{\Omega}_\Box)_v \coloneqq  \{W_\Box((\Box \psi)_\sigma) \colon \sigma^\smallfrown\Box \psi \in \Omega, \sigma^\smallfrown\Box \psi' \not \in \Omega\}
 \end{align*}
\end{fleqn}
for 
\[
W_\Box((\Box \psi)_\sigma) \coloneqq
((\Box \psi)_\sigma \leftrightarrow (r_{\sigma, \sigma^\smallfrown\Box \psi} \rightarrow \psi_{\sigma^\smallfrown\Box \psi})) \wedge  ((\Box \psi)_\sigma \rightarrow \bigwedge_{\sigma^\smallfrown \beta \in \Omega}(r_{\sigma, \sigma^\smallfrown \beta} \rightarrow \psi_{\sigma^\smallfrown \beta})).\]
\textit{The value of the $\Box$ formulas that are undestood as witnessed (i.e., the sequence ending with that formula primed does not belong to $\Omega$) is given following the corresponding definition of the Kripke structure.}

\item 

\begin{fleqn}\begin{align*}
\bullet\ (u\mathcal{W}it^{\Omega}_\Box)_v \coloneqq  \{\neg (\Box \psi)_\sigma\colon \sigma ^\smallfrown\Box \psi' \in \Omega\}.
 \end{align*}
\end{fleqn}
\textit{The value of any $\Box$ formula that is undestood as unwitnessed (i.e., the sequence ending with that formula primed does belong to $\Omega$) is equal to $0$ (following quasi-witnessed completeness).}

\item 

\begin{fleqn}\begin{align*}
\bullet\ (\mathcal{C}trl^{\Omega}_\Box)_v \coloneqq& \{\neg \neg (\Box \psi)_\sigma \rightarrow (r_{\sigma, \sigma^\smallfrown\Box \varphi'} \rightarrow \psi_{\sigma^\smallfrown\Box \varphi'}) \colon \sigma^\smallfrown\Box \varphi'\in \Omega, (\Box \psi)_\sigma \in \mathcal{V}^\Omega_v\}\cup\\
& \{\neg \Delta (r_{\sigma, \sigma^\smallfrown\Box \varphi'} \rightarrow \varphi_{\sigma^\smallfrown\Box \varphi'}) \colon \sigma^\smallfrown\Box \varphi'\in \Omega\}.
 \end{align*}
\end{fleqn}

\textit{This follows properties ($\ref{eq:(P1)}$) and ($\ref{eq:(P2)}$) from the previous page.}

\end{enumerate}

\begin{definition}[c.f. Definition \ref{def:definingFormulas} ]\label{def:definingFormulasValued}
Let $\Upsilon$ be a finite set of modal formulas.
We let\[\Upsilon^\star \coloneqq \bigvee_{\Omega \underset{\text{coh}}{\subseteq} \Sigma}\bigwedge \Upsilon^\Omega_v, \]

for \[\Upsilon^\Omega_v \coloneqq \mathcal{EP}^{\Omega}_v \cup (\mathcal{W}it^{\Omega}_\Diamond)_v \cup (\mathcal{W}it^{\Omega}_\Box)_v \cup (u\mathcal{W}it^{\Omega}_\Box)_v \cup (\mathcal{C}trl^{\Omega}_\Box)_v\]
\end{definition}

We will now let 
\[uWV'(\Upsilon, \Omega) \coloneqq \bigvee_{\sigma^\smallfrown \Box \varphi'\in \Omega} r_{\sigma, \sigma^\smallfrown \Box \varphi'}\rightarrow \varphi_{\sigma^\smallfrown \Box \varphi'}.\]
 
We can now rephrase Proposition \ref{prop:ModToProp} as follows:

 \begin{proposition}\label{prop:ModToPropValued}
Let $\Upsilon = \Gamma \cup \{\varphi\}$ be such that $\Gamma \not \Vdash_{M\Pi} \varphi$. Then 
\[\Gamma_{\langle 0 \rangle}, \Upsilon^\star_v
\not \models_{\Pi_\Delta} \varphi_{\langle 0 \rangle}\]

\end{proposition}
\begin{proof}
The proof proceeds as the one of Proposition \ref{prop:ModToProp}. For if $\Gamma \not \Vdash_{M\Pi} \varphi$, there is a quasi-witnessed valued 
$\mathbf{R}$-tree $\mod{T}$, with root $r$, such that $\Gamma \not \Vdash_{\mod{T},r} \varphi$, and from Lemma \ref{lem:conservativeValuedModel}, we also know that $\Gamma \not \Vdash_{\mod{T}^+_v,\langle 0 \rangle} \varphi$. Let the set $\Omega = HRW(\mod{T}^+_v)$, which is coherent, and let $h$ be the propositional $\Pi_\Delta$-homomorphism into $\mathbf{R}_\Delta$ given by 
\begin{align*}
h(p_\sigma) \coloneqq & e^+_v(\sigma, p) \text{ for }p \in \mathcal{V}(\Upsilon), \sigma \in \Omega,\\
h((\heartsuit\psi)_\sigma) \coloneqq & e^+_v(\sigma, \heartsuit \psi) \text{ for } \heartsuit \in \{\Box, \Diamond\}, \sigma \in \Omega,\\
h(r_{\sigma, \sigma^\smallfrown\alpha}) \coloneqq & R^+_v(\sigma, \sigma^\smallfrown\alpha) \text{ for }  \sigma, \sigma^\smallfrown\alpha \in \Omega.
\end{align*}
Checking that every formula in $\Upsilon^\Omega_v$ is evaluated to $\mathbf{0}$ under $h$ is done as in the proof of Proposition \ref{prop:ModToProp}, following the definition of the hereditarily relevant worlds, and the properties from Observation \ref{obs:propsValued}.
\end{proof}

To go back from the above propositional condition to a standard (valued) product Kripke model is now much simpler than the construction we did in the previous section. It is based on the well-known fact that  given any product homomorphism $h$, for any $k \in \omega$, the map $h^k$ is a product homomorphism too and $h^k(\varphi) = (h(\varphi))^k$.

\begin{proposition}[c.f. Proposition \ref{prop:PropToMod}]
Let $\Upsilon = \Gamma \cup \{\varphi\} \subseteq_{\omega} Fm$ and $\Omega$ a coherent subset of $\Sigma$ such that 

\[\Gamma_{\langle 0 \rangle}, \Upsilon^\star_v
\not \models_{\Pi} \varphi_{\langle 0 \rangle}.\]

Then, $\Gamma \not \Vdash_{M[0,1]_\Pi} \varphi$. 
\end{proposition}
\begin{proof}
From the premises of the proposition, there is a propositional homomorphism $h$ such that $h(\Gamma_{\langle 0 \rangle}) \subseteq \{1\}$, $h(\Upsilon^\star_v)\subseteq \{1\}$ and $h(\varphi_{\langle 0 \rangle}) < 1$. 

Recall, for a sequence $\sigma \in \Sigma$, the definitions of $\muwave{\sigma}$ and of $prods(\sigma)$, given in the proof of Proposition \ref{prop:PropToMod}. Now, for any element $\nu \in prods(\sigma)$ for some $\sigma \in \Sigma$, we let 
\[\Theta(\langle 0 \rangle) = 0, \qquad \Theta(\nu^\smallfrown \heartsuit \varphi) = \Theta(\nu), \qquad \Theta(\nu^\smallfrown (\Box \varphi)_k) \coloneqq \Theta(\nu) + k.\]
The intuition behind this value is the ``power" to which the values at that world will be raised (with respect to the value $1$).

We define the $[0,1]$-valued model $\mod{T}$  as follows:
\begin{itemize}
\item $W \coloneqq \bigcup_{\sigma \in \Omega}prods(\sigma)$,
\item $R(\sigma, \sigma^\smallfrown \alpha) \coloneqq h(r_{\muwave{\sigma}, \muwave{\sigma^\smallfrown \alpha}})^{\Theta(\sigma^\smallfrown \alpha)}$
\item $e(\sigma, p) \coloneqq h(p_{\muwave{\sigma}})^{\Theta(\sigma)}$
\end{itemize}
 
We check below that $e(\sigma, \varphi) = h(\varphi_{\muwave{\sigma}})^{\Theta(\sigma)}$ for any $\sigma \in W$ and $\varphi \in \Upsilon_{\vert \sigma \vert -1}$.
This implies that $e(\langle 0 \rangle, \Gamma) = \{1\}$ and $e(\langle 0 \rangle, \varphi) < 1$, concluding the proposition.

For propositional formulas, the claim follows since, as we remarked before, $h^k(\varphi) = (h(\varphi))^k$.
For formulas beginning with $\Box$, by the definitions of $\mod{T}$ and of $\Sigma$, and Induction Hypothesis\footnote{Recall that in the standard product algebra $[0,1]_\Pi$,  $a^k \rightarrow b^k = (a \rightarrow b)^k$ for any $a,b \in [0,1], k \in \omega$.} we get that  
\begin{align*}
e(\sigma, \Box \psi) =
 \bigwedge_{\sigma^\smallfrown \alpha \in W} R(\sigma, \sigma^\smallfrown\alpha) \rightarrow e(\sigma^\smallfrown \alpha, \psi) & =\\
  \bigwedge_{\sigma^\smallfrown \alpha \in W} h(r_{\muwave{\sigma}, \muwave{\sigma^\smallfrown \alpha}})^{\Theta(\sigma^\smallfrown \alpha)} \rightarrow
h(\psi_{\muwave{\sigma^\smallfrown \alpha}})^{\Theta(\sigma^\smallfrown \alpha)} & =\\
\bigwedge_{\sigma^\smallfrown \heartsuit \chi \in W} (h(r_{\muwave{\sigma}, \muwave{\sigma^\smallfrown \heartsuit \chi}}) \rightarrow
h(\psi_{\muwave{\sigma^\smallfrown \heartsuit \chi}}))^{\Theta(\sigma)} \wedge
 \bigwedge_{\sigma^\smallfrown (\Box \chi)_k \in W} (h(r_{\muwave{\sigma}, \muwave{\sigma^\smallfrown (\Box \chi)_k}}) & \rightarrow
h(\psi_{\muwave{\sigma^\smallfrown (\Box \chi)_k}}))^{\Theta(\sigma^\smallfrown (\Box \chi)_k)}.
\end{align*}

Consider first the case when $\muwave{\sigma}^\smallfrown \Box \psi' \not \in \Omega$. From $(\mathcal{W}it_\Box)_v^\Omega$ we have that 
\begin{equation}\label{eq:BoxCondition}
h((\Box \psi)_{\muwave{\sigma}}) = h(r_{\muwave{\sigma}, \muwave{\sigma}^\smallfrown \Box \psi}) \rightarrow h(\psi_{\muwave{\sigma}^\smallfrown \Box \psi}) \leq h(r_{\muwave{\sigma}, \muwave{\sigma}^\smallfrown \alpha}) \rightarrow h(\psi_{\muwave{\sigma}^\smallfrown \alpha})\text{ for any }\muwave{\sigma}^\smallfrown \alpha \in \Omega.
\end{equation}

If $h((\Box \psi)_{\muwave{\sigma}}) = 0$, from $(\mathcal{W}it^{\Omega}_\Box)_v$ it follows that $h((\Box \psi)_{\muwave{\sigma}}) = h(r_{\muwave{\sigma}, \muwave{\sigma}^\smallfrown \Box \psi} \rightarrow \psi_{\muwave{\sigma}^\smallfrown \Box \psi})$. Thus, also 
 \[\bigwedge_{\sigma^\smallfrown \heartsuit \chi \in W} (h(r_{\muwave{\sigma}, \muwave{\sigma}^\smallfrown \heartsuit \chi}) \rightarrow
h(\psi_{\muwave{\sigma}^\smallfrown \heartsuit \chi}))^{\Theta(\sigma)} \leq h(r_{\muwave{\sigma}, \muwave{\sigma}^\smallfrown \Box \psi}) \rightarrow h(\psi_{\muwave{\sigma}^\smallfrown \Box \psi}) = 0.\]
 Therefore, 
$e(\sigma, \Box \psi) = 0 = h((\Box \psi)_{\muwave{\sigma}})  = h((\Box \psi)_{\muwave{\sigma}})^{\Theta(\sigma)}$.


On the other hand, if $h((\Box \psi)_{\muwave{\sigma}}) > 0$, we know from $(\mathcal{C}trl^\Omega_\Box)_v$ that $h(r_{\muwave{\sigma}, \muwave{\sigma}^\smallfrown \Box \chi'}) \rightarrow h(\psi_{\muwave{\sigma}^\smallfrown \Box \chi'}) = 1$ for any
 $\muwave{\sigma}^\smallfrown \Box \chi' \in \Omega$. Hence 
 $(h(r_{\muwave{\sigma}, \muwave{\sigma^\smallfrown (\Box \chi)_k}}) \rightarrow h(\psi_{\muwave{\sigma^\smallfrown (\Box \chi)_k}}))^{\Theta(\sigma^\smallfrown (\Box \chi)_k)} = 1$ too, for any $\sigma^\smallfrown (\Box \chi)_k \in W$. 
 Any other $\sigma^\smallfrown\alpha \in W$ is of the form $\sigma^\smallfrown\heartsuit \chi$, and so, $\Theta(\sigma^\smallfrown \heartsuit \chi)= \Theta(\sigma)$ by definition. All the above, using Equation (\ref{eq:BoxCondition}), imply that

\begin{align*}
e(\sigma, \Box \psi) =
\bigwedge_{\sigma^\smallfrown \heartsuit \chi \in W} (h(r_{\muwave{\sigma}, \muwave{\sigma}^\smallfrown \heartsuit \chi}) \rightarrow
h(\psi_{\muwave{\sigma}^\smallfrown \heartsuit \chi}))^{\Theta(\sigma)} \wedge 1 =
h((\Box \psi)_{\muwave{\sigma}})^{\Theta(\sigma)}.
\end{align*}

For the second case, assume that $\muwave{\sigma}^\smallfrown \Box \psi' \in \Omega$. Then, on the one hand, by $(u\mathcal{W}it_\Box^\Omega)_v$, necessarily $h((\Box \psi)_{\muwave{\sigma}}) = 0$. On the other, we know from $(\mathcal{C}trl^\Omega_\Box)_v$ that $h( r_{\muwave{\sigma}, \muwave{\sigma}^\smallfrown \Box \psi'}\rightarrow \psi_{\muwave{\sigma}^\smallfrown \Box \psi'}) < 1$.
Lastly, by the definition of $W$, we know that if $\muwave{\sigma}^\smallfrown \Box \psi' \in \Omega$ and $\sigma \in W$, then, for any $k \in \omega$, we have $\sigma^\smallfrown (\Box \psi)_k \in W$ (and clearly, $\Theta(\sigma^\smallfrown (\Box \psi)_k) \geq k$). It follows that

\begin{align*}
\bigwedge_{\sigma^\smallfrown (\Box \chi)_k \in W} (h(r_{\muwave{\sigma}, \muwave{\sigma^\smallfrown (\Box \chi)_k}}) \rightarrow
h(\psi_{\muwave{\sigma^\smallfrown (\Box \chi)_k}}))^{\Theta(\sigma^\smallfrown (\Box \chi)_k)} & \leq \\
\bigwedge_{\sigma^\smallfrown (\Box \psi)_k\colon k \in \omega} (h(r_{\muwave{\sigma}, \muwave{\sigma^\smallfrown (\Box \psi)_k}}) \rightarrow
h(\psi_{\muwave{\sigma^\smallfrown (\Box \psi)_k}}))^k & = 0.
\end{align*}

Therefore, $e(\sigma, \Box \psi) = 0 = h((\Box \psi)_{\muwave{\sigma}}) = h((\Box \psi)_{\muwave{\sigma}})^{\Theta(\muwave{\sigma})}$.

The case for formulas beginning with $\Diamond$ is simpler. Indeed, analogously to the previous case, 
\begin{align*}
e(\sigma, \Diamond \psi) =
 \bigvee_{\sigma^\smallfrown \alpha \in W} R(\sigma, \sigma^\smallfrown\alpha) \with e(\sigma^\smallfrown \alpha, \psi) & =\\
  \bigvee_{\sigma^\smallfrown \alpha \in W} h(r_{\muwave{\sigma}, \muwave{\sigma^\smallfrown \alpha}})^{\Theta(\sigma^\smallfrown \alpha)} \with
h(\psi_{\muwave{\sigma^\smallfrown \alpha}})^{\Theta(\sigma^\smallfrown \alpha)} & =\\
\bigvee_{\sigma^\smallfrown \heartsuit \chi \in W} (h(r_{\muwave{\sigma}, \muwave{\sigma^\smallfrown \heartsuit \chi}}) \with
h(\psi_{\muwave{\sigma^\smallfrown \heartsuit \chi}}))^{\Theta(\sigma)} \vee
 \bigvee_{\sigma^\smallfrown (\Box \chi)_k \in W} (h(r_{\muwave{\sigma}, \muwave{\sigma^\smallfrown (\Box \chi)_k}}) & \with
h(\psi_{\muwave{\sigma^\smallfrown (\Box \chi)_k}}))^{\Theta(\sigma^\smallfrown (\Box \chi)_k)}.
\end{align*}
Observe that $(h(r_{\muwave{\sigma}, \muwave{\sigma^\smallfrown (\Box \chi)_k}}) \with
h(\psi_{\muwave{\sigma^\smallfrown (\Box \chi)_k}}))^{\Theta(\sigma^\smallfrown (\Box \chi)_k)} \leq (h(r_{\muwave{\sigma}, \muwave{\sigma}^\smallfrown \Box \chi'}) \with
h(\psi_{\muwave{\sigma}^\smallfrown \Box \chi'}))^{\Theta(\sigma)}$, so the second disjunction can be ignored. 

By $(\mathcal{W}it_\Diamond^\Omega)_v$ we also know that $h((\Diamond \psi)_{\muwave{\sigma}}) = h(r_{\muwave{\sigma}, \muwave{\sigma}^\smallfrown\Diamond \psi}) \with h(\psi_{\muwave{\sigma}^\smallfrown\Diamond \psi}) \geq h(r_{\muwave{\sigma}, \muwave{\sigma}^\smallfrown\alpha}) \with h(\psi_{\muwave{\sigma}^\smallfrown\alpha})$ for any $\muwave{\sigma}^\smallfrown\alpha \in \Omega$. From the above, we immediately conclude that
\begin{align*}
e(\sigma, \Diamond \psi) =
\bigvee_{\sigma^\smallfrown \heartsuit \chi \in W} (h(r_{\muwave{\sigma}, \muwave{\sigma^\smallfrown \heartsuit \chi}}) \with
h(\psi_{\muwave{\sigma^\smallfrown \heartsuit \chi}}))^{\Theta(\sigma)}  =
 (h(r_{\muwave{\sigma}, \muwave{\sigma}^\smallfrown\Diamond \psi}) \with h(\psi_{\muwave{\sigma}^\smallfrown\Diamond \psi}))^{\Theta(\sigma)}&  =  h((\Diamond \psi)_{\muwave{\sigma}})^{\Theta(\sigma)}. \qedhere
 \end{align*}
 
\end{proof}

\noindent
\textbf{Acknowledgments.} This work has been funded by the ERDF-Project "Knowledge in the Age of Distrust" from the Programme Johannes Amos Comenius under the Ministry of Education, Youth and Sports of the Czech Republic, CZ.02.01.01/00/23\_025/0008711. This work has also received funding from the PIE CSIC 20235AT019 project, and from the European Union’s Horizon 2020 research and innovation programme under the Marie Skłodowska-Curie grant agreement No. 101027914.

 The author would like to express her deep gratitude to the reviewers for their valuable comments and suggestions, which have substantially contributed to improving the readability and presentation of this work.

%
%
%
%
%
%
%


\bibliographystyle{abbrv}
\bibliography{extracted.bib}

\section*{Appendix 1. Some proofs from Section \ref{sec:Prelim}.}\label{sec:app1}

\begin{proof}[Proof of Proposition \ref{prop:compTreesFiniteDepth}]
The proof is analogous to the one from the classical modal logic, simply taking into account the valued accessibility and evaluation. We provide the details for the interested reader. 

We will denote by $[u_1, \ldots u_n]$ a list of $n$ elements, and allow us to denote this by $\sigma$. By $\sigma^\smallfrown u$ we denote the list $[u_1, \ldots, u_n, u]$, and by $\mathrm{l}(\sigma)$, we denote the element $u_n$ (the last one). 

Let us define the tree $\mod{T} \coloneqq \langle W^{\mod{T}}, R^{\mod{T}}, e^{\mod{T}}\rangle$ where 
\begin{itemize}
\item $W^{\mod{T}} \coloneqq \bigcup_{i \in \omega} W_i$, for
\begin{align*}
W_0 \coloneqq& \{[v]\}, \\
W_{i+1} \coloneqq& \{\sigma^\smallfrown u\colon \sigma \in W_i, u \in W, R(\mathrm{l}(\sigma), u) > 0\};
\end{align*}
\item $R^{\mod{T}}(\sigma, \delta) \coloneqq \begin{cases} R(\mathrm{l}(\sigma), \mathrm{l}(\delta)) &\hbox{ if } \delta = \sigma^\smallfrown u \text{ for some }u \in W, \\
0 &\hbox{ otherwise; }\end{cases}$
\item $e^{\mod{T}}(\sigma, p) \coloneqq e(\mathrm{l}(\sigma), p)$.
\end{itemize}

It is first easy to check that $\mod{T}$ is a tree. Indeed, consider the splitting given by the construction of $W^{\mod{T}}$ (in which $W_i \cap W_j = \emptyset$ for any $i \neq j$ by definition, simply because the elements of $W_i$ have exactly length $i+1$). Clearly $W_0$ is a singleton. 
Furthermore, if $R^{\mod{T}}(\sigma, \delta) > 0$, then by definition $\delta = \sigma^\smallfrown u$ for some $u \in W$. Therefore, if $\sigma \in W_i$ for some $i$, necessarily $\delta \in W_{i+1}$.

On the other hand, let $\sigma^\smallfrown u \in W_{i+1}$ for some $i$, so $l(\sigma^\smallfrown u ) = u$. By definition of $W_{i+1}$ follows that $R^{\mod{T}}(\sigma, \sigma^\smallfrown u) = R(l(\sigma),u) > 0$. Furthermore, for any other $\delta$ different from $\sigma$, $R^{\mod{T}}(\delta, \sigma^\smallfrown u) = 0$.

We claim that, for any formula $\psi$, \[e^{\mod{T}}(\sigma, \psi) = e(\mathrm{l}(\sigma), \psi).\] This will conclude the proposition, since in particular the values of $\Gamma$ and $\varphi$ will be preserved from $v$ to the root of $\mod{T}$, $v'$, which is by definition the list $[v]$.

The case of propositional formulas is immediate by the definition of $e^{\mod{T}}$. 
Consider a modal formula $\Box \psi$. By definition and I.H., 
\[e^{\mod{T}}(\sigma, \Box \psi) = \bigwedge_{\sigma^\smallfrown u \in W^{\mod{T}}} R(\mathrm{l}(\sigma), \mathrm{l}(\sigma^\smallfrown u)) \rightarrow e^{\mod{T}}(\sigma^\smallfrown u, \psi) = \bigwedge_{\sigma^\smallfrown u \in W^{\mod{T}}} R(\mathrm{l}(\sigma), u) \rightarrow e(u, \psi)\] the last equality by I.H. 
 
Observe $\sigma^\smallfrown u \in W^{\mod{T}}$ if and only if $R(\mathrm{l}(\sigma), u) > 0$, therefore the above is equal to $\bigwedge_{u \in W} R(\mathrm{l}(\sigma), u) \rightarrow e(u, \psi)$, which is, by definition, $e(\mathrm{l}(\sigma), \Box \psi)$.

The proof for the formulas beginning with a $\Diamond$ operator is analogous.

Checking that $\mod{T}$ preserves crispness, and witnessing and quasi-witnessing conditions is a matter of basic calculations. 
 \end{proof}

\section*{Appendix 2. Some proofs from Section \ref{sec:soundness}.}\label{sec:app2}

\begin{proof}[Proof of Observation \ref{obs:ordermapsto}]
The first claim is obvious by definition. 

For the second claim, the non trivial case is that of $a < b$, which holds only if there is $p \in \mathds{Q}$ such that $a(p) < b(p)$ and $a(s) = b(s)$ for all $s < p$. If $q  < p$, naturally $a_{\mapsfrom q} = b_{\mapsfrom q}$ since $a_{\mapsfrom q}(s) = a(s) = b(s) = b_{\mapsfrom q}(s)$ for all $s \leq q$ (given that $q < p$) and  $a_{\mapsfrom q}(s) = 0 =  b_{\mapsfrom q}(s)$ for all $s > q$. On the other hand, if $p \leq q$, we have that $a_{\mapsfrom q}(p) = a(p) < b(p) = b_{\mapsfrom q}(p)$, while $a_{\mapsfrom q}(s) = a(s) = b(s) = b_{\mapsfrom q}(s)$ for all $s < p$ (since $p \leq q$), hence $a_{\mapsfrom q} \leq b_{\mapsfrom q}$. 

The third claim is trivial if $a = b$, and also if $a = \bot$. Otherwise, by the definition of $_{\mapsfrom}$ 
\[(x + x_{\mapsfrom q})(s) = \begin{cases} 2x(s) &\hbox{ if } s\leq q,\\ x(s) &\hbox {otherwise.}\end{cases}\] The claim follows from this after few basic calculations, using that the order is the lexicographic one.

For the last case, by definition \[a_{\mapsfrom q}[s] + b_{\mapsfrom q}[s] = \begin{cases} a[s] + b[s] &\hbox{if }s \leq q\\ 0 &\hbox{otherwise} \end{cases}\] which clearly coincides with $(a + b)_{\mapsfrom q}$. The analogous holds for the $-$ operation.

\end{proof}

\begin{proof}[Proof of Lemma \ref{lem:infSupmapsfrom}]

Recall that for any $a,b\in \mathbf{R}$, $a + b \leq a$.

To prove the first condition, assume first that $\bigwedge X = \bot$. It is immediate that $\bigwedge_{x \in X}x + (\bigwedge_{x \in X} x)_{\mapsfrom q} = \bot$ too, and also, that also \[\bigwedge_{x \in X} (x + x_{\mapsfrom q}) \leq \bigwedge_{x \in X} x = \bot.\] 

Suppose then that $\bigwedge X = m > \bot$, and $m \in X$. It is first clear that, since $m \leq x$ for all $x \in X$, 
$m + m_{\mapsfrom q} \leq x + x_{\mapsfrom q}$ for all $x \in X$ too. Therefore \[\bigwedge X + (\bigwedge X)_{\mapsfrom q} \leq \bigwedge_{x \in X}(x + x_{\mapsfrom q}).\]
On the other hand, since $m \in X$, it also holds that $\bigwedge_{x \in X}(x + x_{\mapsfrom q}) \leq m + m_{\mapsfrom q} = \bigwedge X + (\bigwedge X)_{\mapsfrom q}$, concluding the proof of the case.

For the supremum case, assume $\bigvee X = m \in X$. Then $x + x_{\mapsfrom q} \leq m + m_{\mapsfrom q}$ for all $x \in X$, and therefore $\bigvee_{x \in X} (x + x_{\mapsfrom q}) \leq m + m_{\mapsfrom q} = \bigvee X + (\bigvee X)_{\mapsfrom q}$. On the other hand, since $m \in X$ it also holds that 
$m + m_{\mapsfrom q} \leq \bigvee_{x \in X} (x + x_{\mapsfrom q})$, concluding the proof.

\end{proof}

\begin{proof}[Proof of Corollary \ref{cor:quasiWitComp}]
Suppose $\Gamma \not \vdash_{K\Pi} \varphi$. Then, there is some model in $K\Pi$ (i.e., crisp) over which $\Gamma \not \Vdash_{\mod{M}} \varphi$. From Prop. \ref{prop:modalToFo} we get $\Gamma^l \not \Vdash_{\mod{M}^\times} \varphi^l$. Moreover, following Obs. \ref{obs:languageModalToFO}, $\mod{M}^\times$ is crisp (in the predicates sense), so it satisfies $\forall x,yR(x,y) \vee \neg R(x,y)$. Therefore, also $\Gamma^l, \forall x,yR(x,y) \vee \neg R(x,y) \not \Vdash_{\mod{M}^\times} \varphi^l$. From completeness of $\forall \Pi$, we get that $\Gamma^l, \forall x,yR(x,y)  \vee \neg R(x,y) \not \vdash_{\forall \Pi} \varphi^l$, and now, using Theorem \ref{Thm:LaMa07}, we know there is some countable, quasi-witnessed $\mathbf{R}$-structure $\mod{N}$ such that $\Gamma^l, \forall x,yR(x,y)  \vee \neg R(x,y) \not \Vdash_{\mod{N}} \varphi^l$. By definition, $\mod{N}$ is crisp, and thus, from Obs. \ref{obs:languageModalToFO}, also $\mod{N}_+$ is crisp (as a Kripke model), and also quasi-witnessed over $\mathbf{R}$. From Prop. \ref{prop:modalToFo}, it follows that $\Gamma \not \Vdash_{\mod{N}_+} \varphi$. Applying Corollary \ref{cor:LocalFiniteDepth} to this latter model concludes the corollary.
\end{proof}

\begin{proof}[Proof of Lemma \ref{lem:inverseEtaPreservation}]
It can be checked by induction on the complexity of the formula. For propositional letters and connectives it follows from the definition of $\eta(e)$. 

Consider the case of formulas beginning by a $\Box$ operator. By definition and applying I.H. we have that 
\[\eta(e)(\sigma, \Box \varphi) = \bigwedge_{\sigma ^\smallfrown \alpha \in \eta(W)} e(\eta^{-1}(\sigma^\smallfrown \alpha), \varphi) \qquad \text{ and }\qquad e(\eta^{-1}(\sigma), \Box \varphi) = \bigwedge_{R\eta^{-1}(\sigma),v} e(v, \varphi).\]

Let us first check that $e(\eta^{-1}(\sigma), \Box \varphi) \leq \eta(e)(\sigma, \Box \varphi)$. 

For $\sigma ^\smallfrown \alpha \in \eta(W)$, there are two cases. First, if $\alpha =  \heartsuit \psi$ for some $ \heartsuit \psi$, then $\eta^{-1}(\sigma ^\smallfrown \heartsuit \psi) = \eta^{-1}(\sigma)_{\heartsuit \psi}$ and $e(\eta^{-1}(\sigma ^\smallfrown \heartsuit \psi), \varphi) = e(\eta^{-1}(\sigma)_{\heartsuit \psi}, \varphi)$. Since by definition $R \eta^{-1}(\sigma),\eta^{-1}(\sigma)_{\heartsuit \psi}$, we get that $\bigwedge_{R\eta^{-1}(\sigma),v} e(v, \varphi) \leq e(\eta^{-1}(\sigma)_{\heartsuit \varphi}, \varphi)$.
 Second, if $\alpha =  u$ for some $u \in W$, then we know that $ \eta^{-1}(\sigma ^\smallfrown u) = u$ and so that also $R \eta^{-1}(\sigma),u$. The rest of the case is as the previous one. 

To check that $\eta(e)(\sigma, \Box \varphi) \leq e(\eta^{-1}(\sigma), \Box \varphi)$, consider any $v$ such that $R\eta^{-1}(\sigma),v$. We will prove that there is $u_v \in \eta(W)$ with $\eta(R)(\sigma, u_v)$ and such that $\eta(e)(u_v, \varphi) = e(v, \varphi)$. By I.H., this would be equivalent to show that there is an $u_v$ with the two previous conditions and such that $e(\eta^{-1}(u_v), \varphi) = e(v,\varphi)$. 

Since $R\eta^{-1}(\sigma),v$, we have that $\eta^{-1}(\sigma) = par(v)$, and so, $\sigma \in \eta(par(v))$. If $v \not \in RW(\mod{T})$, then $\sigma^\smallfrown v \in \eta(v) \subseteq \eta(W)$ and $\eta(R) \sigma,\sigma^\smallfrown v $ following the definition, and moreover $v = \eta^{-1}(\sigma^\smallfrown v)$. Thus $e(\eta^-1(\sigma^\smallfrown v), \varphi) = e(v, \varphi)$, proving the condition above.  

  Otherwise (if $v \in RW(\mod{T})$), we know that $v = \eta^{-1}(\sigma)_{\heartsuit \psi}$ for some modal formula $\heartsuit \psi$. This implies that $\sigma^\smallfrown\heartsuit \psi \in \eta(v)\subseteq \eta(W)$ and $\eta(R)\sigma,\sigma^\smallfrown\heartsuit \psi $ again,  and also, $v = \eta^{-1}(\sigma^\smallfrown\heartsuit \psi)$. Hence, immediately, $e(\eta^{-1}(\sigma^\smallfrown\heartsuit \psi), \varphi) = e(v, \varphi)$. The proof for formulas beginning by $\Diamond$ is done analogously. 

\end{proof}

\begin{proof}[Proof of Observation \ref{obs:freshSuccessorsOnlyOriginalWorlds}]
To prove the first point, we know that if a world $\delta$ was added (fresh) at step $i-1$ then, by definition, it belongs to $W^i_{\langle \sigma, \Box \varphi\rangle}$ for some $\langle \sigma, \Box \varphi\rangle \in Gens_i^\Upsilon(\mod{T})$. Among other things, this implies that $\sigma \in \eta(RW(\mod{T}))$, so $\sigma = \sigma_-$. 
By definition, the above means that $\delta = \rho[\sigma^\smallfrown\Box \varphi/\sigma^\smallfrown\Box \varphi']$ for some $\rho \in W[\mod{T}^i_{\sigma^\smallfrown\Box \varphi}]$, thus for some $\rho = \sigma^\smallfrown\Box \varphi^\smallfrown \delta_2$.
Therefore, $\delta = \sigma^\smallfrown\Box \varphi'^\smallfrown\delta_2$ for $\langle \sigma, \Box \varphi\rangle \in Gens_i^\Upsilon(\mod{T})$.

To prove the second point, following the previous point, $\delta^\smallfrown \alpha \in W^{i-1}\setminus W^i$ if and only if 
$\delta^\smallfrown \alpha = \delta_1^\smallfrown \Box \chi'^\smallfrown \delta_2$ for some $\langle \delta_1, \Box \chi\rangle \in Gens_i^\Upsilon(\mod{T})$. If $\delta = \delta_1$ we are done. Otherwise, necessarily $\delta = \delta_1^\smallfrown \Box \chi'$, and $\alpha = \delta_2$. However, in this case, again using the previous observation, since $\langle \delta_1, \Box \chi\rangle \in Gens_i^\Upsilon(\mod{T}) $ implies $\delta \in W^{i-1}\setminus W^i$ (it is the root of the tree  $\mod{T}^i_{\delta^\smallfrown\Box \chi}$). This contradicts the fact that $\delta \in W^i$, concluding the proof.

To check the second point, again from the first point we get that 
$\delta = \delta_1^\smallfrown\Box \varphi'^\smallfrown\delta_2$ for some $\langle \delta_1, \Box \varphi \rangle \in Gens_i^\Upsilon(\mod{T})$, which implies that $\vert \delta_1 \vert = i+1$ and that $\delta_- = \delta_1^\smallfrown\Box \varphi^\smallfrown\delta_2$.

If ${\delta_-}^\smallfrown \alpha \in W^{i-1}\setminus W^i$, then it would also hold that $\delta_1^\smallfrown\Box \varphi^\smallfrown\delta_2^\smallfrown \alpha = \varrho_1 ^\smallfrown\Box \psi'^\smallfrown\varrho_2$ for some $\langle \varrho_1, \Box \psi \rangle \in Gens_i^\Upsilon(\mod{T})$, and so, with $\vert \varrho_1\vert = i+1$ too. Therefore, $\delta_1 = \varrho_1$, but since $\Box \varphi \neq \Box \psi'$ this leads to a contradiction, forcing that ${\delta_-}^\smallfrown \alpha \in W^i$. Iterating this reasoning we obtain the second claim.

\end{proof}

\begin{proof}[Proof of Lemma \ref{lem:conservativeM+}]
The second claim follows immediately from the iteration of the first one (up to level $md(\Upsilon)$) and the application of the definition of $\mod{T}^{md(\Upsilon)}$ (since $e^{md(\Upsilon)}$ is $\eta(e)$) and of Lemma \ref{lem:inverseEtaPreservation}.

We will prove the first claim by induction in the complexity of the formula. 
For propositional formulas and connectives the proof follows from the definition of $\mod{T}^{i-1}$ itself. 

For formulas beginning with a $\Diamond$ operator, by I.H. and applying the definition of $R^i$ and $R^{i-1}$ and that of the evaluation of $\Diamond$-formulas we have that 
\begin{align*}
e^{i-1}(\sigma, \Diamond \varphi) = \bigvee_{\sigma^\smallfrown\alpha \in W^{i}} e^{i-1}(\sigma^\smallfrown\alpha, \varphi) \vee \bigvee_{\sigma^\smallfrown\alpha \in W^{i-1}\setminus W^i} e^{i-1}(\sigma^\smallfrown\alpha, \varphi) = \\
\bigvee_{\sigma^\smallfrown\alpha \in W^{i}} e^{i}(\sigma^\smallfrown\alpha, \varphi) \vee \bigvee_{\sigma^\smallfrown\alpha \in W^{i-1}\setminus W^i} e^{i-1}(\sigma^\smallfrown\alpha, \varphi) = \\
e^i(\sigma, \Diamond \varphi) \vee \bigvee_{\sigma^\smallfrown\alpha \in W^{i-1}\setminus W^i} e^{i-1}(\sigma^\smallfrown\alpha, \varphi)
\end{align*}

Consider any $\sigma^\smallfrown\alpha \in W^{i-1}\setminus W^i$. Observation \ref{obs:freshSuccessorsOnlyOriginalWorlds} (2) tells us that this is the case if and only if
$\langle \sigma, \Box \psi\rangle \in Gens_i^\Upsilon(\mod{T})$ for $\alpha = \Box \psi'$. Therefore, $\sigma^\smallfrown\alpha_- = \sigma^\smallfrown \Box \psi$, and hence, by construction, $ \sigma^\smallfrown \Box \psi \in W^i$ and $R^i \sigma, \sigma^\smallfrown \Box \psi$. 

On the other hand, from Lemma  \ref{lem:truthLemmaM+} it follows that that $e^{i-1}(\sigma^\smallfrown\alpha, \varphi) \leq e^i({\sigma^\smallfrown\alpha_-, \varphi})$, and so, $e^i(\sigma, \Diamond \varphi) \geq e^{i-1}(\sigma^\smallfrown\alpha, \varphi)$ concluding the proof of the step.

For formulas beginning with a $\Box$,  similarly by I.H. and the definitions we have that 
\[e^{i-1}(\sigma, \Box \varphi) = \bigwedge_{\sigma^\smallfrown\alpha \in W^{i}} e^{i}(\sigma^\smallfrown\alpha, \varphi) \wedge \bigwedge_{\sigma^\smallfrown\alpha \in W^{i-1}\setminus W^i} e^{i-1}(\sigma^\smallfrown\alpha, \varphi).\]
From Lemma \ref{lem:truthLemmaM+} the previous equals to 
\begin{equation}\label{eq:BoxCalculus}
e^{i}(\sigma, \Box \varphi) \wedge \bigwedge_{\sigma^\smallfrown\alpha \in W^{i-1}\setminus W^i} (e^i(\sigma^\smallfrown\alpha_-, \varphi) + e^i(\sigma ^\smallfrown\alpha_-, \varphi)_{\mapsfrom{{\vartriangleleft}e^i(c(\sigma^\smallfrown\alpha), {pr}(\sigma^\smallfrown\alpha))}}).
\end{equation}

As above, for any world $\sigma^\smallfrown\alpha \in W^{i-1}\setminus W^i$, following Observation \ref{obs:freshSuccessorsOnlyOriginalWorlds} (2), we know that
$\langle \sigma, \Box \psi\rangle \in Gens_i^\Upsilon(\mod{T})$ for $\alpha = \Box \psi'$. Hence $\sigma^\smallfrown\alpha_- = \sigma^\smallfrown\Box \psi$ (so in $W^i$ and with $R^i \sigma, \sigma ^\smallfrown\Box \psi$). Not only this, but also, since $\langle \sigma, \Box \psi\rangle \in Gens_i^\Upsilon(\mod{T})$, necessarily $\sigma \in \eta(RW(\mod{T}))$, and henceforth, 
$(\sigma^\smallfrown\alpha)_- = c(\sigma^\smallfrown\alpha) = \sigma^\smallfrown \Box \psi \in \eta(RW(\mod{T}))$
and ${pr}(\sigma^\smallfrown\alpha) = \psi$. 

All the above, and applying I.H. again (since $\sigma ^\smallfrown\Box \psi \in W^{md(\Upsilon)}$) allow us to conclude that
\begin{align*}
e^{i-1}(\sigma, \Box \varphi) = e^i(\sigma, \Box \varphi) \wedge \bigwedge_{\sigma^\smallfrown\alpha \in W^{i-1}\setminus W^i} (e^i(\sigma^\smallfrown\alpha_-, \varphi) + e^i(\sigma ^\smallfrown\alpha_-, \varphi)_{\mapsfrom{{\vartriangleleft}e^i(c(\sigma^\smallfrown\alpha), {pr}(\sigma^\smallfrown\alpha))}}) = \\
e^i(\sigma, \Box \varphi) \wedge \bigwedge_{\Box \psi \in uWits(\sigma, \Upsilon_{\depth(\sigma)})} (e^i(\sigma ^\smallfrown\Box \psi, \varphi) + e^i(\sigma ^\smallfrown\Box \psi, \varphi)_{\mapsfrom{{\vartriangleleft}e^i(\sigma ^\smallfrown\Box \psi, \psi)}}) = \\
\eta(e)(\sigma, \Box \varphi) \wedge \bigwedge_{\Box \psi \in uWits(\sigma, \Upsilon_{\depth(\sigma)})} (\eta(e)(\sigma ^\smallfrown\Box \psi, \varphi) + \eta(e)(\sigma ^\smallfrown\Box \psi, \varphi)_{\mapsfrom{{\vartriangleleft}\eta(e)(\sigma ^\smallfrown\Box \psi, \psi)}})
\end{align*}
If $\eta(e)(\sigma, \Box \varphi) = \bot$, this concludes the proof. Otherwise, since by Corollary \ref{cor:inverseEtaPreservation}, we know that 

\[{\vartriangleleft}\eta(e)(\sigma ^\smallfrown\Box \psi, \psi) < \bigwedge \{{\vartriangleleft}e(\eta^{-1}(\sigma)_{\Box \psi}, \chi) \colon \Box \chi \in \Upsilon_{i}, e(\eta^{-1}(\sigma), \Box \chi) > \bot\},\] in particular we get that
${\vartriangleleft}\eta(e)(\sigma ^\smallfrown\Box \psi, \psi) < {\vartriangleleft}e(\eta^{-1}(\sigma)_{\Box \psi}, \varphi)$. 

Lemma \ref{lem:inverseEtaPreservation} implies that ${\vartriangleleft}e(\eta^{-1}(\sigma)_{\Box \psi}, \varphi) = {\vartriangleleft}\eta(e)(\sigma ^\smallfrown \Box \psi, \varphi)$. Therefore, Observation \ref{obs:closetriangle} (3) allows us to conclude that $\eta(e)(\sigma ^\smallfrown\Box \psi, \varphi)_{\mapsfrom{{\vartriangleleft}\eta(e)(\sigma ^\smallfrown\Box \psi, \psi)}}) = \mathbf{0}$. 

Therefore, 
\[\bigwedge_{\Box \psi \in uWits(\sigma, \Upsilon_{\depth(\sigma)})} (\eta(e)(\sigma ^\smallfrown\Box \psi, \varphi) + \eta(e)(\sigma ^\smallfrown\Box \psi, \varphi)_{\mapsfrom{{\vartriangleleft}\eta(e)(\sigma ^\smallfrown\Box \psi, \psi)}}) = \bigwedge_{\Box \psi \in uWits(\sigma, \Upsilon_{\depth(\sigma)})} \eta(e)(\sigma ^\smallfrown\Box \psi, \varphi)\]
Since by definition $\eta(e)(\sigma, \Box \varphi) \leq \bigwedge_{\Box \psi \in uWits(\sigma, \Upsilon_{\depth(\sigma)})} \eta(e)(\sigma ^\smallfrown\Box \psi, \varphi)$, this concludes the proof of the step. 

\end{proof}

\begin{proof}[Proof of Lemma \ref{lem:awitnesses}]
First observe that according to the choice of the world $\eta^{-1}(\underline{\sigma})_{\Diamond \varphi}$, from Lemmas \ref{lem:inverseEtaPreservation} and \ref{lem:conservativeM+} it is immediate that

\[e^+(\underline{\sigma}, \Diamond \varphi) = e(\eta^{-1}(\underline{\sigma}), \Diamond \varphi) = e(\eta^{-1}(\underline{\sigma})_{\Diamond \varphi}, \varphi) = e^+(\underline{\sigma}^ \smallfrown\Diamond \varphi, \varphi)\]

From Observation \ref{obs:freshSuccessorsOnlyOriginalWorlds}(3) it immediately follows that $\underline{\sigma}^\smallfrown\Box \psi' \not \in W^+$, so $\langle \underline{\sigma}, \Box \psi\rangle \not \in Gens_{\depth(\sigma)}$. Therefore, similarly, 
\[
e^+(\underline{\sigma}, \Box \psi) = e(\eta^{-1}(\underline{\sigma}), \Box \psi) = e(\eta^{-1}(\underline{\sigma})_{\Box \psi}, \psi) = e^+(\underline{\sigma}^ \smallfrown\Box \psi, \psi)\]

We will use the above two properties, and also resort without further notice to the facts that ${\sigma_-}^\smallfrown\heartsuit \chi = (\sigma^\smallfrown\heartsuit \chi)_-$, that $c(\sigma) = c(\sigma^\smallfrown\heartsuit \chi)$ and that ${pr}(\sigma) = {pr}(\sigma^\smallfrown\heartsuit \chi)$.

Consider then $\heartsuit \chi$ to be either $\Diamond \varphi$ or $\Box \psi$ as in the premises of the Lemma. 
We will prove by induction on the number of primed elements in $\sigma$  that 
\begin{equation}\label{eq:inawitnesses}
e^+(\sigma, \heartsuit \chi) = e^+(\sigma^\smallfrown\heartsuit \chi, \chi).
\end{equation}


To prove the base case of Equation (\ref{eq:inawitnesses}), assume $\sigma_- = \underline{\sigma}$ (namely, if $\sigma$ has only one primed element). 

Applying the definition of $a$ and the above equalities we get that
\[a_{\sigma, \heartsuit \chi}  = e^+(\sigma_-, \heartsuit \chi)_{\mapsfrom {\vartriangleleft}e^+(c(\sigma), {pr}(\sigma))} = e^+({\sigma_-}^\smallfrown \heartsuit \chi, \chi)_{\mapsfrom {\vartriangleleft}e^+(c(\sigma^ \smallfrown \heartsuit \chi), {pr}(\sigma^ \smallfrown \heartsuit \chi)} = a_{\sigma^\smallfrown\heartsuit \chi, \chi}.\]

Now following Proposition \ref{prop:valueProperties} (1) and again the above equalities, we conclude that
\[e^+(\sigma, \heartsuit \chi) = e^+(\sigma_-, \heartsuit \chi) + a_{\sigma, \heartsuit \chi} = e^+({\sigma_-}^ \smallfrown \heartsuit \chi, \chi) + a_{\sigma^\smallfrown\heartsuit \chi, \chi} = e^+(\sigma ^ \smallfrown \heartsuit \chi, \chi).\]

The inductive step of Equation \ref{eq:inawitnesses} is proven analogously: first, by  I.H. and the definition of $a$ we get
\begin{align*}
a_{\sigma, \heartsuit \chi} = e^+(\sigma_-, \heartsuit \chi)_{\mapsfrom {\vartriangleleft}e^+(c(\sigma), {pr}(\sigma))} = 
e^+({\sigma_-}^\smallfrown\heartsuit \chi, \chi)_{\mapsfrom {\vartriangleleft}e^+(c(\sigma^\smallfrown\heartsuit \chi), {pr}(\sigma^\smallfrown\heartsuit \chi))} = 
a_{\sigma^\smallfrown\heartsuit \chi, \chi}.
\end{align*}
The previous, together with Proposition \ref{prop:valueProperties} (1) and I.H. allows us to conclude the proof of Equation \ref{eq:inawitnesses}, since 
\[e^+(\sigma, \heartsuit \chi) = e^+(\sigma_-, \heartsuit \chi) + a_{\sigma, \heartsuit \chi} = e^+({\sigma_-}^ \smallfrown \heartsuit \chi, \chi) + a_{\sigma^\smallfrown\heartsuit \chi, \chi)} = e^+(\sigma^ \smallfrown \heartsuit \chi, \chi).\]

We have already proven that $a_{\sigma, \heartsuit \chi} = a_{\sigma^\smallfrown\heartsuit \chi, \chi}$ at each step. On the other hand, the fact that
$a_{\sigma^\smallfrown\Box \chi, \chi} \leq a_{\sigma^\smallfrown \delta, \chi}$ for any $\sigma ^\smallfrown \delta \in W^+$ (and the dual case for the $\Diamond$ formulas) follows easily from the definition of $a$ and of the evaluation of $\Box \chi$, plus Observation \ref{obs:ordermapsto} (2) and the fact that $c(\sigma) = c(\sigma^\smallfrown \delta)$ and ${pr}(\sigma) = {pr}(\sigma^\smallfrown \delta)$.

\end{proof}

\section*{Appendix 3. Proofs from Section \ref{sec:completeness}.}\label{sec:app3}

\begin{proof}[Proof of Proposition \ref{prop:truthLemmaOP}]
We will prove it by induction on the complexity of the formula. Only the step of the implication will need of the assumption of order preservation.

For the $\Delta$ operation, it is immediate that $fg(\Delta \varphi) = \Delta fg(\varphi) = \Delta (f(\varphi) \cdot g(\varphi))$ by I.H. Since $\Delta$ distributes over $\cdot$ in $[0,1]$, the previous is equal to $\Delta f(\varphi) \cdot \Delta g(\varphi)$, and so, to $f (\Delta \varphi) \cdot g(\Delta \varphi)$. 

The step for the $\with$ operation is also fairly direct by conmutativity of the product. By I.H., 
\[
fg(\varphi \with \chi)  =  fg(\varphi) \cdot  fg(\chi) =
 f(\varphi) \cdot g(\varphi) \cdot f(\chi) \cdot g(\chi).
\]
By commutativity again, the previous is equal to  
$f(\varphi)\cdot f(\psi) \cdot g(\varphi) \cdot g(\psi)$ and so, to $f(\varphi \with \psi) \cdot g(\varphi \with \psi).$

The step for the $\rightarrow$ operation will be proven in two cases.
\begin{enumerate} 
\item Suppose $fg(\varphi) \rightarrow  fg(\chi) < 1$. This holds if and only if $fg(\varphi) > fg(\chi)$, and hence, by I.H. and by the definition of the implication in $[0,1]_{\Pi_\Delta}$, 
\[1 > \frac{f(\chi) \cdot g(\chi)}{f(\varphi)\cdot g(\varphi)} = \frac{f(\chi)}{f(\varphi)}\cdot \frac{g(\chi)}{g(\varphi)}\]
This implies that one of the two fractions is necessarily below $1$. 

Assume wlog that $\frac{f(\chi)}{f(\varphi)} < 1$. Then $f(\chi) < f(\varphi)$, and by definition, $f(\varphi \rightarrow \psi) = \frac{f(\chi)}{f(\varphi)}$. Moreover, by (contraposition of) condition (2) of Definition \ref{def:orderPresDeformation}, $h(\chi) \leq h(\varphi)$, and so, again by the same condition, $g(\chi) \leq g(\varphi)$. Therefore $g(\varphi \rightarrow \psi) = \frac{g(\chi)}{g(\varphi)}$. We conclude that 
\[\frac{f(\chi) \cdot g(\chi)}{f(\varphi)\cdot g(\varphi)} = \frac{f(\chi)}{f(\varphi)}\cdot \frac{g(\chi)}{g(\varphi)} = f(\varphi \rightarrow \chi) \cdot g(\varphi \rightarrow \chi).\]

\item Suppose that, otherwise, $fg(\varphi) \rightarrow  fg(\chi) = 1$, so, by I.H. $f(\varphi)\cdot g(\varphi) \leq f(\chi)\cdot g(\chi)$. 

We have three subcases:
 \begin{itemize}
\item If $f(\varphi) = f(\chi)$, then either a) this value is $0$, and by condition (1) of Definition \ref{def:orderPresDeformation}, also $g(\varphi) = g(\chi) = 0$, or b) by the assumption of the case and monotonicity of the product it is necessary that $g(\varphi) \leq g(\chi)$. In both cases, $f(\varphi \rightarrow \chi) = g(\varphi \rightarrow \chi) = 1$, and their product also equals 1, proving the subcase.
\item If $f(\varphi) < f(\psi)$, by (contraposition of) condition (2) of Definition \ref{def:orderPresDeformation} $g(\varphi) \leq g(\psi)$ and so, $f(\varphi \rightarrow \psi) = g(\varphi \rightarrow \psi) = 1$ proving the subcase as before. 
\item Suppose, with a view to contradiction, that $0 \leq f(\psi) < f(\varphi)$. By condition (2) of Definition \ref{def:orderPresDeformation} (applied twice), also $g(\psi) \leq g(\varphi)$. Thus, to get that $f(\varphi)\cdot g(\varphi) \leq f(\psi)\cdot g(\psi)$, it would be necessary that $f(\varphi)\cdot g(\varphi) = f(\psi)\cdot g(\psi) = 0$. However, in the subcase we have that $0 < f(\varphi)$, and hence, following condition (1) of Definition \ref{def:orderPresDeformation}, $0 \neq h(\varphi)$ and $0 \neq g(\varphi)$ too. Therefore, $f(\varphi)\cdot g(\varphi) > 0$, reaching a contradiction. \qedhere
\end{itemize} 

\end{enumerate}
\end{proof}

\section*{Appendix 4. Some proofs from Section \ref{sec:valued}.}\label{sec:app4}

\begin{proof}[Proof of Lemma \ref{lem:conservativeValuedModel}]
The lemma is immediate for propositional formulas. For what concerns the modal step, by definition 
\[e^{i-1}_v(\sigma, \Box \psi) = \bigwedge_{\delta \in W^{i-1}_v} R^{i-1}_v(\sigma, \delta) \rightarrow e^{i-1}_v(\delta, \psi)\]
By I.H., this equals $\bigwedge \{R^{i-1}_v(\sigma, \delta) \rightarrow e^{i}_v(\delta, \psi)\colon \delta \in W^i\} \wedge \bigwedge \{R^{i-1}_v(\sigma, \delta)  \rightarrow  e^{i-1}_v(\delta, \psi) \colon \delta \in W^{i-1}\setminus W^i\}$, namely $e^i_v(\sigma, \Box \psi) \wedge \bigwedge \{R^{i-1}_v(\sigma, \delta)  \rightarrow  e^{i-1}_v(\delta, \psi) \colon \delta \in W^{i-1}\setminus W^i\}$.

Suppose $e^i_v(\sigma, \Box \psi) > \bot$, since otherwise the case is immediate, and let $\delta \in W^{i-1}\setminus W^i$. From Definition \ref{def:mod+Valued}, either $R^{i-1}_v(\sigma, \delta) = \bot$ (and so, $R^{i-1}_v(\sigma, \delta)  \rightarrow  e^{i-1}_v(\delta, \psi) = \mathbf{0}$, therefore not contributing towards the infimum) or $\delta = \sigma^\smallfrown \Box \chi'$ for $\langle \sigma, \Box \chi \rangle \in (Gens_i^\Upsilon(\mod{T}))_v$. In this latter case, we know by definition that 
$\sigma \in \eta^v(W)$, and so, by I.H., $e^i_v(\sigma, \Box \psi) =  e(\eta^{-1}(\sigma), \Box \psi)$. It is also immediate that $\eta^{-1}(\sigma)_{\Box \chi} = \eta^{-1}(\sigma^\smallfrown \Box \chi) = \eta^{-1}(\delta_-)$. Therefore, since by definition  
$R^{i-1}_v(\sigma, \delta) = r_{\eta^{-1}(\sigma), \Box \chi}$, we know from Lemma \ref{lem:belowWorldValued} that 
\[R^{i-1}_v(\sigma, \delta) < R(\eta^{-1}(\sigma), \eta^{-1}(\sigma)_{\Box \chi}) \wedge \bigwedge \{e(\eta^{-1}(\sigma)_{\Box \chi}, \phi) \colon \Box \phi \in \Upsilon_{\depth(\sigma)}, e(\eta^{-1}(\sigma), \Box \phi)>\bot\}.\] 
In particular, $R^{i-1}_v(\sigma, \delta) < e(\eta^{-1}(\delta_-), \psi) = e^{md(\Upsilon)}_v(\delta_-, \psi)$. 

On the other hand, from the previous lemma, we know that $e^{i-1}_v(\delta, \psi) = e^i_v(\delta_-, \psi)$ for any $\delta \in W^{i-1}\setminus W^i\}$. Since $\delta_- \in W^{md(\Upsilon)}_v$, by I.H., we get that $e^{i-1}_v(\delta, \psi) = e^{md(\Upsilon)}_v(\delta_-, \psi)$. Together with the above inequality we conclude that $R^{i-1}_v(\sigma, \delta) < e^{i-1}_v(\delta, \psi)$. Therefore, $\bigwedge \{R^{i-1}_v(\sigma, \delta)  \rightarrow  e^{i-1}_v(\delta, \psi) \colon \delta \in W^{i-1}\setminus W^i\} = \mathbf{0}$ (i.e., the top element), and thus $e^{i-1}_v(\sigma, \Box \psi) = e^i_v(\sigma, \Box \psi)$, proving the case for formulas beginning by $\Box$.

For the $\Diamond$ case, the proof is similar but simpler. As before, by definition, 
\[e^{i-1}_v(\sigma, \Diamond \psi) = e^i_v(\sigma_-, \Diamond \psi) \vee \bigvee \{R^{i-1}_v(\sigma, \delta)  \with  e^{i-1}_v(\delta, \psi) \colon \delta \in W^{i-1}\setminus W^i\}\]

Again, for $\delta \in W^{i-1}\setminus W^i$, either $R^{i-1}_v(\sigma, \delta) = \bot$ (and so, $R^{i-1}_v(\sigma, \delta)  \with e^{i-1}_v(\delta, \psi) = \bot$ (therefore not contributing towards the supremum) or $\delta = \sigma^\smallfrown \Box \chi'$ for $\langle \sigma, \Box \chi \rangle \in (Gens_i^\Upsilon(\mod{T}))_v$. In this case, in particular, $R^{i-1}_v(\sigma, \delta) < R(\eta^{-1}(\sigma), \eta^{-1}(\sigma)_{\Box \chi}) = R^{md(\Upsilon)}_v(\sigma, \delta_-)$.

Therefore, in combination with Lemma \ref{lem:conservativeValuedModel}, we conclude that 
\begin{align*}
\bigvee \{R^{i-1}_v(\sigma, \delta)  \with  e^{i-1}_v(\delta, \psi) \colon \delta \in W^{i-1}\setminus W^i\} &\leq  \\
 \bigvee \{R^{md(\Upsilon)}_v(\sigma, \delta_-)  \with  e^{md(\Upsilon)}_v(\delta_-, \psi) \colon \delta \in W^{i-1}\setminus W^i\}& \leq  \\ 
e^{i-1}(\sigma, \Diamond \psi).&  \qedhere
\end{align*}
\end{proof}
\end{document}